\documentclass[pdflatex,sn-mathphys,Numbered]{sn-jnl}

\usepackage{graphicx}%
\usepackage{multirow}%
\usepackage{amsmath}%
\usepackage{amsthm}%
\usepackage{mathrsfs}%
\usepackage[title]{appendix}%
\usepackage{xcolor}%
\usepackage{textcomp}%
\usepackage{manyfoot}%
\usepackage{booktabs}%
\usepackage{algorithm}%
\usepackage{algorithmicx}%
\usepackage{algpseudocode}%
\usepackage{listings}%

\makeatletter
\newtheorem*{rep@theorem}{\rep@title}
\newcommand{\newreptheorem}[2]{%
\newenvironment{rep#1}[1]{%
 \def\rep@title{#2 \ref{##1}}%
 \begin{rep@theorem}}%
 {\end{rep@theorem}}}
\makeatother

\theoremstyle{thmstyleone}
\newtheorem{theorem}{Theorem}[section]
\newtheorem{proposition}[theorem]{Proposition}
\newtheorem{corollary}[theorem]{Corollary}
\newtheorem{lemma}[theorem]{Lemma}

\newreptheorem{theorem}{Theorem}
\newreptheorem{corollary}{Corollary}

\theoremstyle{thmstyletwo}%

\theoremstyle{thmstylethree}%

\usepackage{tikz-cd}
\usepackage[all,cmtip,color,matrix,arrow,2cell]{xy}
\UseAllTwocells
\usetikzlibrary{graphs,decorations.pathmorphing,decorations.markings}

\usepackage[charter]{mathdesign}
\DeclareFontFamily{U}{dutchcal}{\skewchar\font=45 }
\DeclareFontShape{U}{dutchcal}{m}{n}{<-> s*[1.0] dutchcal-r}{}
\DeclareFontShape{U}{dutchcal}{b}{n}{<-> s*[1.0] dutchcal-b}{}
\DeclareMathAlphabet{\mathlcal}{U}{dutchcal}{m}{n}
\SetMathAlphabet{\mathlcal}{bold}{U}{dutchcal}{b}{n}

\DeclareSymbolFont{cmcal}{OMS}{cmsy}{m}{n}
\SetSymbolFont{cmcal}{bold}{OMS}{cmsy}{b}{n}
\DeclareSymbolFontAlphabet{\mathcal}{cmcal} 

\DeclareMathAlphabet{\mathpzc}{OT1}{pzc}{m}{it}

\numberwithin{equation}{section}

\def\I{\mathbb{I}}
\def\tensor{\circledast}
\def\unit{\textit{1}}
\def\wtensor{\text{ }\!\widetilde{\circledast}\text{ }\!}

\DeclareMathOperator{\id}{id}
\DeclareMathOperator{\op}{op}

\DeclareMathOperator{\Hom}{Hom}

\usepackage{mathtools}
\newlength{\myeqskip}
\AtBeginDocument{%
  \setlength\abovedisplayskip{\myeqskip}%
  \setlength\belowdisplayskip{\myeqskip}%
  \setlength\abovedisplayshortskip{\myeqskip-\baselineskip}%
  \setlength\belowdisplayshortskip{\myeqskip}
}

\begin{document}

\title[Tensor Enriched Categorical Generalization of the Eilenberg-Watts Theorem]
{Tensor Enriched Categorical Generalization of the Eilenberg-Watts Theorem}


\author*[1]{\fnm{Jaehyeok} \sur{Lee}}\email{jhlee1228lee@ajou.ac.kr}

\affil*[1]{
  \orgdiv{Research Institute of Basic Sciences},
  \orgname{Ajou University},
  \orgaddress{\street{206, Worldcup-ro, Yeongtong-gu}, \city{Suwon-si, Gyeonggi-do}, \postcode{16499},
  \country{Republic of Korea}}}


\abstract{
  Let $b$, $b'$ be commutative monoids in a B\'{e}nabou cosmos.
  Motivated by six-functor formalisms in algebraic geometry,
  we prove that the category of commutative monoids over $b\otimes b'$
  is equivalent to the category of cocontinuous lax monoidal enriched functors
  between the monoidal enriched categories of right modules over $b$, $b'$.
}

\keywords{Eilenberg-Watts theorem, enriched category, symmetric monoidal category, tensorial strength, six-functor formalism}

\pacs[Statements and Declarations]{
  This work was supported by the National Research Foundation of Korea(NRF) grant funded by the Korea government(MSIT)(RS-2024-00334558).
}

\pacs[MSC Classification]{13B10, 13C60, 16D10, 16D20, 16D90, 18D20, 18D25, 18M05}


\pacs[Acknowledgments]{The author is deeply grateful to Jae-Suk Park who always had passionate discussions with the author. He gave helpful comments and introduced relevant references which the author needed to obtain the results of this paper. The author would also like to thank two anonymous reviewers who carefully read through the details\!\!, gave important feedbacks about the content\!\!, and suggested the overall presentation of the paper.}

\pacs[ORCID]{0000-0003-3357-0088}

\maketitle

\section{Introduction} \label{sec Intro}
A six-functor formalism is an assignment
$\mathcal{X}\mapsto \mathcal{D}(\mathcal{X})$
which sends each geometric object $\mathcal{X}$
to a category $\mathcal{D}(\mathcal{X})$,
such that for each morphism
$f:\mathcal{X}'\to\mathcal{X}$
of geometric objects
we have six functors
$f^*$, $f_*$,
$\otimes$, $\underline{\textit{Hom}}$,
$f_!$, $f^!$
between
$\mathcal{D}(\mathcal{X})$, $\mathcal{D}(\mathcal{X}')$
satisfying certain properties.
The category $\mathcal{D}(\mathcal{X})$
is often a triangulated category or a stable $\infty$-category.
Such formalism was first formulated during the development of \'{e}tale cohomology groups of schemes,
as it contains the structural properties of cohomology groups.
Six-functor formalisms are realized, or expected to exist
in various contexts throughout geometry and topology.
We refer to \cite{Scholze2022}
for a detailed explanation, examples and recent developments
of six-functor formalisms.

Consider the case of affine schemes $\mathcal{X}$
and the categories $\textit{QCoh}(\mathcal{X})$ of their quasi-coherent sheaves.
For each morphism
$f:\mathcal{X}'\to \mathcal{X}$
of affine schemes, we have a pair of adjunctions
$f^*\dashv f_*$
,
$f_!\dashv f^!$
between $\textit{QCoh}(\mathcal{X})$, $\textit{QCoh}(\mathcal{X}')$
and a natural isomorphism
$f_!\cong f_*$.
Furthermore,
$f^*:\textit{QCoh}(\mathcal{X})\to \textit{QCoh}(\mathcal{X}')$
is a symmetric monoidal functor
and the lax symmetric monoidal adjunction
$f^*\dashv f_*$ satisfies the projection formula.
For each pair of morphisms
$\mathcal{X}\leftarrow \tilde{\mathcal{X}}\to \mathcal{X}'$
of affine schemes,
which is a correspondence from $\mathcal{X}$ to $\mathcal{X}'$,
we associate a cocontinuous lax symmetric monoidal functor
$\textit{QCoh}(\mathcal{X})\to \textit{QCoh}(\mathcal{X}')$
as follows.
\begin{equation} \label{eq Intro motivation}
  \vcenter{\hbox{
    \xymatrix@R=15pt@C=15pt{
      \text{ }
      &\tilde{\mathcal{X}}
      \ar[dl]_-{f}
      \ar[dr]^-{f'}
      &\text{ }
      \\
      \mathcal{X}
      &\text{ }
      &\mathcal{X}'
    }
  }}
  \qquad\mapsto\qquad
  \vcenter{\hbox{
    \xymatrix{
      \textit{QCoh}(\mathcal{X})
      \ar[r]^-{f'_*f^*}
      &\textit{QCoh}(\mathcal{X}')
    }
  }}
\end{equation}
The six functors
associated to each morphism of affine schemes
are additive functors
(i.e., enriched in the closed symmetric monoidal category $\textit{Ab}$ of abelian groups),
and their symmetric monoidal structures are compatible with additive structures.
It turns out that if we impose this compatibility condition,
(\ref{eq Intro motivation})
becomes an equivalence of categories.
This is an immediate consequence of the result of this paper.

The main result of this paper states that the assignment
(\ref{eq Intro motivation})
defines an equivalence of categories,
not only for affine schemes,
but for commutative monoids in an abstract symmetric monoidal category.
A \emph{B\'{e}nabou cosmos}
$(\mathcal{V},\otimes,c)$
is a closed symmetric monoidal category
whose underlying category is complete and cocomplete.
Such monoidal category $\mathcal{V}$ has enough properties
so that a large amount of the theory of ordinary categories extends 
to the theory of categories enriched in $\mathcal{V}$:
see \cite[p. 134]{Street2006}.
We present the main result.

\begin{reptheorem}{thm Intro MainThm}
  Let $b$, $b'$ be commutative monoids
  in a B\'{e}nabou cosmos $\mathcal{V}$
  and let
  $\textit{Mod}_b$, $\textit{Mod}_{b'}$
  be the symmetric monoidal $\mathcal{V}$-enriched categories of right modules over $b$, $b'$.
  \begin{enumerate}
    \item 
    We have a fully faithful left adjoint functor
    \begin{equation*}
      \vcenter{\hbox{
        \xymatrix@R=0pt{
          \textit{Comm}_{b\otimes b'}
          \ar[r]
          &\mathcal{V}\text{-}\textit{Funct}^{\textit{lax}}(\textit{Mod}_b,\textit{Mod}_{b'}\!)
          \\
          b
          \xrightarrow[]{f}
          \tilde{b}
          \xleftarrow[]{f'}
          b'
          \ar@{|->}[r]
          &f'_*f^*:
          \textit{Mod}_b\to \textit{Mod}_{b'}
        }
      }}
    \end{equation*}
    from the category of commutative monoids
    in $\mathcal{V}$ over $b\otimes b'$
    to the category of lax symmetric monoidal $\mathcal{V}$-enriched functors
    $\textit{Mod}_b\to \textit{Mod}_{b'}$.

    \item
    The essential image of this functor
    is the coreflective full subcategory
    $\mathcal{V}\text{-}\textit{Funct}^{\textit{lax}}_{\textit{cocon}}(\textit{Mod}_b,\textit{Mod}_{b'}\!)$
    of cocontinuous lax symmetric monoidal $\mathcal{V}$-enriched functors
    $\textit{Mod}_b\to \textit{Mod}_{b'}$,
    and we have an adjoint equivalence of categories
    \begin{equation*}
      \vcenter{\hbox{
        \xymatrix@R=0pt{
          \textit{Comm}_{b\otimes b'}
          \ar@<0.5ex>[r]^-{\simeq}
          &\text{ }
          \mathcal{V}\text{-}\textit{Funct}^{\textit{lax}}_{\textit{cocon}}(\textit{Mod}_b,\textit{Mod}_{b'}\!)
          .
          \ar@<0.5ex>[l]^-{\simeq}
        }
      }}
    \end{equation*}
  \end{enumerate}
\end{reptheorem}

The adjoint equivalence of categories
in Theorem~\ref{thm Intro MainThm}
for the case $\mathcal{V}=\textit{Ab}$ tells us that
affine schemes and their morphisms can be understood
in terms of their associated six functors.
As this point of view towards affine schemes is immediately generalized to 
commutative monoids in an arbitrary B\'{e}nabou cosmos $\mathcal{V}$,
the author expects that certain topics in algebraic geometry
can be developed over any $\mathcal{V}$.

Brandenburg and Chirvasitu
\cite[Theorem 3.4.3]{Brandenburg2014}
showed that given a quasi-compact quasi-separated scheme
$\mathcal{X}$ and an arbitrary scheme $\mathcal{X}'$,
the assignment $f\mapsto f^*$
defines an equivalence 
from the discrete category
of morphisms $\mathcal{X}'\to \mathcal{X}$
of schemes
to the category of cocontinuous strong symmetric monoidal functors
$\textit{QCoh}(\mathcal{X})\to \textit{QCoh}(\mathcal{X}')$.
The following, which is a consequence of Theorem~\ref{thm Intro MainThm},
generalizes the particular case of
\cite[Theorem 3.4.3]{Brandenburg2014}
when schemes $\mathcal{X}$, $\mathcal{X}'$ are affine.

\begin{reptheorem}{thm Intro MainThm2}
  Let $b$, $b'$ be commutative monoids
  in a B\'{e}nabou cosmos $\mathcal{V}$
  and let
  $\textit{Mod}_b$, $\textit{Mod}_{b'}$
  be the symmetric monoidal $\mathcal{V}$-enriched categories of right modules over $b$, $b'$.
  We have an equivalence of categories
  \begin{equation*}
    \xymatrix@R=0pt{
      \Hom_{\textit{Comm}(\mathcal{V})}(b,b')
      \ar[r]^-{\simeq}
      &\mathcal{V}\text{-}\textit{Funct}^{\textit{strong}}_{\textit{cocon}}(\textit{Mod}_b,\textit{Mod}_{b'}\!)
      \\
      f:b\to b'
      \text{ }
      \ar@{|->}[r]^-{\simeq}
      &\text{ }
      f^*:\textit{Mod}_b\to \textit{Mod}_{b'}
    }
  \end{equation*}
  from the discrete category of morphisms
  $b\to b'$ of commutative monoids in $\mathcal{V}$
  to the category of cocontinuous strong symmetric monoidal $\mathcal{V}$-enriched functors
  $\textit{Mod}_b\to \textit{Mod}_{b'}$.
\end{reptheorem}

The well-known theorem of 
Eilenberg \cite{Eilenberg1960} and Watts \cite{Watts1960}
states that every cocontinuous functor
$\textit{Mod}_R\to\textit{Mod}_{R'}$
between the categories of right modules over rings $R$, $R'$
is naturally isomorphic to
$-\otimes_RN:\textit{Mod}_R\to\textit{Mod}_{R'}$
for some $(R,R')$-bimodule $N$.
We introduce the following generalization of the Eilenberg-Watts theorem,
which serves as a midpoint towards Theorem~\ref{thm Intro MainThm}.

\begin{reptheorem}{thm Intro EnEW}
  Let $b$ be a monoid in a B\'{e}nabou cosmos $\mathcal{V}$
  and let $\mathcal{C}$ be a cocomplete $\mathcal{V}$-enriched category.
  \begin{enumerate}
    \item 
    We have a fully faithful left adjoint functor
    \begin{equation*}
      \vcenter{\hbox{
        \xymatrix@R=0pt{
          \text{ }\!\!_{b}\mathcal{C}
          \ar[r]
          &\mathcal{V}\textit{-Funct}(\textit{Mod}_b,\mathcal{C})
          \\
          X
          \ar@{|->}[r]
          &-\tensor X:\textit{Mod}_b\to \mathcal{C}
        }
      }}
    \end{equation*}
    from the category of left $b$-module objects in $\mathcal{C}$
    to the category of $\mathcal{V}$-enriched functors
    $\textit{Mod}_b\to\mathcal{C}$.

    \item
    The essential image of this functor
    is the coreflective full subcategory
    $\mathcal{V}\textit{-Funct}_{\textit{cocon}}(\textit{Mod}_b,\mathcal{C})$
    of cocontinuous $\mathcal{V}$-enriched functors
    $\textit{Mod}_b\to \mathcal{C}$,
    and we have an adjoint equivalence of categories
    \begin{equation*}
      \vcenter{\hbox{
        \xymatrix@R=0pt{
          \text{ }\!\!_{b}\mathcal{C}
          \ar@<0.5ex>[r]^-{\simeq}
          &\mathcal{V}\textit{-Funct}_{\textit{cocon}}(\textit{Mod}_b,\mathcal{C}).
          \ar@<0.5ex>[l]^-{\simeq}
        }
      }}
    \end{equation*}
  \end{enumerate}
\end{reptheorem}

Although there are some minor improvements,
the statement and the proof of Theorem~\ref{thm Intro EnEW}
are mostly due to the work of Hovey in
\cite[\textsection 1-\textsection 2]{Hovey2009}.
We are going to follow Hovey's approach and
present a detailed proof of Theorem~\ref{thm Intro EnEW},
which was a necessary step for the author to obtain Theorem~\ref{thm Intro MainThm}.

We list some generalizations of the Eilenberg-Watts theorem
presented in the literature,
and explain how they are related to 
the results in this paper.
\begin{itemize}
  \item 
  In the theory of categories enriched in a B\'{e}nabou cosmos,
  the Eilenberg-Watts theorem holds in great generality
  where rings, modules and bimodules can be replaced by 
  small enriched categories, enriched presheaves and enriched profunctors.
  \footnote{
    It is stated as Theorem 2.5 in the nLab page,
    https://ncatlab.org/nlab/show/Eilenberg-Watts+theorem
  }
  A variant of this is the equivalence in Theorem~\ref{thm Intro EnEW}.
  
  \item
  In \cite[Proposition 1.1]{Hovey2009}
  Hovey shows that for each $\mathcal{V}$-enriched functor
  $\mathcal{F}:\textit{Mod}_b\to \textit{Mod}_{b'}$,
  we have a $(b,b')$-bimodule $\mathcal{F}(b)$
  and a natural transformation
  $\lambda^{\mathcal{F}}:-\tensor \mathcal{F}(b)\Rightarrow\mathcal{F}$.
  In \cite[Theorem 2.3]{Hovey2009}
  he also gives a characterization of $\mathcal{F}$
  whose associated $\lambda^{\mathcal{F}}$ is a natural isomorphism.
  These results,
  together with an additional property that
  $\lambda^{\mathcal{F}}$ is enriched in $\mathcal{V}$,
  are presented with details in this paper as
  Proposition~\ref{prop EnEWThm defining lambda}.
  
  \item 
  In \cite{Nyman2010}
  Nyman studies obstructions to
  an $R$-linear cocontinuous functor
  $\mathcal{F}:\textit{QCoh}(\mathcal{X})\to \textit{QCoh}(\mathcal{X}')$
  between the categories of quasi-coherent sheaves 
  on $R$-schemes $\mathcal{X}$, $\mathcal{X}'$
  being isomorphic to tensoring with a bimodule.
  Main difference with our results 
  is that we are interested in characterizing
  which lax symmetric monoidal
  functor $\mathcal{F}$
  is isomorphic to $f'_*f^*$
  for some correspondence from
  $\mathcal{X}$ to $\mathcal{X}'$.
  Another difference is that we consider $\mathcal{X}$, $\mathcal{X}'$
  to be affine schemes,
  but defined internal to a B\'{e}nabou cosmos $\mathcal{V}$.

  \item
  In \cite[Theorem 3.4]{Femic2016}
  Femi\'{c} proves the Eilenberg-Watts theorem for
  $2$-categories of the representation categories
  $\mathscr{C}\textit{-Mod}$
  of finite tensor categories $\mathscr{C}$.
  The author does not see any precise connection
  between \cite[Theorem 3.4]{Femic2016}
  and the results presented in this paper.
\end{itemize}

\textbf{Summary.}
We summarize the content as follows.
The author would like to mention that 
the paper is partially expository
and contains details of proofs.

\begin{itemize}
  \item 
  In \textsection\!~\ref{sec Preliminaries}
  we fix a B\'{e}nabou cosmos $(\mathcal{V},\otimes,c)$
  and explain notations, terminologies
  as well as elementary categorical tools
  used in the paper.
  Every enriched structure is considered as enriched in $\mathcal{V}$.

  \item
  In \textsection\!~\ref{sec EnEWthm}
  we present a detailed proof of Theorem~\ref{thm Intro EnEW},
  following the approach introduced in \cite[\textsection 1-\textsection 2]{Hovey2009}.
  Let $b$ be a monoid in $\mathcal{V}$
  and let $\mathcal{C}$ be a cocomplete enriched category.
  For each enriched functor
  $\mathcal{F}:\textit{Mod}_b\to \mathcal{C}$,
  the image $\mathcal{F}(b)$
  of the right $b$-module $b$ has a canonical left $b$-action.
  The main step is to show that
  we have an enriched natural transformation
  \begin{equation} \label{eq Intro lambdaF}
    \lambda^{\mathcal{F}}:
    -\tensor \mathcal{F}(b)\Rightarrow\mathcal{F}
  \end{equation}
  between enriched functors $\textit{Mod}_b\to\mathcal{C}$,
  which is an enriched natural isomorphism
  if and only if
  $\mathcal{F}$ is cocontinuous as an enriched functor.
  After proving Theorem~\ref{thm Intro EnEW},
  we explain the result when
  $\mathcal{C}=\textit{Mod}_{b'}$
  for another monoid $b'$ in $\mathcal{V}$.

  \item
  In \textsection\!~\ref{sec SixFunctors}
  we establish the six functors
  associated to each morphism of commutative monoids
  $f:b\to b'$ in $\mathcal{V}$.
  We describe the closed symmetric monoidal
  enriched category 
  $(\textit{Mod}_b,\tensor,b)$,
  the pair of enriched adjunctions
  $f^*\dashv f_*$, $f_!\dashv f^!$,
  the enriched natural isomorphism
  $f_!\cong f_*$
  and the lax symmetric monoidal enriched adjunction
  $f^*\dashv f_*$ satisfying the projection formula.

  \item
  In \textsection\!~\ref{sec MainResult}
  we prove Theorem~\ref{thm Intro MainThm}
  as well as Theorem~\ref{thm Intro MainThm2}.
  Let $b$, $b'$ be commutative monoids in $\mathcal{V}$.
  Each lax symmetric monoidal enriched functor
  $\mathcal{F}:\textit{Mod}_b\to\textit{Mod}_{b'}$
  determines a pair of morphisms
  $f:b\to \mathcal{F}(b)$,
  $f':b'\to\mathcal{F}(b)$
  of commutative monoids in $\mathcal{V}$.
  The main step is to show that
  the enriched natural transformation (\ref{eq Intro lambdaF})
  becomes a monoidal enriched natural transformation
  $f'_*f^*\Rightarrow\mathcal{F}$.

  \item
  In \textsection\!~\ref{sec DayConvolution}
  we provide another method of 
  obtaining the equivalence of categories
  in Theorem~\ref{thm Intro MainThm}.
  The main observation is that
  cocontinuous lax symmetric monoidal enriched functors
  $\textit{Mod}_b\to\textit{Mod}_{b'}$
  are identified with commutative monoids in the category
  of cocontinuous enriched functors
  $\textit{Mod}_b\to \textit{Mod}_{b'}$
  with respect to Day convolution.

  \item 
  In the appendix \textsection\!~\ref{sec Appendix},
  we provide details which were omitted in proofs.
\end{itemize}

The author would like to end this section
with an explanation of the title of the paper.
By `Tensor enriched categorical'
it means two categorical structures.
One is symmetric monoidal enriched categorical structures.
The other is tensored enriched categorical structures,
such as tensorial strengths associated to
enriched functors between tensored enriched categories.

\section{Preliminaries}
\label{sec Preliminaries}
A \emph{symmetric monoidal category}
is a tuple $(\mathcal{C},\tensor,\unit;a,\imath,\jmath,s)$
of a category $\mathcal{C}$,
a functor $\tensor:\mathcal{C}\times \mathcal{C}\to \mathcal{C}$,
a distinguished object $\unit$ in $\mathcal{C}$ called the \emph{unit object},
and isomorphisms 
\begin{equation*}
  Z\tensor (X\tensor Y)\xrightarrow[\cong]{a_{Z,X,Y}}(Z\tensor X)\tensor Y,
  \qquad
  \begin{aligned}
    &\unit\tensor X\xrightarrow[\cong]{\imath_X}X,
    \\
    &X\tensor \unit\xrightarrow[\cong]{\jmath_X}X,
  \end{aligned}
  \qquad
  X\tensor Y\xrightarrow[\cong]{s_{X,Y}}Y\tensor X
  ,
  \qquad
  Z,X,Y\in\mathcal{C}
\end{equation*}
in $\mathcal{C}$ called the 
\emph{coherence isomorphisms}
which are natural in variables $X$, $Y$, $Z\in \mathcal{C}$
and satisfy the pentagon, unity, symmetry and hexagon axioms.
We denote a symmetric monoidal category
$(\mathcal{C},\tensor,\unit;a,\imath,\jmath,s)$
as $\mathcal{C}$ when the rest of the structures are clear in the given context.
We say a symmetric monoidal category
$\mathcal{C}$ is \emph{closed}
if for each object $X$ in $\mathcal{C}$,
the functor $-\tensor X:\mathcal{C}\to \mathcal{C}$
admits a right adjoint $[X,-]:\mathcal{C}\to \mathcal{C}$.

Let $\mathcal{C}$ be a symmetric monoidal category.
A diagram in $\mathcal{C}$ is called \emph{formal}
when it only consists of the coherence isomorphisms
$a$, $\imath$, $\jmath$, $s$,
their inverses
and monoidal products of these isomorphisms with identity morphisms.
The coherence theorem, proved by S. Mac Lane in \cite{Mac1963},
states that every formal diagram in $\mathcal{C}$ commutes.
We will use this fact freely.

Let $\mathcal{C}$, $\mathcal{D}$ be symmetric monoidal categories.
A \emph{lax (symmetric) monoidal functor}
$(\mathcal{F},\mathcal{F}^0,\mathcal{F}^2):\mathcal{C}\to \mathcal{D}$
is a tuple of a functor 
$\mathcal{F}:\mathcal{C}\to \mathcal{D}$,
a morphism
$\mathcal{F}^0:\unit\to\mathcal{F}(\unit)$
in $\mathcal{D}$
and a collection of morphisms
\begin{equation*}
  \mathcal{F}^2_{X,Y}:
  \mathcal{F}(X)\tensor \mathcal{F}(Y)
  \to \mathcal{F}(X\tensor Y)
  ,
  \qquad
  X,Y\in \mathcal{C}
\end{equation*}
in $\mathcal{D}$
which is natural in variables $X$, $Y$,
such that the morphisms $\mathcal{F}^0$, $\mathcal{F}^2$
are compatible with coherence isomorphisms
$a$, $\imath$, $\jmath$, $s$ of $\mathcal{C}$, $\mathcal{D}$.
We denote a lax monoidal functor
$(\mathcal{F},\mathcal{F}^0,\mathcal{F}^2)$
as $\mathcal{F}$,
and the morphism $\mathcal{F}^2_{X,Y}$ as $\mathcal{F}^2$
when the meaning is clear in context.
A \emph{colax (symmetric) monoidal functor} 
$(\mathcal{F},\mathring{\mathcal{F}}^0,\mathring{\mathcal{F}}^2):\mathcal{C}\to \mathcal{D}$
is similarly defined, but the morphisms
$\mathring{\mathcal{F}}^0:\mathcal{F}(\unit)\to \unit$,
$\mathring{\mathcal{F}}^2_{X,Y}:\mathcal{F}(X\tensor Y)\to \mathcal{F}(X)\tensor \mathcal{F}(Y)$
in $\mathcal{D}$ are in the opposite direction.
Again, we denote a colax monoidal functor
$(\mathcal{F},\mathring{\mathcal{F}}^0,\mathring{\mathcal{F}}^2)$
as $\mathcal{F}$,
and the morphism $\mathring{\mathcal{F}}^2_{X,Y}$ as $\mathring{\mathcal{F}}^2$
when the meaning is clear in context.
A \emph{strong (symmetric) monoidal functor}
from $\mathcal{C}$ to $\mathcal{D}$
is a lax monoidal functor
$\mathcal{F}:\mathcal{C}\to \mathcal{D}$
whose morphisms
$\mathcal{F}^0$, $\mathcal{F}^2$ in $\mathcal{D}$
are isomorphisms.

Let $\mathcal{F}$, $\mathcal{G}:\mathcal{C}\to \mathcal{D}$
be lax monoidal functors between symmetric monoidal categories $\mathcal{C}$, $\mathcal{D}$.
A \emph{monoidal natural transformation}
from $\mathcal{F}$ to $\mathcal{G}$
is a natural transformation $\mathcal{F}\Rightarrow\mathcal{G}$
which is compatible with $\mathcal{F}^0$, $\mathcal{G}^0$ and $\mathcal{F}^2$, $\mathcal{G}^2$.
For detailed explanations of symmetric monoidal categories,
lax/colax/strong monoidal functors and monoidal natural transformations,
see \cite[\textsection 1.2]{Johnson2021}.

We refer
\cite{Kelly2005}, \cite[\textsection 6]{Borceux1994}
for the basics of enriched category theory.
A \emph{B\'{e}nabou cosmos}
is a closed symmetric monoidal category
$(\mathcal{V},\otimes,c;a,\imath,\jmath,s)$
whose underlying category
$\mathcal{V}$ is complete and cocomplete.
Some examples of B\'{e}nabou cosmoi are
the closed symmetric monoidal categories
\begin{equation*}
  \textit{Set}
  ,\text{ }
  \textit{Cat}
  ,\text{ }
  \textit{sSet}
  ,\text{ }
  \textit{Mod}_R
  ,\text{ }
  \textit{Ch}_R
  ,\text{ }
  \textit{CGWH}
  ,\text{ }
  \textit{Sp}^{\Sigma}
  ,\text{ }
  \textit{Ban}
\end{equation*}
of
(small) sets,
categories,
simplicial sets,
modules over a commutative ring $R$,
(co)chain complexes over a commutative ring $R$,
compactly generated weakly Hausdorff spaces,
symmetric spectra,
Banach spaces and linear contractions.
Every Grothendieck topos is also a B\'{e}nabou cosmos.

Throughout this paper,
we fix a B\'{e}nabou cosmos $\mathcal{V}$.
Objects in the fixed B\'{e}nabou cosmos $\mathcal{V}$
are denoted with small letters $w$, $x$, $y$, $z\in \mathcal{V}$.
For instance, the unit object is denoted as $c\in\mathcal{V}$.
Every enriched structure in this paper
is considered as enriched in $\mathcal{V}$.
We use the same notations for
enriched categorical structures and their underlying (ordinary) category structures.

Let $\mathcal{C}$ be an enriched category
and let $X$, $Y$, $Z$ be objects in $\mathcal{C}$.
By a morphism in $\mathcal{C}$
we mean a morphism $X\to Y$ in the underlying category of $\mathcal{C}$.
We denote the Hom object from $X$ to $Y$ as
$\mathcal{C}(X,Y)\in \mathcal{V}$,
and the identity, composition morphisms in $\mathcal{V}$ as
$\I_X:X\to X$, $\mu_{X,Y,Z}:\mathcal{C}(Y,Z)\otimes \mathcal{C}(X,Y)\to \mathcal{C}(X,Z)$.
Let $\mathcal{D}$ be another enriched category
and let $\mathcal{F}:\mathcal{C}\to \mathcal{D}$
be an enriched functor.
We denote the induced morphism between Hom objects as
\begin{equation*}
  \mathcal{F}_{X,Y}:\mathcal{C}(X,Y)\to \mathcal{D}(\mathcal{F}(X),\mathcal{F}(Y)).
\end{equation*}
Let $\mathcal{G}:\mathcal{C}\to \mathcal{D}$
be another enriched functor and let
$\xi:\mathcal{F}\Rightarrow\mathcal{G}$
be an enriched natural transformation.
The component of $\xi$ at each object $X$ in $\mathcal{C}$
is denoted as $\xi_X:\mathcal{F}(X)\to \mathcal{G}(X)$.

A \emph{(symmetric) monoidal enriched category}
is a tuple $(\mathcal{C},\tensor,\unit;a,\imath,\jmath,s)$
of an enriched category $\mathcal{C}$,
an enriched functor $\tensor:\mathcal{C}\times \mathcal{C}\to \mathcal{C}$,
a distinguished object $\unit$ in $\mathcal{C}$
and enriched natural isomorphisms
$a$, $\imath$, $\jmath$, $s$
satisfying the pentagon, unity, symmetry and hexagon axioms.
Again, we denote a monoidal enriched category
$(\mathcal{C},\tensor,\unit;a,\imath,\jmath,s)$
as $\mathcal{C}$ when there is no confusion in meaning about the rest of the structures.
A monoidal enriched category $\mathcal{C}$ is called \emph{closed}
if for each object $X$ in $\mathcal{C}$
the enriched functor $-\tensor X:\mathcal{C}\to \mathcal{C}$
has a right enriched adjoint $[X,-]:\mathcal{C}\to \mathcal{C}$.

We can enhance $\mathcal{V}$ to a closed monoidal enriched category
which is complete and cocomplete as an enriched category.
Let $x$, $y$ be objects in $\mathcal{V}$.
The Hom object is given by $\mathcal{V}(x,y)=[x,y]$,
and the enriched functor $-\otimes x:\mathcal{V}\to \mathcal{V}$
is left adjoint to the Hom enriched functor
$\mathcal{V}(x,-):\mathcal{V}\to \mathcal{V}$.

Let $\mathcal{C}$, $\mathcal{D}$
be monoidal enriched categories.
A \emph{lax (symmetric) monoidal enriched functor}
$(\mathcal{F},\mathcal{F}^0,\mathcal{F}^2):\mathcal{C}\to \mathcal{D}$
is a tuple of an enriched functor 
$\mathcal{F}:\mathcal{C}\to \mathcal{D}$,
a morphism
$\mathcal{F}^0:\unit\to\mathcal{F}(\unit)$
in $\mathcal{D}$
and a collection of morphisms
$\mathcal{F}^2_{X,Y}:
\mathcal{F}(X)\tensor \mathcal{F}(Y)
\to \mathcal{F}(X\tensor Y)$
in $\mathcal{D}$
which is enriched natural in variables $X$, $Y\in \mathcal{C}$,
such that $\mathcal{F}^0$, $\mathcal{F}^2$
are compatible with coherence isomorphisms of $\mathcal{C}$, $\mathcal{D}$.
We also denote a lax monoidal enriched functor
$(\mathcal{F},\mathcal{F}^0,\mathcal{F}^2)$
as $\mathcal{F}$,
and the morphism $\mathcal{F}^2_{X,Y}$ as $\mathcal{F}^2$.
We omit the explanation of \emph{colax/strong (symmetric) monoidal enriched functors}.

Let $\mathcal{F}$, $\mathcal{G}:\mathcal{C}\to \mathcal{D}$
be lax monoidal enriched functors between
monoidal enriched categories $\mathcal{C}$, $\mathcal{D}$.
A \emph{monoidal enriched natural transformation}
from $\mathcal{F}$ to $\mathcal{G}$
is an enriched natural transformation $\mathcal{F}\Rightarrow\mathcal{G}$
which is compatible with $\mathcal{F}^0$, $\mathcal{G}^0$ and $\mathcal{F}^2$, $\mathcal{G}^2$.

A \emph{lax (symmetric) monoidal enriched adjunction} 
is an adjunction internal to the $2$-category of
monoidal enriched categories,
lax monoidal enriched functors
and monoidal enriched natural transformations.
A pair of lax monoidal enriched functors
$\mathcal{F}:\mathcal{C}\to \mathcal{D}$,
$\mathcal{G}:\mathcal{D}\to \mathcal{C}$
between monoidal enriched categories $\mathcal{C}$, $\mathcal{D}$
forms a lax monoidal enriched adjunction
$\mathcal{F}\dashv \mathcal{G}$
if there exists unit, counit monoidal enriched natural transformations
$\eta:\id_{\mathcal{C}}\Rightarrow \mathcal{G}\mathcal{F}$,
$\varepsilon:\mathcal{F}\mathcal{G}\Rightarrow \id_{\mathcal{D}}$
which satisfy the triangle identities.
We say a lax monoidal enriched adjunction
$\mathcal{F}\dashv \mathcal{G}$
\emph{satisfies the projection formula} if
the following composition is an isomorphism
for all $X\in\mathcal{C}$, $A\in\mathcal{D}$.
\begin{equation}\label{eq Preliminary projformuladef}
  \xymatrix@C=35pt{
    X\tensor \mathcal{G}(A)
    \ar[r]^-{\eta_X\tensor I_{\mathcal{G}(A)}}
    &\mathcal{G}\mathcal{F}(X)\tensor \mathcal{G}(A)
    \ar[r]^-{\mathcal{G}^2_{\mathcal{F}(X),A}}
    &\mathcal{G}(\mathcal{F}(X)\tensor A)
  }
\end{equation}

The next lemma is an instance of a result
by Kelly
\cite[Proposition 1.3]{Kelly1974}
on doctrinal adjunctions, obtained by considering the 2-monad on the 2-category of enriched categories whose algebras are
monoidal enriched categories.

\begin{lemma} \label{lem Intro laxtensorCadj}
  Let $\mathcal{C}$, $\mathcal{D}$ be monoidal enriched categories
  and let
  \begin{equation*}
    \vcenter{\hbox{
      \xymatrix{
        \mathcal{C}
        \ar@<1ex>[r]^-{\mathcal{F}}
        \ar@{}[r]|-{\bot}
        &\mathcal{D}
        \ar@<1ex>[l]^-{\mathcal{G}}
      }
    }}
  \end{equation*}
  be an enriched adjunction
  with unit, counit enriched natural transformations
  $\eta$, $\varepsilon$.
  Suppose the right enriched adjoint
  $\mathcal{G}:\mathcal{D}\to\mathcal{C}$
  has a structure of a lax monoidal enriched functor.
  \begin{enumerate}
    \item 
    The left enriched adjoint
    $\mathcal{F}:\mathcal{C}\to \mathcal{D}$
    becomes a colax monoidal enriched functor,
    whose morphisms $\mathring{\mathcal{F}}^0$, $\mathring{\mathcal{F}}^2$
    are given as follows. Let $X$, $Y\in \mathcal{C}$.
    \begin{equation*}
      \vcenter{\hbox{
        \xymatrix@C=30pt{
          \mathcal{F}(\textit{1})
          \ar[d]_-{\mathcal{F}(\mathcal{G}^0)}
          \ar@/^1pc/[dr]^-{\mathring{\mathcal{F}}^0}
          \\
          \mathcal{F}\mathcal{G}(\textit{1})
          \ar[r]^-{\varepsilon_{\textit{1}}}
          &\textit{1}
        }
      }}
      \qquad\quad
      \vcenter{\hbox{
        \xymatrix@R=20pt@C=40pt{
          \mathcal{F}(X\tensor Y)
          \ar[d]_-{\mathcal{F}(\eta_X\tensor \eta_Y)}
          \ar@/^2pc/[ddr]^-{\mathring{\mathcal{F}}^2_{X,Y}}
          \\
          \mathcal{F}(\mathcal{G}\mathcal{F}(X)\tensor \mathcal{G}\mathcal{F}(Y))
          \ar[d]_-{\mathcal{F}(\mathcal{G}^2_{\mathcal{F}(X),\mathcal{F}(Y)})}
          \\
          \mathcal{F}\mathcal{G}(\mathcal{F}(X)\tensor \mathcal{F}(Y))
          \ar[r]^-{\varepsilon_{\mathcal{F}(X)\tensor\mathcal{F}(Y)}}
          &\mathcal{F}(X)\tensor \mathcal{F}(Y)
        }
      }}
    \end{equation*}

    \item
    If the morphisms
    $\mathring{\mathcal{F}}^0$, $\mathring{\mathcal{F}}^2$
    are isomorphisms,
    then $\mathcal{F}:\mathcal{C}\to \mathcal{D}$
    is a strong monoidal enriched functor with
    $\mathcal{F}^0=(\mathring{\mathcal{F}}^0)^{-1}$,
    $\mathcal{F}^2=(\mathring{\mathcal{F}}^2)^{-1}$
    and $\mathcal{F}\dashv \mathcal{G}$ is a lax monoidal enriched adjunction.
    \qed
  \end{enumerate}
\end{lemma}

For the rest of this section,
we recall the properties of tensorial strengths 
associated to enriched functors
between tensored enriched categories.
See \cite[\textsection 3]{SegrtRatkovic2013} for a detailed explanation.
An enriched category $\mathcal{C}$ is called \emph{tensored}
if for each object $X$ in $\mathcal{C}$
the Hom enriched functor
$\mathcal{C}(X,-):\mathcal{C}\to \mathcal{V}$
admits a left enriched adjoint
$-\otimes X:\mathcal{V}\to \mathcal{C}$.
We denote the components of
the adjunction unit, counit enriched natural transformations
at $z\in\mathcal{V}$, $Y\in \mathcal{C}$ as
\begin{equation}\label{eq Preliminary tensoredenrichedcat}
  \vcenter{\hbox{
    \xymatrix{
      \mathcal{V}
      \ar@<1ex>[r]^-{-\otimes X}
      \ar@{}[r]|-{\bot}
      &\mathcal{C}
      \ar@<1ex>[l]^-{\mathcal{C}(X,-)}
    }
  }}
  \qquad
  \textit{Cv}_{z,X}:z\to \mathcal{C}(X,z\otimes X)
  ,
  \qquad
  \textit{Ev}_{X,Y}:\mathcal{C}(X,Y)\otimes X\to Y
  .
\end{equation}
We have an isomorphism
$\imath_X:c\otimes X\xrightarrow{\cong} X$ in $\mathcal{C}$
which corresponds to the morphism
$\I_X:c\to \mathcal{C}(X,X)$ in $\mathcal{V}$.

Let 
$\mathcal{F}:\mathcal{C}\to \mathcal{D}$
be an enriched functor
between tensored enriched categories
$\mathcal{C}$, $\mathcal{D}$.
The \emph{tensorial strength} associated to $\mathcal{F}$
is the enriched natural transformation
$t^{\mathcal{F}}:
-\otimes \mathcal{F}(-)\Rightarrow \mathcal{F}(-\otimes -)$
between enriched functors from $\mathcal{V}\times \mathcal{C}$ to $\mathcal{D}$,
whose component at 
$z\in \mathcal{V}$, $X\in \mathcal{C}$ is the unique morphism
$t^{\mathcal{F}}_{z,X}:z\otimes \mathcal{F}(X)\to \mathcal{F}(z\otimes X)$
in $\mathcal{D}$ which corresponds to the 
morphism
\begin{equation*}
  \vcenter{\hbox{
    \xymatrix@C=35pt{
      z
      \ar[r]^-{\textit{Cv}_{z,X}}
      &\mathcal{C}(X,z\otimes X)
      \ar[r]^-{\mathcal{F}_{X,z\otimes X}}
      &\mathcal{D}(\mathcal{F}(X),\mathcal{F}(z\otimes X))
    }
  }}
\end{equation*}
in $\mathcal{V}$.
We say the given enriched functor
$\mathcal{F}$
\emph{preserves tensored objects}
if the associated tensorial strength
$t^{\mathcal{F}}$ is an enriched natural isomorphism.

Let $X$ be an object in a tensored enriched category $\mathcal{C}$.
We denote the tensorial strength associated to the cocontinuous enriched functor
$-\otimes X:\mathcal{V}\to \mathcal{C}$ as
\begin{equation*}
  a_{w,z,X}
  =t^{-\otimes X}_{w,z}
  :w\otimes (z\otimes X)\xrightarrow[]{\cong}(w\otimes z)\otimes X,
  \qquad
  w,z\in \mathcal{V}
\end{equation*}
which is an isomorphism in $\mathcal{C}$.
We often omit $a_{w,z,X}$
and simply write as $w\otimes z\otimes X$.

Let $\mathcal{C}$, $\mathcal{D}$ be tensored enriched categories.
For each enriched functor
$\mathcal{F}:\mathcal{C}\to \mathcal{D}$,
the associated tensorial strength $t^{\mathcal{F}}$
satisfies the following relations for all $w$, $z\in \mathcal{V}$, $X\in \mathcal{C}$.
\begin{equation}\label{eq Intro strengthCfunct}
  \hspace*{-0.5cm}
  \vcenter{\hbox{
    \xymatrix@R=20pt@C=10pt{
      c\otimes \mathcal{F}(X)
      \ar[r]^(0.5){t^{\mathcal{F}}_{c,X}}
      \ar@/_1pc/[dr]_-{\imath_{\mathcal{F}(X)}}^-{\cong}
      &\mathcal{F}(c\otimes X)
      \ar[d]^-{\mathcal{F}(\imath_X)}_-{\cong}
      \\
      \text{ }
      &\mathcal{F}(X)
    }
  }}
  \text{ }
  \vcenter{\hbox{
    \xymatrix@R=20pt@C=18pt{
      w\otimes (z\otimes \mathcal{F}(X))
      \ar[d]_-{a_{w,z,\mathcal{F}(X)}}^-{\cong}
      \ar[r]^(0.5){\I_w\otimes t^{\mathcal{F}}_{z,X}}
      &w\otimes \mathcal{F}(z\otimes X)
      \ar[r]^-{t^{\mathcal{F}}_{w,z\otimes X}}
      &\mathcal{F}(w\otimes (z\otimes X))
      \ar[d]^-{\mathcal{F}(a_{w,z,X})}_-{\cong}
      \\
      (w\otimes z)\otimes \mathcal{F}(X)
      \ar[rr]^-{t^{\mathcal{F}}_{w\otimes z,X}}
      &\text{ }
      &\mathcal{F}((w\otimes z)\otimes X)
    }
  }}
\end{equation}
Conversely, suppose we have a functor
$\mathcal{F}:\mathcal{C}\to \mathcal{D}$ between the underlying categories
and a natural transformation
$t^{\mathcal{F}}:-\otimes \mathcal{F}(-)\Rightarrow \mathcal{F}(-\otimes -)$
between functors from $\mathcal{V}\times \mathcal{C}$ to $\mathcal{D}$
that satisfies the relations (\ref{eq Intro strengthCfunct}).
Then there is a unique enrichment of $\mathcal{F}$
whose associated tensorial strength is equal to
the given natural transformation $t^{\mathcal{F}}_{z,X}$.

Let $\mathcal{F}$, $\mathcal{G}:\mathcal{C}\to \mathcal{D}$
be enriched functors between tensored enriched categories $\mathcal{C}$, $\mathcal{D}$.
For each enriched natural transformation
$\xi:\mathcal{F}\Rightarrow \mathcal{G}$,
the following relation holds for $z\in \mathcal{V}$, $X\in\mathcal{C}$.
\begin{equation} \label{eq Intro strengthCnat}
  \vcenter{\hbox{
    \xymatrix@R=20pt{
      z\otimes \mathcal{F}(X)
      \ar[d]_-{\I_z\otimes \xi_X}
      \ar[r]^-{t^{\mathcal{F}}_{z,X}}
      &\mathcal{F}(z\otimes X)
      \ar[d]^-{\xi_{z\otimes X}}
      \\
      z\otimes \mathcal{G}(X)
      \ar[r]^-{t^{\mathcal{G}}_{z,X}}
      &\mathcal{G}(z\otimes X)
    }
  }}
\end{equation}
Conversely, a natural transformation
$\xi:\mathcal{F}\Rightarrow\mathcal{G}$
between the underlying functors
becomes an enriched natural transformation
if and only if it satisfies the relation (\ref{eq Intro strengthCnat}) for all 
$z\in \mathcal{V}$, $X\in\mathcal{C}$.

\section{Eilenberg-Watts theorem in the enriched setting}
\label{sec EnEWthm}

A \emph{monoid} in $\mathcal{V}$
is a tuple $(b,u_b,m_b)$
of an object $b$ in $\mathcal{V}$
and unit, product morphisms $u_b:c\to b$, $m_b:b\otimes b\to b$
in $\mathcal{V}$ satisfying the associativity, unital relations.
We will abbreviate the unit, product morphisms
and denote a monoid $(b,u_b,m_b)$ in $\mathcal{V}$ as $b$.
Throughout this section,
we fix a monoid $b$ in $\mathcal{V}$
and denote
\begin{equation*}
  \textit{Mod}_b
\end{equation*}
as the complete, cocomplete enriched category
of right $b$-modules.
An object in $\textit{Mod}_b$
is a pair $(x,\gamma_x)$
of an object $x$ in $\mathcal{V}$,
and a morphism
$\gamma_x:x\otimes b\to x$
in $\mathcal{V}$ that satisfies the right $b$-action relations.
We will omit the right $b$-action and simply denote
a right $b$-module $(x,\gamma_x)$ as $x$.
Let $x$, $y$ be right $b$-modules.
The Hom object is given by the equalizer
\begin{equation*}
  \xymatrix@C=40pt{
    \textit{Mod}_{b}(x,y)
    \ar@{^{(}->}[r]^-{\textit{eq}_{x,y}}
    &\mathcal{V}(x,y)
    \ar@<0.5ex>[rr]^-{\mathcal{V}(\gamma_x,\I_y)}
    \ar@<-0.5ex>[rr]_-{\mathcal{V}(\I_x,\gamma_y)\circ (-\otimes b)_{x,y}}
    &\text{ }
    &\mathcal{V}(x\otimes b,y)
  }
\end{equation*}
and for $z\in \mathcal{V}$,
the tensored object with $x$ is 
$z\otimes x$ whose right $b$-action is
$\I_z\otimes \gamma_x$.
We can see $b$ itself as a right $b$-module,
whose right $b$-action is the product morphism $m_b$.
The right $b$-action morphism $\gamma_x:x\otimes b\to x$ of $x$
becomes a morphism in $\textit{Mod}_b$.

For the rest of this section,
let $\mathcal{C}$ be a cocomplete enriched category.
A \emph{left $b$-module object in $\mathcal{C}$}
is a pair $(X,\rho_X)$
of an object $X$ in $\mathcal{C}$,
and a morphism $\rho_X:b\otimes X\to X$ in $\mathcal{C}$
that satisfies the left $b$-action relations.
We omit the left $b$-action morphism and denote
a left $b$-module object $(X,\rho_X)$ in $\mathcal{C}$
as $X$.
The morphism in $\mathcal{V}$
corresponding to $\rho_X:b\otimes X\to X$
is denoted as 
$\bar{\rho}_X:b\to \mathcal{C}(X,X)$.
A morphism $X\to Y$
of left $b$-module objects in $\mathcal{C}$
is a morphism in $\mathcal{C}$
which is compatible with left $b$-actions $\rho_X$, $\rho_Y$.
We denote 
\begin{equation*}
  \text{ }\!\!_{b}\mathcal{C}
\end{equation*}
as the category of left $b$-module objects in $\mathcal{C}$.

\begin{lemma}
  Let $X$ be an object in $\text{ }\!\!_b\mathcal{C}$.
  We have an enriched adjunction
  \begin{equation} \label{eq EnEWThm bX adjunction}
    \vcenter{\hbox{
      \xymatrix@C=40pt{
        \textit{Mod}_b
        \ar@<1.2ex>[r]^-{-\tensor X}
        \ar@{}[r]|-{\bot}
        &\mathcal{C}
        \ar@<1.2ex>[l]^-{\mathcal{C}(X,-)}
      }
    }}
  \end{equation}
  which is a version of a Hom-tensor adjunction.
  \begin{enumerate}
    \item 
    The left enriched adjoint
    $-\tensor X:\textit{Mod}_b\to \mathcal{C}$
    sends each $z\in\textit{Mod}_b$
    to the following coequalizer in $\mathcal{C}$.
    \begin{equation*}
      \xymatrix@C=40pt{
        z\otimes b\otimes X
        \ar@<0.5ex>[r]^-{\gamma_z\otimes \I_X}
        \ar@<-0.5ex>[r]_-{\I_z\otimes \rho_X}
        &z\otimes X
        \ar@{->>}[r]^-{\textit{cq}_{z,X}}
        &z\tensor X
      }
    \end{equation*}
  
    \item 
    The right enriched adjoint
    $\mathcal{C}(X,-):\mathcal{C}\to\textit{Mod}_b$
    sends each $Y\in \mathcal{C}$ to the object
    $\mathcal{C}(X,Y)$ in $\mathcal{V}$ equipped with
    the right $b$-action morphism

    \begin{equation*}
      \gamma_{\mathcal{C}(X,Y)}:\!
      \xymatrix@C=40pt{
        \mathcal{C}(X,Y)\otimes b
        \ar[r]^-{\I_{\mathcal{C}(X,Y)}\otimes \bar{\rho}_X}
        &\mathcal{C}(X,Y)\otimes \mathcal{C}(X,X)
        \ar[r]^-{\mu_{X,X,Y}}
        &\mathcal{C}(X,Y)
        .
      }
    \end{equation*}
  \end{enumerate}
\end{lemma}
\begin{proof}
  One can check that
  $-\tensor X:\textit{Mod}_b\to \mathcal{C}$
  is a well-defined enriched functor,
  whose tensorial strength 
  $a_{w,z,X}:=t^{-\tensor X}_{w,z}$
  at $w\in \mathcal{V}$,
  $z\in \textit{Mod}_b$
  is the (unique) isomorphism
  in $\mathcal{C}$ satisfying the relation below.
  \begin{equation*}
    \xymatrix@R=22pt@C=40pt{
      w\otimes (z\otimes X)
      \ar@{->>}[d]_-{\I_w\otimes \textit{cq}_{z,X}}
      \ar[r]^-{a_{w,z,X}}_-{\cong}
      &(w\otimes z)\otimes X
      \ar@{->>}[d]^-{\textit{cq}_{w\otimes z,X}}
      \\
      w\otimes (z\tensor X)
      \ar@{.>}[r]^-{\exists!\text{ } a_{w,z,X}}_-{\cong}
      &(w\otimes z)\tensor X
    }
  \end{equation*}
  We also leave for the readers to check that
  $\mathcal{C}(X,-):\mathcal{C}\to \textit{Mod}_b$
  is a well-defined enriched functor.
  The tensorial strength
  $t^{\mathcal{C}(X,-)}_{w,Y}:
  w\otimes \mathcal{C}(X,Y)
  \to \mathcal{C}(X,w\otimes Y)$
  associated to the enriched functor
  $\mathcal{C}(X,-):\mathcal{C}\to \mathcal{V}$
  at $w\in \mathcal{V}$, $Y\in\mathcal{C}$
  becomes a morphism in $\textit{Mod}_b$,
  which is the tensorial strength
  associated to 
  $\mathcal{C}(X,-):\mathcal{C}\to \textit{Mod}_b$.
  The component of the adjunction unit at each $z\in\textit{Mod}_b$
  is the morphism
  $\eta_z:z\to \mathcal{C}(X,z\tensor X)$
  in $\textit{Mod}_b$,
  whose underlying morphism in $\mathcal{V}$ is
  \begin{equation*}
    \eta_z:
    \xymatrix@C=40pt{
      z
      \ar[r]^-{\textit{Cv}_{z,X}}
      &\mathcal{C}(X,z\otimes X)
      \ar[r]^-{\mathcal{C}(\I_X,\textit{cq}_{z,X})}
      &\mathcal{C}(X,z\tensor X)
      .
    }
  \end{equation*}
  The component of the adjunction counit
  at each $Y\in\mathcal{C}$ is the unique morphism
  $\varepsilon_Y:\mathcal{C}(X,Y)\tensor X\to Y$ in $\mathcal{C}$
  which satisfies the following relation.
  \begin{equation*}
    \xymatrix@R=22pt@C=40pt{
      \mathcal{C}(X,Y)\otimes X
      \ar@{->>}[d]_-{\textit{cq}_{\mathcal{C}(X,Y),X}}
      \ar@/^1.2pc/[dr]^-{\textit{Ev}_{X,Y}}
      \\
      \mathcal{C}(X,Y)\tensor X
      \ar@{.>}[r]^-{\exists!\text{ }\varepsilon_Y}
      &Y
    }
  \end{equation*}
  It is straight forward to show that $\eta$, $\varepsilon$
  are well-defined enriched natural transformations.
  One can verify the triangle identities of $\eta$, $\varepsilon$
  by pre-composing suitable epimorphisms in $\mathcal{C}$
  and reducing them to the triangle identities 
  of $\textit{Cv}$, $\textit{Ev}$
  (recall the enriched adjunction
  $-\otimes X\dashv \mathcal{C}(X,-)$
  in (\ref{eq Preliminary tensoredenrichedcat})).
\end{proof}

For each $X\in\text{ }\!\!_b\mathcal{C}$,
we have an isomorphism
$\imath_X:b \tensor X\xrightarrow[]{\cong}X$
in $\mathcal{C}$
which is the unique morphism satisfying the following relation.
\begin{equation} \label{eq EnEWThm imathbbX}
  \vcenter{\hbox{
    \xymatrix@R=22pt@C=40pt{
      b\otimes X
      \ar@/^1pc/[dr]^-{\rho_X}
      \ar@{->>}[d]_-{\textit{cq}_{b,X}}
      \\
      b\tensor X
      \ar@{.>}[r]^-{\exists!\text{ }\imath_X}_-{\cong}
      &X
    }
  }}
\end{equation}

Each morphism $X\to Y$ in $\text{ }\!\!_b\mathcal{C}$
determines an enriched natural transformation
$-\tensor X\Rightarrow -\tensor Y$
between enriched functors from
$\textit{Mod}_b$ to $\mathcal{C}$.
Thus we obtain a functor
\begin{equation}\label{eq EnEWThm leftadj}
  \vcenter{\hbox{
    \xymatrix@R=0pt{
      \text{ }\!\!_b\mathcal{C}
      \text{ }
      \ar[r]
      &\text{ }
      \mathcal{V}\textit{-Funct}(\textit{Mod}_b,\mathcal{C})
      \\
      X
      \text{ }
      \ar@{|->}[r]
      &\text{ }
      -\tensor X
      :\textit{Mod}_b\to \mathcal{C}
    }
  }}
\end{equation}
from $\text{ }\!\!_b\mathcal{C}$
to the category of enriched functors $\textit{Mod}_b\to \mathcal{C}$.
Meanwhile, we have the functor
\begin{equation} \label{eq EnEWThm rightadj}
  \vcenter{\hbox{
    \xymatrix@R=0pt{
      \mathcal{V}\textit{-Funct}(\textit{Mod}_b,\mathcal{C})
      \text{ }
      \ar[r]
      &\text{ }
      \text{ }\!\!_b\mathcal{C}
      \\
      \mathcal{F}:\textit{Mod}_b\to \mathcal{C}
      \text{ }\text{ }
      \ar@{|->}[r]
      &\text{ }\text{ }
      \mathcal{F}(b)
    }
  }}
\end{equation}
of evaluating at $b\in\textit{Mod}_b$.
Each enriched functor $\mathcal{F}:\textit{Mod}_b\to \mathcal{C}$
determines a left $b$-module object $\mathcal{F}(b)$ in $\mathcal{C}$
whose left $b$-action is
\begin{equation*}
  \xymatrix{
    \rho_{\mathcal{F}(b)}:
    b\otimes \mathcal{F}(b)
    \ar[r]^-{t^{\mathcal{F}}_{b,b}}
    &\mathcal{F}(b\otimes b)
    \ar[r]^-{\mathcal{F}(m_b)}
    &\mathcal{F}(b)
  }
\end{equation*}
and each enriched natural transformation 
$\xi:\mathcal{F}\Rightarrow\mathcal{G}$
between enriched functors
$\mathcal{F}$, $\mathcal{G}:\textit{Mod}_b\to \mathcal{C}$
determines a morphism
$\xi_b:\mathcal{F}(b)\to \mathcal{G}(b)$
in $\text{ }\!\!_b\mathcal{C}$.

\begin{proposition} \label{prop EnEWThm defining lambda}
  Let $b$ be a monoid in $\mathcal{V}$
  and let $\mathcal{C}$ be a cocomplete enriched category.
  For each enriched functor
  $\mathcal{F}:\textit{Mod}_b\to \mathcal{C}$
  we have an enriched natural transformation
  \begin{equation*}
    \lambda^{\mathcal{F}}:
    -\tensor \mathcal{F}(b)
    \Rightarrow
    \mathcal{F}
  \end{equation*}
  between enriched functors from $\textit{Mod}_b$ to $\mathcal{C}$,
  whose component at $z\in \textit{Mod}_b$
  is the unique morphism
  $\lambda^{\mathcal{F}}_z
  :z\tensor \mathcal{F}(b)\to \mathcal{F}(z)$
  in $\mathcal{C}$ which makes the following diagram commutative.
  \begin{equation} \label{eq EnEWThm defining lambda}
    \vcenter{\hbox{
      \xymatrix@R=20pt@C=70pt{
        z\otimes b\otimes \mathcal{F}(b)
        \ar@<-0.5ex>[d]_-{\gamma_z\otimes \I_{\mathcal{F}(b)}}
        \ar@<0.5ex>[d]^-{\I_z\otimes \rho_{\mathcal{F}(b)}}
        \ar[r]^-{t^{\mathcal{F}}_{z\otimes b,b}}
        &\mathcal{F}(z\otimes b\otimes b)
        \ar@<-0.5ex>[d]_-{\mathcal{F}(\gamma_z\otimes \I_{b})}
        \ar@<0.5ex>[d]^-{\mathcal{F}(\I_z\otimes m_b)}
        \\
        z\otimes \mathcal{F}(b)
        \ar@{->>}[d]_-{\textit{cq}_{z,\mathcal{F}(b)}}
        \ar[r]^-{t^{\mathcal{F}}_{z,b}}
        &\mathcal{F}(z\otimes b)
        \ar[d]^-{\mathcal{F}(\gamma_{z})}
        \\
        z\tensor \mathcal{F}(b)
        \ar@{.>}[r]^-{\exists!\text{ }\lambda^{\mathcal{F}}_z}
        &\mathcal{F}(z)
      }
    }}
  \end{equation}
  The component of $\lambda^{\mathcal{F}}$ at $b\in\textit{Mod}_b$
  is the isomorphism
  $\imath_{\mathcal{F}(b)}:
  b\tensor\mathcal{F}(b)
  \xrightarrow[]{\cong} \mathcal{F}(b)$
  in $\mathcal{C}$
  described in (\ref{eq EnEWThm imathbbX}).
  Moreover, the following are equivalent:
  \begin{enumerate}
    \item
    The enriched functor $\mathcal{F}:\textit{Mod}_b\to \mathcal{C}$
    is a left enriched adjoint;

    \item
    The enriched functor $\mathcal{F}:\textit{Mod}_b\to \mathcal{C}$
    is cocontinuous.
    
    \item 
    The enriched functor $\mathcal{F}:\textit{Mod}_b\to \mathcal{C}$
    preserves tensored objects,
    and the underlying functor of $\mathcal{F}$
    preserves coequalizers.

    \item 
    The enriched natural transformation
    $\lambda^{\mathcal{F}}:-\tensor\mathcal{F}(b)\Rightarrow\mathcal{F}$
    is an enriched natural isomorphism.
  \end{enumerate}
\end{proposition}
\begin{proof}
  For each $z\in \textit{Mod}_b$, 
  we have the diagram
  \begin{equation*}
    \xymatrix@R=18pt{
      z\otimes b\otimes \mathcal{F}(b)
      \ar[dd]_-{\gamma_z\otimes \I_{\mathcal{F}(b)}}
      \ar@{=}[r]
      &z\otimes b\otimes \mathcal{F}(b)
      \ar[dd]^-{t^{\mathcal{F}}_{z\otimes b,b}}
      \ar@{=}[r]
      &z\otimes b\otimes \mathcal{F}(b)
      \ar[d]^-{\I_z\otimes t^{\mathcal{F}}_{b,b}}
      \ar@{=}[r]
      &z\otimes b\otimes \mathcal{F}(b)
      \ar[dd]^-{I_z\otimes \rho_{\mathcal{F}(b)}}
      \\
      \text{ }
      &\text{ }
      &z\otimes \mathcal{F}(b\otimes b)
      \ar@/^0.5pc/[dl]|-{t^{\mathcal{F}}_{z,b\otimes b}}
      \ar[d]^-{\I_z\otimes \mathcal{F}(m_b)}
      &\text{ }
      \\
      z\otimes \mathcal{F}(b)
      \ar[d]_-{t^{\mathcal{F}}_{z,b}}
      &\mathcal{F}(z\otimes b\otimes b)
      \ar@/^0.5pc/[dl]|-{\mathcal{F}(\gamma_z\otimes \I_{b})}
      \ar@/_0.5pc/[dr]|-{\mathcal{F}(\I_z\otimes m_b)}
      &z\otimes \mathcal{F}(b)
      \ar[d]^-{t^{\mathcal{F}}_{z,b}}
      \ar@{=}[r]
      &z\otimes \mathcal{F}(b)
      \ar[d]^-{t^{\mathcal{F}}_{z,b}}
      \\
      \mathcal{F}(z\otimes b)
      \ar[d]_-{\mathcal{F}(\gamma_{z})}
      &\text{ }
      &\mathcal{F}(z\otimes b)
      \ar[d]^-{\mathcal{F}(\gamma_{z})}
      \ar@{=}[r]
      &\mathcal{F}(z\otimes b)
      \ar[d]^-{\mathcal{F}(\gamma_{z})}
      \\
      \mathcal{F}(z)
      \ar@{=}[rr]
      &\text{ }
      &\mathcal{F}(z)
      \ar@{=}[r]
      &\mathcal{F}(z)
    }
  \end{equation*}
  which shows that 
  the morphism $\lambda^{\mathcal{F}}_z$ in $\mathcal{C}$
  is well-defined, and the diagram (\ref{eq EnEWThm defining lambda}) commutes.
  To conclude that $\lambda^{\mathcal{F}}$ is an enriched natural transformation,
  it suffices to show that the following diagram commutes for
  $w\in \mathcal{V}$, $z\in\textit{Mod}_b$.
  \begin{equation*}
    \xymatrix@R=22pt@C=40pt{
      w\otimes (z\tensor \mathcal{F}(b))
      \ar[d]_-{\I_w\otimes \lambda^{\mathcal{F}}_{z}}
      \ar[r]^-{a_{w,z,\mathcal{F}(b)}}_-{\cong}
      &(w\otimes z)\tensor \mathcal{F}(b)
      \ar[d]^-{\lambda^{\mathcal{F}}_{w\otimes z}}
      \\
      w\otimes \mathcal{F}(z)
      \ar[r]^-{t^{\mathcal{F}}_{w,z}}
      &\mathcal{F}(w\otimes z)
    }
  \end{equation*}
  We obtain this by right-cancelling the epimorphism 
  $\I_w\otimes \textit{cq}_{z,\mathcal{F}(b)}$
  in the diagram below.
  \begin{equation*}
    \xymatrix@R=24pt@C=17pt{
      w\otimes (z\otimes \mathcal{F}(b))
      \ar@{->>}[d]_-{\I_w\otimes \textit{cq}_{z,\mathcal{F}(b)}}
      \ar@{=}[r]
      &w\otimes (z\otimes \mathcal{F}(b))
      \ar[dd]^-{a_{w,z,\mathcal{F}(b)}}_-{\cong}
      \ar@{=}[r]
      &w\otimes (z\otimes \mathcal{F}(b))
      \ar[d]^-{\I_w\otimes t^{\mathcal{F}}_{z,b}}
      \ar@{=}[r]
      &w\otimes (z\otimes \mathcal{F}(b))
      \ar@{->>}[d]^-{\I_w\otimes \textit{cq}_{z,\mathcal{F}(b)}}
      \\
      w\otimes (z\tensor \mathcal{F}(b))
      \ar[dd]_-{a_{w,z,\mathcal{F}(b)}}^-{\cong}
      &\text{ }
      &w\otimes \mathcal{F}(z\otimes b)
      \ar[d]^-{t^{\mathcal{F}}_{w,z\otimes b}}
      \ar@/^0.5pc/[dr]|(0.6){I_w\otimes \mathcal{F}(\gamma_{z})}
      &w\otimes (z\tensor \mathcal{F}(b))
      \ar[d]^-{\I_w\otimes \lambda^{\mathcal{F}}_{z}}
      \\
      \text{ }
      &(w\otimes z)\otimes \mathcal{F}(b)
      \ar@/^0.5pc/@{->>}[dl]|(0.45){\textit{cq}_{w\otimes z,\mathcal{F}(b)}}
      \ar[d]^-{t^{\mathcal{F}}_{w\otimes z,b}}
      &\mathcal{F}(w\otimes (z\otimes b))
      \ar@/^0.5pc/[dl]^(0.45){\mathcal{F}(a_{w,z,b})}_(0.45){\cong}
      \ar[dd]^-{\mathcal{F}(\I_w\otimes \gamma_z)}
      &w\otimes \mathcal{F}(z)
      \ar[dd]^-{t^{\mathcal{F}}_{w,z}}
      \\
      (w\otimes z)\tensor \mathcal{F}(b)
      \ar[d]_-{\lambda^{\mathcal{F}}_{w\otimes z}}
      &\mathcal{F}((w\otimes z)\otimes b)
      \ar[d]^-{\mathcal{F}(\gamma_{w\otimes z})}
      &\text{ }
      &\text{ }
      \\
      \mathcal{F}(w\otimes z)
      \ar@{=}[r]
      &\mathcal{F}(w\otimes z)
      \ar@{=}[r]
      &\mathcal{F}(w\otimes z)
      \ar@{=}[r]
      &\mathcal{F}(w\otimes z)
    }
  \end{equation*}
  This shows that $\lambda^{\mathcal{F}}$
  is a well-defined enriched natural transformation.
  From the defining relation
  (\ref{eq EnEWThm imathbbX}) of $\imath_{\mathcal{F}(b)}$,
  we see that
  $\lambda^{\mathcal{F}}_{b}
  =\imath_{\mathcal{F}(b)}:
  b\tensor\mathcal{F}(b)
  \xrightarrow[]{\cong} \mathcal{F}(b)$.

  We are left to show that the statements 1-4 are equivalent.
  The implications $1\Rightarrow 2\Rightarrow 3$ and $4\Rightarrow1$ are immediate.
  We prove $3\Rightarrow 4$ as follows.
  Since the enriched functor $\mathcal{F}:\textit{Mod}_b\to \mathcal{C}$
  preserves tensored objects,
  the top and middle horizontal morphisms in the diagram (\ref{eq EnEWThm defining lambda})
  are isomorphisms in $\mathcal{C}$.
  Since the underlying functor of $\mathcal{F}$
  preserves coequalizers,
  the right vertical morphisms in the diagram (\ref{eq EnEWThm defining lambda})
  is also a coequalizer diagram in $\mathcal{C}$.
  Thus we obtain that $\lambda^{\mathcal{F}}_{z}$
  is an isomorphism in $\mathcal{C}$ for every
  $z\in \textit{Mod}_b$.
\end{proof}

\begin{theorem} \label{thm Intro EnEW}
  Let $b$ be a monoid in $\mathcal{V}$
  and let $\mathcal{C}$ be a cocomplete enriched category.
  \begin{enumerate}
    \item
    We have a fully faithful left adjoint functor
    \begin{equation} \label{eq Intro EnEW}
      \vcenter{\hbox{
        \xymatrix@R=0pt{
          \text{ }\!\!_{b}\mathcal{C}
          \ar[r]
          &\mathcal{V}\textit{-Funct}(\textit{Mod}_b,\mathcal{C})
          \\
          X
          \ar@{|->}[r]
          &-\tensor X:\textit{Mod}_b\to \mathcal{C}
        }
      }}
    \end{equation}
    from the category of left $b$-module objects in $\mathcal{C}$
    to the category of enriched functors
    $\textit{Mod}_b\to\mathcal{C}$.

    \item
    The essential image of the functor (\ref{eq Intro EnEW})
    is the coreflective full subcategory
    $\mathcal{V}\textit{-Funct}_{\textit{cocon}}(\textit{Mod}_b,\mathcal{C})$
    of cocontinuous enriched functors
    $\textit{Mod}_b\to \mathcal{C}$,
    and we have an adjoint equivalence of categories    
    \begin{equation*}
      \vcenter{\hbox{
        \xymatrix@R=0pt{
          \text{ }\!\!_{b}\mathcal{C}
          \ar@<0.5ex>[r]^-{\simeq}
          &\mathcal{V}\textit{-Funct}_{\textit{cocon}}(\textit{Mod}_b,\mathcal{C}).
          \ar@<0.5ex>[l]^-{\simeq}
        }
      }}
    \end{equation*}
  \end{enumerate}
\end{theorem}
\begin{proof}
  We first show statement 1.
  We claim that the functor
  $\text{ }\!\!_b\mathcal{C}\to \mathcal{V}\textit{-Funct}(\textit{Mod}_b,\mathcal{C})$
  in (\ref{eq EnEWThm leftadj})
  is left adjoint to the functor
  $\mathcal{V}\textit{-Funct}(\textit{Mod}_b,\mathcal{C})\to \text{ }\!\!_b\mathcal{C}$
  of evaluating at $b\in\textit{Mod}_b$
  explained in (\ref{eq EnEWThm rightadj}).
  \begin{equation}\label{eq EnEWProof adjunction}
    \vcenter{\hbox{
      \xymatrix@C=40pt{
        \text{ }\!\!_b\mathcal{C}
        \ar@<0.7ex>[r]^-{(\ref{eq EnEWThm leftadj})}
        \text{ }
        &\text{ }
        \mathcal{V}\textit{-Funct}(\textit{Mod}_b,\mathcal{C})
        \ar@<0.7ex>[l]^-{(\ref{eq EnEWThm rightadj})}
      }
    }}
  \end{equation}
  For each $X\in \text{ }\!\!_b\mathcal{C}$,
  the functor (\ref{eq EnEWThm rightadj})
  sends the enriched functor $-\tensor X:\textit{Mod}_b\to\mathcal{C}$
  to $b\tensor X\in \text{ }\!\!_b\mathcal{C}$
  whose left $b$-action morphism is
  $\rho_{b\tensor X}=(m_b\tensor \I_X)\circ a_{b,b,X}:
  b\otimes (b \tensor X)
  \to b\tensor X$.
  The isomorphism
  $\imath_X:b\tensor X\xrightarrow{\cong} X$
  in $\mathcal{C}$
  introduced in (\ref{eq EnEWThm imathbbX})
  becomes an isomorphism
  in $\text{ }\!\!_b\mathcal{C}$,
  and one can check that it is natural in variable $X\in\text{ }\!\!_b\mathcal{C}$.
  The component of the unit of the adjunction (\ref{eq EnEWProof adjunction})
  at $X\in \text{ }\!\!_b\mathcal{C}$
  is the isomorphism
  $\imath_X^{-1}
  :X\xrightarrow{\cong}b\tensor X$ in $\text{ }\!\!_b\mathcal{C}$.
  
  The component of the counit of the adjunction
  (\ref{eq EnEWProof adjunction})
  at each enriched functor
  $\mathcal{F}:\textit{Mod}_b\to\mathcal{C}$
  is the enriched natural transformation
  $\lambda^{\mathcal{F}}:-\tensor\mathcal{F}(b)\Rightarrow \mathcal{F}$
  defined in Proposition~\ref{prop EnEWThm defining lambda}.
  We can check that $\lambda^{\mathcal{F}}$ is natural in variable $\mathcal{F}$ as follows:
  for each enriched natural transformation
  $\xi:\mathcal{F}\Rightarrow\mathcal{G}$
  between enriched functors
  $\mathcal{F}$, $\mathcal{G}:\textit{Mod}_b\to \mathcal{C}$
  and for each $z\in\textit{Mod}_b$,
  we have the relation
  \begin{equation} \label{eq EnEWProof lambda naturalinF}
    \vcenter{\hbox{
      \xymatrix@R=20pt@C=40pt{
        z\tensor \mathcal{F}(b)
        \ar[d]_-{\I_{z}\tensor \xi_{b}}
        \ar[r]^-{\lambda^{\mathcal{F}}_{z}}
        &\mathcal{F}(z)
        \ar[d]^-{\xi_{z}}
        \\
        z\tensor \mathcal{G}(b)
        \ar[r]^-{\lambda^{\mathcal{G}}_{z}}
        &\mathcal{G}(z)
      }
    }}
  \end{equation}
  which is obtained by right-cancelling the epimorphism
  $\textit{cq}_{z,\mathcal{F}(b)}$ in the diagram below.
  \begin{equation*}
    \xymatrix@R=20pt@C=35pt{
      z\otimes \mathcal{F}(b)
      \ar@{->>}[d]_-{\textit{cq}_{z,\mathcal{F}(b)}}
      \ar@{=}[r]
      &z\otimes \mathcal{F}(b)
      \ar[d]^-{t^{\mathcal{F}}_{z,b}}
      \ar@{=}[r]
      &z\otimes \mathcal{F}(b)
      \ar[d]^-{\I_z\otimes \xi_{b}}
      \ar@{=}[r]
      &z\otimes \mathcal{F}(b)
      \ar@{->>}[d]^-{\textit{cq}_{z,\mathcal{F}(b)}}
      \\
      z\tensor \mathcal{F}(b)
      \ar[d]_-{\lambda^{\mathcal{F}}_{z}}
      &\mathcal{F}(z\otimes b)
      \ar@/^0.5pc/[dl]|(0.45){\mathcal{F}(\gamma_{z})}
      \ar[d]^-{\xi_{z\otimes b}}
      &z\otimes \mathcal{G}(b)
      \ar@/^0.5pc/[dl]|(0.45){t^{\mathcal{G}}_{z,b}}
      \ar@/_0.5pc/@{->>}[dr]|-{\textit{cq}_{z,\mathcal{G}(b)}}
      &z\tensor \mathcal{F}(b)
      \ar[d]^-{\I_{z}\tensor\xi_{b}}
      \\
      \mathcal{F}(z)
      \ar[d]_-{\xi_{z}}
      &\mathcal{G}(z\otimes b)
      \ar[d]^-{\mathcal{G}(\gamma_{z})}
      &\text{ }
      &z\tensor \mathcal{G}(b)
      \ar[d]^-{\lambda^{\mathcal{G}}_{z}}
      \\
      \mathcal{G}(z)
      \ar@{=}[r]
      &\mathcal{G}(z)
      \ar@{=}[rr]
      &\text{ }
      &\mathcal{G}(z)
    }
  \end{equation*}
  Next, we show that the unit, counit satisfy the triangle identities.
  \begin{equation}
    \label{eq EnEWProof triangle identities}
    \vcenter{\hbox{
      \xymatrix@R=20pt@C=35pt{
        \mathcal{F}(b)
        \ar@/_1pc/@{=}[dr]
        \ar[r]^-{\imath_{\mathcal{F}(b)}^{-1}}_-{\cong}
        &b\tensor \mathcal{F}(b)
        \ar[d]^-{\lambda^{\mathcal{F}}_{b}}
        &-\tensor X
        \ar@/_1pc/@{=}[dr]
        \ar@2{->}[r]^-{-\tensor \imath_X^{-1}}_-{\cong}
        &-\tensor (b\tensor X)
        \ar@2{->}[d]^-{\lambda^{-\tensor X}}
        \\
        \text{ }
        &\mathcal{F}(b)
        &\text{ }
        &-\tensor X
      }
    }}
  \end{equation}
  The triangle identity on the left side of (\ref{eq EnEWProof triangle identities}) is already verified in Proposition~\ref{prop EnEWThm defining lambda}.
  We obtain the other triangle identity by right-cancelling the epimorphisms
  $\I_z\tensor \textit{cq}_{b,X}$
  and 
  $\textit{cq}_{z,b\tensor X}$
  in the diagram below.
  Let $z=(z,\gamma_z)\in \textit{Mod}_b$.
  \begin{equation*}
    \xymatrix@R=18pt@C=22pt{
      z\otimes (b\otimes X)
      \ar@{->>}[d]_-{\I_z\otimes \textit{cq}_{b,X}}
      \ar@{=}[r]
      &z\otimes (b\otimes X)
      \ar[d]^-{a_{z,b,X}}_-{\cong}
      \ar@{=}[rr]
      &\text{ }
      &z\otimes (b\otimes X)
      \ar[dd]_-{\I_z\otimes \rho_X}
      \ar@{=}[r]
      &z\otimes (b\otimes X)
      \ar@{->>}[d]^-{\I_z\otimes \textit{cq}_{b,X}}
      \\
      z\otimes (b\tensor X)
      \ar@{->>}[d]_-{\textit{cq}_{z,b\tensor X}}
      \ar@/^0.5pc/[dr]^(0.45){a_{z,b,X}}_(0.55){\cong}
      &(z\otimes b)\otimes X
      \ar@{->>}[d]^-{\textit{cq}_{z\otimes b,X}}
      \ar@/^0.5pc/[dr]^(0.55){\gamma_z\otimes \I_X}
      &\text{ }
      &\text{ }
      &z\otimes (b\tensor X)
      \ar@/_0.5pc/[dl]_(0.5){\I_z\otimes \imath_{X}}^-{\cong}
      \ar@{->>}[d]^-{\textit{cq}_{z,b\tensor X}}
      \\
      z\tensor (b\tensor X)
      \ar[d]_-{\lambda^{-\tensor X}_{z}}
      &(z\otimes b)\tensor X
      \ar[d]^-{\gamma_{z}\tensor \I_{X}}
      &z\otimes X
      \ar@{->>}[d]^-{\textit{cq}_{z,X}}
      &z\otimes X
      \ar@{->>}[d]^-{\textit{cq}_{z,X}}
      &z\tensor (b\tensor X)
      \ar[d]^-{\I_{z}\tensor\imath_{X}}_-{\cong}
      \\
      z\tensor X
      \ar@{=}[r]
      &z\tensor X
      \ar@{=}[r]
      &z\tensor X
      \ar@{=}[r]
      &z\tensor X
      \ar@{=}[r]
      &z\tensor X
    }
  \end{equation*}
  This shows that we have the adjunction
  (\ref{eq EnEWProof adjunction})
  as we claimed,
  whose adjunction unit is a natural isomorphism.
  Thus the left adjoint is fully faithful,
  and we obtain statement 1.

  Next we show statement 2.
  Let $\mathcal{F}:\textit{Mod}_b\to\mathcal{C}$
  be an enriched functor.
  The component $\lambda^{\mathcal{F}}$
  of the counit of the adjunction (\ref{eq EnEWProof adjunction})
  at $\mathcal{F}$ is an isomorphism in $\mathcal{V}\textit{-Funct}(\textit{Mod}_b,\mathcal{C})$
  if and only if $\mathcal{F}$
  is a cocontinuous enriched functor,
  by Proposition~\ref{prop EnEWThm defining lambda}.
  We conclude that 
  the essential image of the left adjoint
  $\text{ }\!\!_b\mathcal{C}
  \to \mathcal{V}\textit{-Funct}(\textit{Mod}_b,\mathcal{C})$
  in (\ref{eq EnEWProof adjunction})
  is the coreflective full subcategory
  $\mathcal{V}\textit{-Funct}_{\textit{cocon}}(\textit{Mod}_b,\mathcal{C})$
  of cocontinous enriched functors
  $\textit{Mod}_b\to\mathcal{C}$.
\end{proof}

For the rest of this section,
we examine the case when $\mathcal{C}$
is the enriched category of right modules over another monoid in $\mathcal{V}$.
By considering this special case,
we obtain a result
more similar to the original Eilenberg-Watts theorem.
This is also a preparation
for the contents in the upcoming sections.

Let $b'$ be another monoid in $\mathcal{V}$.
A left $b$-module object in $\textit{Mod}_{b'}$
is a \emph{$(b,b')$-bimodule},
which is a right $b'$-module $(x,\gamma'_x:x\otimes b'\to x)$
equipped with a left $b$-action morphism
$\rho_x:b\otimes x\to x$ in $\textit{Mod}_{b'}$.
We mostly omit
the right $b'$-action $\gamma'_x$
as well as
the left $b$-action $\rho_x$,
and denote a $(b,b')$-bimodule $(x,\gamma'_x,\rho_x)$ as $x$.
The category of $(b,b')$-bimodules
is denoted as
$\text{ }\!\!_b\textit{Mod}_{b'}$,
and for each $x\in \text{ }\!\!_b\textit{Mod}_{b'}$
the associated enriched adjunction (\ref{eq EnEWThm bX adjunction})
is denoted as follows.
\begin{equation} \label{eq EnEWThm bimodule adjunction}
  \vcenter{\hbox{
    \xymatrix@C=60pt{
      \textit{Mod}_b
      \ar@<1.2ex>[r]^-{-\tensor x}
      \ar@{}[r]|-{\bot}
      &\textit{Mod}_{b'}
      \ar@<1.2ex>[l]^-{\textit{Mod}_{b'}(x,-)}
    }
  }}
\end{equation}
The left enriched adjoint
$-\tensor x:\textit{Mod}_b\to \textit{Mod}_{b'}$
sends each $z\in \textit{Mod}_b$ to the coequalizer
\begin{equation*}
  \xymatrix@C=40pt{
    z\otimes b\otimes x
    \ar@<0.5ex>[r]^-{\gamma_z\otimes \I_x}
    \ar@<-0.5ex>[r]_-{\I_z\otimes \rho_x}
    &z\otimes x
    \ar@{->>}[r]^-{\textit{cq}_{z,x}}
    &z\tensor x
  }
\end{equation*}
in $\textit{Mod}_{b'}$,
and the tensorial strength at
$w\in \mathcal{V}$,
$z\in\textit{Mod}_b$
is the unique isomorphism
$a_{w,z,x}:=t^{-\tensor x}_{w,z}$ in $\textit{Mod}_{b'}$
which satisfies the relation below.
\begin{equation*}
  \xymatrix@R=18pt@C=40pt{
    w\otimes (z\otimes x)
    \ar@{->>}[d]_-{\I_w\otimes \textit{cq}_{z,x}}
    \ar[r]^-{a_{w,z,x}}_-{\cong}
    &(w\otimes z)\otimes x
    \ar@{->>}[d]^-{\textit{cq}_{w\otimes z,x}}
    \\
    w\otimes (z \tensor x)
    \ar@{.>}[r]^-{\exists!\text{ }a_{w,z,x}}_-{\cong}
    &(w\otimes z)\tensor x
  }
\end{equation*}
The following is a special case of Theorem~\ref{thm Intro EnEW},
which exhibits the relation between
bimodules and (cocontinuous) enriched functors
between enriched categories of right modules.



\begin{corollary} \label{cor EnEWthm bimodule}
  Let $b$, $b'$ be monoids
  in a B\'{e}nabou cosmos $\mathcal{V}$.
  \begin{enumerate}
    \item 
    We have a fully faithful left adjoint functor
    \begin{equation*}
      \vcenter{\hbox{
        \xymatrix@R=0pt{
          \text{ }\!\!_b\textit{Mod}_{b'}
          \ar[r]
          &\mathcal{V}\textit{-Funct}(\textit{Mod}_b,\textit{Mod}_{b'}\!)
          \\
          x
          \text{ }
          \ar@{|->}[r]
          &\text{ }
          -\tensor x:
          \textit{Mod}_b\to \textit{Mod}_{b'}
        }
      }}
    \end{equation*}
    from the category of $(b,b')$-bimodules
    to the category of enriched functors
    $\textit{Mod}_b\to \textit{Mod}_{b'}$.

    \item
    The essential image of this functor
    is the coreflective full subcategory
    $\mathcal{V}\textit{-Funct}_{\textit{cocon}}(\textit{Mod}_b,\textit{Mod}_{b'}\!)$
    of cocontinuous enriched functors
    $\textit{Mod}_b\to \textit{Mod}_{b'}$,
    and we have an adjoint equivalence of categories
    \begin{equation*}
      \vcenter{\hbox{
        \xymatrix@C=30pt{
          \text{ }\!\!_b\textit{Mod}_{b'}
          \ar@<0.5ex>[r]^-{\simeq}
          &\text{ }
          \mathcal{V}\textit{-Funct}_{\textit{cocon}}(\textit{Mod}_b,\textit{Mod}_{b'}\!).
          \ar@<0.5ex>[l]^-{\simeq}
        }
      }}
    \end{equation*}
  \end{enumerate}
\end{corollary}
\begin{proof}
  Substitute $\mathcal{C}=\textit{Mod}_{b'}$
  in Theorem~\ref{thm Intro EnEW}.
\end{proof}

\section{The six functors} \label{sec SixFunctors}
In this section, we describe the six functors
associated to each morphism of commutative monoids in $\mathcal{V}$.
A monoid $b$ in $\mathcal{V}$ is \emph{commutative}
if the product $m_b:b\otimes b\to b$ satisfies the commutativity relation.
Throughout this section, we fix a commutative monoid $b$ in $\mathcal{V}$.

Let $x$ be an object in $\textit{Mod}_b$.
We have a morphism
$\mathring{\gamma}_x:b\otimes x\to x$
in $\textit{Mod}_b$
whose underlying morphism in $\mathcal{V}$ is
$\mathring{\gamma}_x
:b\otimes x\xrightarrow[\cong]{s_{b,x}}x\otimes b\xrightarrow{\gamma_x} x$,
and $x$ equipped with $\mathring{\gamma}_x$
becomes a $(b,b)$-bimodule.
We denote the associated enriched adjunction
(\ref{eq EnEWThm bimodule adjunction}) as follows.
\begin{equation} \label{eq SixFunctors module adjunction}
  \vcenter{\hbox{
    \xymatrix@C=40pt{
      \textit{Mod}_b
      \ar@<1.2ex>[r]^-{-\tensor x}
      \ar@{}[r]|-{\bot}
      &\textit{Mod}_b
      \ar@<1.2ex>[l]^-{[x,-]}
    }
  }}
\end{equation}
The left enriched adjoint
$-\tensor x:\textit{Mod}_b\to \textit{Mod}_b$
sends each $z\in \textit{Mod}_b$ to the coequalizer
\begin{equation*}
  \xymatrix@C=40pt{
    z\otimes b\otimes x
    \ar@<0.5ex>[r]^-{\gamma_z\otimes \I_x}
    \ar@<-0.5ex>[r]_-{\I_z\otimes \mathring{\gamma}_x}
    &z\otimes x
    \ar@{->>}[r]^-{\textit{cq}_{z,x}}
    &z\tensor x
  }
\end{equation*}
in $\textit{Mod}_b$,
and the tensorial strength at
$w\in\mathcal{V}$, $z\in \textit{Mod}_b$
is the unique isomorphism 
$a_{w,z,x}:=t^{-\tensor x}_{w,z}$
in $\textit{Mod}_b$ which satisfies the following relation.
\begin{equation*}
  \xymatrix@C=45pt{
    w\otimes (z\otimes x)
    \ar@{->>}[d]_-{\I_w\otimes \textit{cq}_{z,x}}
    \ar[r]^-{a_{w,z,x}}_-{\cong}
    &(w\otimes z)\otimes x
    \ar@{->>}[d]^-{\textit{cq}_{w\otimes z,x}}
    \\
    w\otimes (z\tensor x)
    \ar@{.>}[r]^-{\exists!\text{ }a_{w,z,x}}_-{\cong}
    &(w\otimes z)\tensor x
  }
\end{equation*}

\begin{proposition} \label{prop SixFunctors MbtensorCcat}
  Let $b$ be a commutative monoid in a B\'{e}nabou cosmos $\mathcal{V}$.
  We have the closed monoidal enriched category
  $(\textit{Mod}_b,\tensor,b)$
  of right $b$-modules,
  which is complete and cocomplete as an enriched category.
  Let $x$, $y$, $z\in\textit{Mod}_b$.
  \begin{enumerate}
    \item 
    The coherence isomorphisms
    $a_{z,x,y}$, $\imath_x$, $\jmath_x$
    are unique isomorphisms in $\textit{Mod}_b$
    satisfying the following relations.
    \begin{equation*}
      \hspace*{-0.5cm}
      \vcenter{\hbox{
        \xymatrix@R=23pt@C=22pt{
          z\otimes (x\tensor y)
          \ar@{->>}[d]_-{\textit{cq}_{z,x\tensor y}}
          \ar[rr]^-{a_{z,x,y}}_-{\cong}
          &\text{ }
          &(z\otimes x)\tensor y
          \ar@{->>}[d]^-{\textit{cq}_{z,x}\tensor \I_{y}}
          &b\otimes x
          \ar@{->>}[d]_-{\textit{cq}_{b,x}}
          \ar@/^1pc/[dr]^-{\mathring{\gamma}_{x}}
          &\text{ }
          &x\otimes b
          \ar@{->>}[d]_-{\textit{cq}_{x,b}}
          \ar@/^1pc/[dr]^-{\gamma_{x}}
          &\text{ }
          \\
          z\tensor (x\tensor y)
          \ar@{.>}[rr]^-{\exists!\text{ }a_{z,x,y}}_-{\cong}
          &\text{ }
          &(z\tensor x)\tensor y
          &b\tensor x
          \ar@{.>}[r]^-{\exists!\text{ }\imath_{x}}_-{\cong}
          &x
          &x\tensor b
          \ar@{.>}[r]^-{\exists!\text{ }\jmath_{x}}_-{\cong}
          &x
        }
      }}
    \end{equation*}

    \item 
    The coherence isomorphism
    $s_{x,y}:x\tensor y\xrightarrow{\cong}y\tensor x$
    is the unique isomorphism in $\textit{Mod}_b$,
    whose underlying isomorphism
    in $\mathcal{V}$ satisfies the relation below.
    \begin{equation*}
      \xymatrix@R=23pt@C=40pt{
        x\otimes y
        \ar@{->>}[d]_-{\textit{cq}_{x,y}}
        \ar[r]^-{s_{x,y}}_-{\cong}
        &y\otimes x
        \ar@{->>}[d]^-{\textit{cq}_{y,x}}
        \\
        x\tensor y
        \ar@{.>}[r]^-{\exists!\text{ }s_{x,y}}_-{\cong}
        &y\tensor x
      }
    \end{equation*}
  \end{enumerate}
  In particular, the underlying closed symmetric monoidal category
  $\textit{Mod}_b$ of right $b$-modules is also a B\'{e}nabou cosmos.
  \end{proposition}
\begin{proof}
  The isomorphisms $a_{z,x,y}$, $s_{x,y}$ in $\textit{Mod}_b$
  are well-defined:
  see \textsection\!~\ref{subsec Appendix MbtensorCcat}.
  We leave for the readers to check that the isomorphisms
  $\imath_x$, $\jmath_x$ in $\textit{Mod}_b$ are well-defined,
  and the inverses of their underlying isomorphisms in $\mathcal{V}$ are
  \begin{equation*}
    \begin{aligned}
      \imath_x^{-1}
      &:
      \xymatrix@C=30pt{
        x
        \ar[r]^-{\imath^{-1}_x}_-{\cong}
        &c\otimes x
        \ar[r]^-{u_b\tensor \I_x}
        &b\otimes x
        \ar@{->>}[r]^-{\textit{cq}_{b,x}}
        &b\tensor x
        ,
      }
      \\
      \jmath_x^{-1}
      &:
      \xymatrix@C=30pt{
        x
        \ar[r]^-{\jmath^{-1}_x}_-{\cong}
        &x\otimes c
        \ar[r]^-{\I_x\tensor u_b}
        &x\otimes b
        \ar@{->>}[r]^-{\textit{cq}_{x,b}}
        &x\tensor b
        .
      }
    \end{aligned}
  \end{equation*}
  Now we verify the pentagon, unital, symmetry and hexagon axioms.
  The symmetry and unital axioms are explained in 
  \textsection\!~\ref{subsec Appendix MbtensorCcat}.
  One can check the pentagon axiom as follows.
  \begin{itemize}
    \item 
    Let $w\in \mathcal{V}$.
    The tensorial strengths associated to the enriched functors
    $-\otimes y:\mathcal{V}\to \textit{Mod}_b$
    and
    $-\tensor y:\textit{Mod}_b\to \textit{Mod}_b$
    satisfy the relations (\ref{eq Intro strengthCfunct}).
    In particular, we have the pentagon axioms for the quadruples
    $w$, $z$, $x\in \mathcal{V}$, $y\in\textit{Mod}_b$
    and
    $w$, $z\in \mathcal{V}$, $x$, $y\in \textit{Mod}_b$.

    \item 
    Let $w\in\textit{Mod}_b$.
    By pre-composing suitable epimorphisms,
    the pentagon axiom for the quadruple
    $w$, $z$, $x$, $y\in\textit{Mod}_b$
    reduces to the pentagon axiom of
    $w\in \mathcal{V}$, $z$, $x$, $y\in\textit{Mod}_b$,
    and further reduces to that of 
    $w$, $z\in\mathcal{V}$, $x$, $y\in \textit{Mod}_b$.
  \end{itemize}
  Similarly, we obtain the hexagon axiom in $\textit{Mod}_b$
  as it reduces to the hexagon relation in $\mathcal{V}$
  after pre-composing suitable epimorphisms.
  We omit the details.
  We conclude that $(\textit{Mod}_b,\tensor,b)$
  is a symmetric monoidal category.

  For each $z\in\textit{Mod}_b$,
  we also have an enriched functor
  $z\tensor -:\textit{Mod}_b\to \textit{Mod}_b$
  whose tensorial strength at
  $w\in \mathcal{V}$,
  $x\in\textit{Mod}_b$ is
  the unique isomorphism in $\textit{Mod}_b$
  satisfying the following relation.
  \begin{equation} \label{eq SixFunctors zbtensorb- strength}
    \vcenter{\hbox{
      \xymatrix@R=24pt@C=40pt{
        w\otimes z\otimes x
        \ar@{->>}[d]_-{\I_w\otimes \textit{cq}_{z,x}}
        \ar[r]^-{s_{w,z}\otimes \I_{x}}_-{\cong}
        &z\otimes w\otimes x
        \ar@{->>}[d]^-{\textit{cq}_{z,w\otimes x}}
        \\
        w\otimes (z\tensor x)
        \ar[r]^-{t^{z\tensor-}_{w,x}}_-{\cong}
        &z\tensor (w\otimes x)
      }
    }}
  \end{equation}
  The tensorial strength $t^{z\tensor -}_{w,x}$
  is explicitly given as
  \begin{equation*}
    t^{z\tensor -}_{w,x}:
    \xymatrix@C=30pt{
      w\otimes (z\tensor x)
      \ar[r]^(0.48){\I_w\otimes s_{z,x}}_-{\cong}
      &w\otimes (x\tensor z)
      \ar[r]^-{a_{w,x,z}}_-{\cong}
      &(w\otimes x)\tensor z
      \ar[r]^-{s_{w\otimes x,z}}_-{\cong}
      &z\tensor (w\otimes x)
      .
    }
  \end{equation*}
  The collection of enriched functors
  $z\tensor -$, $-\tensor x$
  give rise to an enriched functor
  $-\tensor -:\textit{Mod}_b\times\textit{Mod}_b\to \textit{Mod}_b$,
  and the collection of enriched adjunctions
  in (\ref{eq SixFunctors module adjunction})
  uniquely determines the enriched functor
  $[-,-]:\textit{Mod}_b^{\op}\times \textit{Mod}_b\to \textit{Mod}_b$.
  We leave for the readers to check that the coherence isomorphisms
  $a$, $\imath$, $\jmath$, $s$ of $\textit{Mod}_b$
  are enriched natural in every variable.
  We conclude that 
  $(\textit{Mod}_b,\tensor,b)$ is a closed monoidal enriched category.
\end{proof}

A morphism $b\to b'$ of commutative monoids in $\mathcal{V}$
is a morphism in $\mathcal{V}$ which is compatible with unit, product morphisms of $b$, $b'$.
We denote
$\textit{Comm}(\mathcal{V})$
as the category of commutative monoids in $\mathcal{V}$.
For each morphism $f:b\to b'$ in $\textit{Comm}(\mathcal{V})$,
we can see $b'$ as a $(b,b')$-bimodule and also as a $(b',b)$-bimodule
where $b'$-actions are the product $m_{b'}$
and $b$-actions are induced by $f$.

\begin{lemma} \label{lem SixFunctors f*f! Cadjuctions}
  Let $f:b\to b'$ be a morphism in $\textit{Comm}(\mathcal{V})$.
  We have a pair of enriched adjunctions
  \begin{equation*}
    \vcenter{\hbox{
      \xymatrix@C=60pt{
        \textit{Mod}_b
        \ar@<1.2ex>[r]^-{f^*=\text{ }\!-\text{ }\!\tensor b'}
        \ar@{}[r]|-{\bot}
        &\textit{Mod}_{b'}
        \ar@<1.2ex>[l]^-{f_*}
      }
    }}
    \qquad\quad
    \vcenter{\hbox{
      \xymatrix@C=60pt{
        \textit{Mod}_{b'}
        \ar@<1.2ex>[r]^-{f_!=\text{ }\!-\text{ }\!\tensor'b'}
        \ar@{}[r]|-{\bot}
        &\textit{Mod}_b
        \ar@<1.2ex>[l]^-{f^!=\text{ }\!\textit{Mod}_b(b',-)}
      }
    }}
  \end{equation*}
  and an enriched natural isomorphism
  $f_!\cong f_*:\textit{Mod}_{b'}\to \textit{Mod}_b$.
  \begin{enumerate}
    \item 
    $f^*=-\tensor b':\textit{Mod}_b\to \textit{Mod}_{b'}$
    is the left enriched adjoint
    of the enriched adjunction
    (\ref{eq EnEWThm bimodule adjunction})
    associated to the $(b,b')$-bimodule $b'$.
    
    \item 
    $f_*:\textit{Mod}_{b'}\to \textit{Mod}_b$
    sends each $y\in \textit{Mod}_{b'}$ to
    $f_*(y)=y\in\textit{Mod}_b$ whose right $b$-action is
    \begin{equation*}
      \vcenter{\hbox{
        \xymatrix{
          \gamma_y:y\otimes b
          \ar[r]^-{\I_y\otimes f}
          &y\otimes b'
          \ar[r]^-{\gamma'_y}
          &y
        }
      }}
    \end{equation*}
    and the tensorial strength of $f_*$ is the identity.
    
    \item 
    $f_!\dashv f^!$
    is the enriched adjunction (\ref{eq EnEWThm bimodule adjunction})
    associated to the $(b',b)$-bimodule $b'$.
  \end{enumerate}
\end{lemma}
\begin{proof}
  We leave for the readers to check that
  the enriched functor $f_*:\textit{Mod}_{b'}\to \textit{Mod}_b$
  is well-defined.
  For each $x\in \textit{Mod}_b$,
  the right $b'$-action morphism of
  $f^*(x)=x\tensor b'$
  is the unique morphism
  $\gamma'_{x\tensor b'}$
  in $\mathcal{V}$ satisfying the relation below.
  \begin{equation*}
    \xymatrix@R=20pt@C=50pt{
      x\otimes b'\otimes b'
      \ar@{->>}[d]_-{\textit{cq}_{x,b'}\otimes \I_{b'}}
      \ar[r]^-{\I_x\otimes m_{b'}}
      &x\otimes b'
      \ar@{->>}[d]^-{\textit{cq}_{x,b'}}
      \\
      (x\tensor b')\otimes b'
      \ar@{.>}[r]^(0.5){\exists!\text{ }\gamma'_{x\tensor b'}}
      &x\tensor b'
    }
  \end{equation*}
  One can check that the enriched functor
  $f_*f^*:\textit{Mod}_b\to \textit{Mod}_b$
  is equal to the left enriched adjoint
  $-\tensor b':\textit{Mod}_b\to \textit{Mod}_b$
  of the enriched adjunction
  (\ref{eq SixFunctors module adjunction})
  associated to the right $b$-module $b'$.
  
  Let us explain the enriched adjunction
  $f^*\dashv f_*$ in detail.
  The component 
  $\eta^*_x:x\to f_*f^*(x)$
  of the adjunction unit $\eta^*$
  at each $x\in\textit{Mod}_b$ is the morphism
  $\eta^*_x:x\xrightarrow[\cong]{\jmath_x^{-1}}x\tensor b\xrightarrow{\I_x\tensor f}x\tensor b'$
  in $\textit{Mod}_b$.
  The component 
  $\varepsilon^*_{y}:f^*f_*(y)\to y$
  of the adjunction counit $\varepsilon^*$
  at each $y\in \textit{Mod}_{b'}$
  is the unique morphism
  $\varepsilon^*_y:y\tensor b'\to y$
  in $\textit{Mod}_{b'}$
  whose underlying morphism in $\mathcal{V}$
  satsifies the relation below.
  \begin{equation*}
    \xymatrix@R=20pt@C=40pt{
      y\otimes b'
      \ar@{->>}[d]_-{\textit{cq}_{y,b'}}
      \ar@/^1pc/[dr]^-{\gamma'_y}
      \\
      y\tensor b'
      \ar@{.>}[r]^-{\exists!\text{ }\varepsilon^*_b}
      &y
    }
  \end{equation*}
  One can check that
  $\eta^*$, $\varepsilon^*$
  are well-defined enriched natural transformations,
  and satisfy one of the triangle identities
  $f_*\varepsilon^*\circ \eta^*f^*=\I_{f^*}$.
  We obtain the other triangle identity
  $\varepsilon^*f^*\circ f^*\eta^*=\I_{f^*}$
  by right-cancelling the epimorphism
  $\textit{cq}_{x,b'}$
  in the following diagram in $\mathcal{V}$.
  Let $x\in\textit{Mod}_b$.
  \begin{equation*}
    \xymatrix@R=22pt@C=20pt{
      x\otimes b'
      \ar@{->>}[d]_-{\textit{cq}_{x,b'}}
      \ar@{=}[r]
      &x\otimes b'
      \ar[ddd]_-{\jmath_{x}^{-1}\otimes \I_{b'}}^-{\cong}
      \ar@{=}[r]
      &x\otimes b'
      \ar[d]_-{\jmath_x^{-1}\otimes \I_{b'}}^-{\cong}
      \ar@{=}[rr]
      &\text{ }
      &x\otimes b'
      \ar@/_0.5pc/[ddl]_-{\I_x\otimes \imath_{b'}^{-1}}^-{\cong}
      \ar@{=}[ddddd]
      \\
      x\tensor b'
      \ar@/_1.8pc/@<-2.5ex>[dddd]|-{f^*(\eta^*_{x})}
      \ar[ddd]^-{\jmath_{x}^{-1}\tensor\I_{b'}}_-{\cong}
      &\text{ }
      &(x\otimes c)\otimes b'
      \ar[d]_-{(\I_x\otimes u_b)\otimes \I_{b'}}
      \ar@/^0.5pc/[dr]^(0.6){a^{-1}_{x,c,b'}}_(0.6){\cong}
      \ar@/^2pc/@<2.5ex>[dd]|(0.7){(\I_x\otimes u_{b'})\otimes \I_{b'}}
      &\text{ }
      &\text{ }
      \\
      \text{ }
      &\text{ }
      &(x\otimes b)\otimes b'
      \ar@/_1pc/@{->>}[dl]|(0.55){\textit{cq}_{x,b}\otimes \I_{b'}}
      \ar[d]_-{(\I_x\otimes f)\otimes \I_{b'}}
      &x\otimes (c\otimes b')
      \ar[dd]^-{\I_x\otimes (u_{b'}\otimes \I_{b'})}
      &\text{ }
      \\
      \text{ }
      &(x\tensor b)\otimes b'
      \ar@/_1pc/@{->>}[dl]|(0.55){\textit{cq}_{x\tensor b,b'}}
      \ar[d]^-{(\I_x\tensor f)\otimes \I_{b'}}
      &(x\otimes b')\otimes b'
      \ar@/^1pc/@{->>}[dl]|(0.45){\textit{cq}_{x,b'}\otimes \I_{b'}}
      \ar@/_0.5pc/[dr]_-{a^{-1}_{x,b',b'}}^-{\cong}
      &\text{ }
      &\text{ }
      \\
      (x\tensor b)\tensor b'
      \ar[d]^-{(\I_x\tensor f)\tensor\I_{b'}}
      &(x\tensor b')\otimes b'
      \ar@/^1pc/@{->>}[dl]|(0.45){\textit{cq}_{x\tensor b',b'}}
      \ar[dd]^-{\gamma'_{x\tensor b'}}
      &\text{ }
      &x\otimes (b'\otimes b')
      \ar@/_0.5pc/[dr]_-{\I_x\otimes m_{b'}}
      &\text{ }
      \\
      (x\tensor b')\tensor b'
      \ar[d]_-{\varepsilon^*_{f^*(x)}}
      &\text{ }
      &\text{ }
      &\text{ }
      &x\otimes b'
      \ar@{->>}[d]^-{\textit{cq}_{x,b'}}
      \\
      x\tensor b'
      \ar@{=}[r]
      &x\tensor b'
      \ar@{=}[rrr]
      &\text{ }
      &\text{ }
      &x\tensor b'
    }
  \end{equation*}
  This shows that the enriched adjunction $f^*\dashv f_*$
  is well-defined.
  We also have an equality of enriched functors
  $f_!=f_*(-\tensor' b'):\textit{Mod}_{b'}\to \textit{Mod}_b$
  and an enriched natural isomorphism
  $f_*(\jmath'):f_!\cong f_*$.
\end{proof}

Let $f:b\to b'$ be a morphism in $\textit{Comm}(\mathcal{V})$.
Then we can see $b'$
as a commutative monoid in 
$\textit{Mod}_b$.
The unit morphism in $\textit{Mod}_b$ is
$u_{b'/b}=f:b\to b'$,
and the product morphism
$m_{b'/b}:b'\tensor b'\to b'$
is the unique morphism in $\textit{Mod}_b$
whose underlying morphism in $\mathcal{V}$
satisfies the following relation.
\begin{equation*}
  \xymatrix@R=22pt@C=40pt{
    b'\otimes b'
    \ar@{->>}[d]_-{\textit{cq}_{b',b'}}
    \ar@/^1pc/[dr]^-{m_{b'}}
    \\
    b'\tensor b'
    \ar@{.>}[r]^-{\exists!\text{ }m_{b'/b}}
    &b'
  }
\end{equation*}
Conversely, for each commutative monoid
$(b',u_{b'/b},m_{b'/b})$
in $\textit{Mod}_b$,
we have a morphism
$f=u_{b'/b}:b\to b'$ in $\textit{Comm}(\mathcal{V})$,
where the unit, product morphisms of $b'$ are
\begin{equation*}
  \xymatrix{
    u_{b'}:c
    \ar[r]^-{u_b}
    &b
    \ar[r]^-{u_{b'/b}}
    &b',
  }
  \qquad
  \xymatrix{
    m_{b'}:
    b'\otimes b'
    \ar@{->>}[r]^-{\textit{cq}_{b',b'}}
    &b'\tensor b'
    \ar[r]^-{m_{b'/b}}
    &b'
    .
  }
\end{equation*}
Let us denote
$\textit{Comm}_b$
as the coslice category of $\textit{Comm}(\mathcal{V})$ under $b$
and call it as the \emph{category of commutative monoids in $\mathcal{V}$ over $b$}.
One can check that we have an isomorphism of categories
\begin{equation} \label{eq SixFunctors Commb=CommMb}
  \textit{Comm}_{b}
  \cong
  \textit{Comm}(\textit{Mod}_b).
\end{equation}

\begin{proposition} \label{prop SixFunctors laxtensorCadj}
  Let $f:b\to b'$
  be a morphism in $\textit{Comm}(\mathcal{V})$.
  The enriched adjunction
  \begin{equation*}
    \vcenter{\hbox{
      \xymatrix@C=40pt{
        \textit{Mod}_b
        \ar@<1.2ex>[r]^-{f^*}
        \ar@{}[r]|-{\bot}
        &\textit{Mod}_{b'}
        \ar@<1.2ex>[l]^-{f_*}
      }
    }}
  \end{equation*}
  in Lemma~\ref{lem SixFunctors f*f! Cadjuctions}
  becomes a lax monoidal enriched adjunction
  satisfying the projection formula.
  \begin{enumerate}
    \item 
    We have an isomorphism
    $(f^*)^0=\imath_{b'}^{-1}
    :b'\xrightarrow{\cong}f^*(b)$
    in $\textit{Mod}_{b'}$,
    and for each pair
    $x$, $y\in\textit{Mod}_b$
    we have a unique isomorphism
    $(f^*)^2_{x,y}:
    f^*(x)\tensor' f^*(y)\xrightarrow{\cong}
    f^*(x\tensor y)$
    in $\textit{Mod}_{b'}$
    whose underlying isomorphism in $\mathcal{V}$
    satisfies the relation below.
    \begin{equation*}
      \xymatrix@R=23pt@C=60pt{
        (x\tensor b')\otimes (y\tensor b')
        \ar@{->>}[d]_-{\textit{cq}_{x\tensor b',y\tensor b'}}
        \ar@{->>}[r]^-{\textit{cq}'_{x\tensor b',y\tensor b'}}
        &(x\tensor b')\tensor' (y\tensor b')
        \ar@{.>}[dd]^-{\exists!\text{ }(f^*)^2_{x,y}}_-{\cong}
        \\
        x\tensor b'\tensor y\tensor b'
        \ar[d]_-{\I_x\tensor s_{b',y}\tensor\I_{b'}}^-{\cong}
        &\text{ }
        \\
        x\tensor y\tensor b'\tensor b'
        \ar[r]^-{\I_{x\tensor y}\tensor m_{b'/b}}
        &x\tensor y\tensor b'
      }
    \end{equation*}
    
    \item 
    We have a morphism
    $(f_*)^0=f:b\to f_*(b')$
    in $\textit{Mod}_b$,
    and for each pair
    $x$, $y\in\textit{Mod}_{b'}$
    we have a unique morphism
    $(f_*)^2_{x,y}:
    f_*(x)\tensor f_*(y)\to f_*(x\tensor' y)$
    in $\textit{Mod}_b$
    whose underlying morphism in $\mathcal{V}$
    satisfies the relation below.
    \begin{equation*}
      \xymatrix@R=23pt@C=50pt{
        x\otimes y
        \ar@{->>}[r]^-{\textit{cq}_{x,y}}
        \ar@/_1pc/@{->>}[dr]_-{\textit{cq}'_{x,y}}
        &x\tensor y
        \ar@{.>}[d]^-{\exists!\text{ }(f_*)^2_{x,y}}
        \\
        \text{ }
        &x\tensor' y
      }
    \end{equation*}
  \end{enumerate}
\end{proposition}
\begin{proof}
  We first show that 
  $f_*:\textit{Mod}_{b'}\to \textit{Mod}_b$
  together with morphisms
  $(f_*)^0$, $(f_*)^2$
  is a lax monoidal enriched functor.
  Let $x$, $y\in\textit{Mod}_{b'}$.
  We leave for the readers to check that
  the morphism 
  $(f_*)^2$ in $\textit{Mod}_b$ is well-defined.
  The morphism $(f_*)^2_{x,y}$
  is enriched natural in variable $x\in\textit{Mod}_{b'}$:
  see \textsection\!~\ref{subsec Appendix laxtensorCadj}.
  We leave for the readers to check that
  $(f_*)^2$ is compatible with
  coherence morphisms $\imath$, $\jmath$, $s$ of $\textit{Mod}_b$
  and $\imath'$, $\jmath'$, $s'$ of $\textit{Mod}_{b'}$.
  In particular,
  $(f_*)^2_{x,y}$ is also enriched natural in variable $y\in\textit{Mod}_{b'}$.
  Furthermore,
  $(f_*)^2$
  is compatible with coherence isomorphisms
  $a$ of $\textit{Mod}_b$ and $a'$ of $\textit{Mod}_{b'}$:
  see \textsection\!~\ref{subsec Appendix laxtensorCadj}.
  This shows that $f_*$ is a well-defined lax monoidal enriched functor.

  Now we show that $f^*\dashv f_*$ is a lax monoidal enriched adjunction.
  Let $x$, $y\in \textit{Mod}_b$.
  One can check that 
  \begin{equation} \label{eq gamma'opgamma'}
    \begin{aligned}
      \gamma'_{x\tensor b'}
      &:
      \xymatrix@C=30pt{
        (x\tensor b')\otimes b'
        \ar@{->>}[r]^-{\textit{cq}_{x\tensor b',b'}}
        &x\tensor b'\tensor b'
        \ar[r]^-{\I_x\tensor m_{b'/b}}
        &x\tensor b'
        ,
      }
      \\
      \mathring{\gamma}'_{y\tensor b'}
      &:
      \xymatrix@C=30pt{
        b'\otimes (y\tensor b')
        \ar@{->>}[r]^-{\textit{cq}_{b',y\tensor b'}}
        &b'\tensor y\tensor b'
        \ar[r]^-{s_{b',y}\tensor \I_{b'}}_-{\cong}
        &y\tensor b'\tensor b'
        \ar[r]^-{\I_y\tensor m_{b'/b}}
        &y\tensor b'
      }
    \end{aligned}
  \end{equation}
  and that the morphism
  $(f^*)^2_{x,y}$ in $\mathcal{V}$ is well-defined.
  Recall the unit, counit 
  $\eta^*$, $\varepsilon^*$
  of the enriched adjunction
  $f^*\dashv f_*$
  described in the proof of Lemma~\ref{lem SixFunctors f*f! Cadjuctions}.
  By Lemma~\ref{lem Intro laxtensorCadj},
  the left enriched adjoint $f^*$
  becomes a colax monoidal enriched functor
  whose morphisms
  $(\mathring{f}^*)^0$, $(\mathring{f}^*)^2$
  in $\textit{Mod}_{b'}$
  are given as follows.
  \begin{equation*}
    \vcenter{\hbox{
      \xymatrix@C=35pt{
        f^*(b)
        \ar[d]_-{f^*\bigl((f_*)^0\bigr)}
        \ar@/^1pc/[dr]^-{(\mathring{f}^*)^0}
        \\
        f^*f_*(b')
        \ar[r]^-{\varepsilon^*_{b'}}
        &b'
      }
    }}
    \qquad
    \vcenter{\hbox{
      \xymatrix@R=30pt@C=40pt{
        f^*(x\tensor y)
        \ar[d]_-{f^*(\eta^*_x \tensor \eta^*_y)}
        \ar@/^2pc/[ddr]^-{(\mathring{f}^*)^2_{x,y}}
        \\
        f^*(f_*f^*(x)\tensor f_*f^*(y))
        \ar[d]_-{f^*\bigl((f_*)^2_{f^*(x),f^*(y)}\bigr)}
        \\
        f^*f_*(f^*(x)\tensor' f^*(y))
        \ar[r]^-{\varepsilon^*_{f^*(x)\tensor'f^*(y)}}
        &f^*(x)\tensor' f^*(y)
      }      
    }}
  \end{equation*}
  To conclude that
  $f^*\dashv f_*$
  is a lax monoidal enriched adjunction,
  it suffices to show that the morphisms
  $(\mathring{f}^*)^0$,
  $(\mathring{f}^*)^2$
  in $\mathcal{V}$ 
  are the inverses of
  $(f^*)^0$,
  $(f^*)^2$.
  We have
  \begin{equation*}
    \varepsilon^*_{b'}
    =m_{b'/b}:b'\tensor b'\to b'
    ,
    \qquad\quad
    (\mathring{f}^*)^0
    =\imath_{b'}
    :b\tensor b'\xrightarrow[]{\cong}b'
  \end{equation*}
  which shows that the morphisms
  $(\mathring{f}^*)^0$,
  $(f^*)^0$
  in $\mathcal{V}$ are inverse to each other.
  The morphisms 
  $(\mathring{f}^*)^2$,
  $(f^*)^2$
  in $\mathcal{V}$
  are also inverse to each other:
  see \textsection\!~\ref{subsec Appendix laxtensorCadj}.
  We conclude that $f^*\dashv f_*$ is a lax monoidal enriched adjunction.
  Moreover, $f^*\dashv f_*$ satisfies the projection formula:
  see \textsection\!~\ref{subsec Appendix laxtensorCadj}.
\end{proof}

Let
$b$, $b'$ be objects in $\textit{Comm}(\mathcal{V})$.
We have $b\otimes b'\in \textit{Comm}(\mathcal{V})$
which is a coproduct of $b$, $b'$
in the catgory $\textit{Comm}(\mathcal{V})$.
Every object in $\textit{Comm}_{b\otimes b'}$
is uniquely written as
$(f,f'):b\otimes b'\xrightarrow{f\otimes f'}\tilde{b}\otimes \tilde{b}\xrightarrow{m_{\tilde{b}}}\tilde{b}$
for some pair of morphisms
\begin{equation*}
  \vcenter{\hbox{
    \xymatrix@R=13pt@C=13pt{
      \text{ }
      &\tilde{b}
      &\text{ }
      \\
      b
      \ar[ur]^-{f}
      &\text{ }
      &b'
      \ar[ul]_-{f'}
    }
  }}
\end{equation*}
in $\textit{Comm}(\mathcal{V})$.
An object 
$(f,f'):b\otimes b'\to\tilde{b}$
in $\textit{Comm}_{b\otimes b'}$,
which we also denoted as
$b\xrightarrow{f}\tilde{b}\xleftarrow{f'}b'$,
determines a cocontinuous lax monoidal enriched functor
\begin{equation*}
  f'_*f^*:\textit{Mod}_b\to \textit{Mod}_{b'}
\end{equation*}
whose underlying cocontinuous enriched functor is equal to
the functor
$-\tensor \tilde{b}:\textit{Mod}_b\to \textit{Mod}_{b'}$
associated to the $(b,b')$-bimodule $\tilde{b}$.
We have
$(f'_*f^*)^0:
b'\xrightarrow{f'}\tilde{b}\xrightarrow[\cong]{\imath_{\tilde{b}}^{-1}}b\tensor \tilde{b}$
and for each pair $x$, $y\in\textit{Mod}_b$,
the underlying morphism of
$(f'_*f^*)^2_{x,y}
:f'_*f^*(x)\tensor'f'_*f^*(y)\to f'_*f^*(x\tensor y)$ in $\mathcal{V}$
is the unique morphism satisfying the following relation.
\begin{equation} \label{eq SixFunctors f'*f* tensor coherence}
  \vcenter{\hbox{
    \xymatrix@R=25pt@C=60pt{
      (x\tensor \tilde{b})\otimes (y\tensor \tilde{b})
      \ar@{->>}[d]_-{\textit{cq}_{x\tensor \tilde{b},y\tensor \tilde{b}}}
      \ar@{->>}[r]^-{\textit{cq}'_{x\tensor \tilde{b},y\tensor \tilde{b}}}
      &(x\tensor \tilde{b})\tensor' (y\tensor \tilde{b})
      \ar@{.>}[dd]^-{(f'_*f^*)^2_{x,y}}_-{\exists!}
      \\
      x\tensor \tilde{b}\tensor y\tensor \tilde{b}
      \ar[d]_-{\I_x\tensor s_{\tilde{b},y}\tensor\I_{\tilde{b}}}^-{\cong}
      &\text{ }
      \\
      x\tensor y\tensor \tilde{b}\tensor \tilde{b}
      \ar[r]^-{\I_{x\tensor y}\tensor m_{\tilde{b}/b}}
      &x\tensor y\tensor \tilde{b}
    }
  }}
\end{equation}
The construction
$\bigl(b\xrightarrow{f}\tilde{b}\xleftarrow{f'}b'\bigr)\mapsto f'_*f^*$
is functorial,
and we obtain a functor
\begin{equation} \label{eq SixFunctors MainThm leftadj}
  \vcenter{\hbox{
    \xymatrix@R=0pt@C=35pt{
      \textit{Comm}_{b\otimes b'}
      \ar[r]
      &\mathcal{V}\textit{-Funct}^{\textit{lax}}(\textit{Mod}_b,\textit{Mod}_{b'}\!)
      \\
      \bigl(b
      \xrightarrow[]{f}
      \tilde{b}
      \xleftarrow[]{f'}
      b'\bigr)
      \ar@{|->}[r]
      &f'_*f^*:\textit{Mod}_b\to\textit{Mod}_{b'}
    }
  }}
\end{equation}
from the category of commutative monoids in $\mathcal{V}$ over $b\otimes b'$
to the category of lax monoidal enriched functors
$\textit{Mod}_b\to \textit{Mod}_{b'}$.

\section{Main result} \label{sec MainResult}
The goal of this section is to prove Theorem~\ref{thm Intro MainThm} and Theorem~\ref{thm Intro MainThm2}.
Throughout this section
we fix commutative monoids
$b$, $b'$ in $\mathcal{V}$.
Let
$\mathcal{F}:\textit{Mod}_b\to \textit{Mod}_{b'}$
be a lax monoidal enriched functor.
We denote the lax monoidal coherence morphisms in $\textit{Mod}_{b'}$
of $\mathcal{F}$ as
\begin{equation*}
  \mathcal{F}^0:
  b'\to \mathcal{F}(b)
  ,
  \qquad
  \mathcal{F}^2_{x,y}:
  \mathcal{F}(x)\tensor'\mathcal{F}(y)
  \to \mathcal{F}(x\tensor y)
  ,
  \qquad
  x,y\in\textit{Mod}_b
  .
\end{equation*}
The following diagrams commute for all
$w\in \mathcal{V}$,
$x$, $y\in\textit{Mod}_b$
as the morphism
$\mathcal{F}^2_{x,y}$
is enriched natural in variables $x$, $y$.
\begin{equation} \label{eq MainResult Fxbyb Cnaturalxb}
  \vcenter{\hbox{
    \xymatrix@R=20pt@C=50pt{
      w\otimes (\mathcal{F}(x)\tensor' \mathcal{F}(y))
      \ar[d]_-{\I_w\otimes \mathcal{F}^2_{x,y}}
      \ar[r]^-{a'_{w,\mathcal{F}(x),\mathcal{F}(y)}}_-{\cong}
      &(w\otimes \mathcal{F}(x))\tensor' \mathcal{F}(y)
      \ar[d]^-{t^{\mathcal{F}}_{w,x}\tensor' \I_{\mathcal{F}(y)}}
      \\
      w\otimes \mathcal{F}(x\tensor y)
      \ar[d]_-{t^{\mathcal{F}}_{w,x\tensor y}}
      &\mathcal{F}(w\otimes x)\tensor' \mathcal{F}(y)
      \ar[d]^-{\mathcal{F}^2_{w\otimes x,y}}
      \\
      \mathcal{F}(w\otimes (x\tensor y))
      \ar[r]^-{\mathcal{F}(a_{w,x,y})}_-{\cong}
      &\mathcal{F}((w\otimes x)\tensor y)
    }
  }}
\end{equation}
\begin{equation} \label{eq MainResult Fxbyb Cnaturalyb}
  \vcenter{\hbox{
    \xymatrix@R=20pt@C=60pt{
      w\otimes (\mathcal{F}(x)\tensor' \mathcal{F}(y))
      \ar[d]_-{\I_w\otimes \mathcal{F}^2_{x,y}}
      \ar[r]^-{t^{\mathcal{F}(x)\tensor -}_{w,\mathcal{F}(y)}}_-{\cong}
      &\mathcal{F}(x)\tensor' (w\otimes \mathcal{F}(y))
      \ar[d]^-{\I_{\mathcal{F}(x)}\tensor' t^{\mathcal{F}}_{w,y}}
      \\
      w\otimes \mathcal{F}(x\tensor y)
      \ar[d]_-{t^{\mathcal{F}}_{w,x\tensor y}}
      &\mathcal{F}(x)\tensor' \mathcal{F}(w\otimes y)
      \ar[d]^-{\mathcal{F}^2_{x,w\otimes y}}
      \\
      \mathcal{F}(w\otimes (x\tensor y))
      \ar[r]^-{\mathcal{F}(t^{x\tensor -}_{w,y})}_-{\cong}
      &\mathcal{F}(x\tensor (w\otimes y))
    }
  }}
\end{equation}
We also have a commutative monoid
$\mathcal{F}(b)$ in $\textit{Mod}_{b'}$
whose unit, product morphisms are
\begin{equation*}
    u_{\mathcal{F}(b)/b'}:
    \xymatrix@C=15pt{
      b'
      \ar[r]^-{\mathcal{F}^0}
      &\mathcal{F}(b)
      ,
    }
    \qquad
    m_{\mathcal{F}(b)/b'}:
    \xymatrix{
      \mathcal{F}(b)\tensor' \mathcal{F}(b)
      \ar[r]^-{\mathcal{F}^2_{b,b}}
      &\mathcal{F}(b\tensor b)
      \ar[r]^-{\mathcal{F}(\jmath_b)}_-{\cong}
      &\mathcal{F}(b)
      .
    }
\end{equation*}

\begin{lemma} \label{lem MainResult MainThm rightadj}
  We have the functor
  of evaluating at $b\in\textit{Mod}_b$
  \begin{equation} \label{eq MainThm rightadj}
    \vcenter{\hbox{
      \xymatrix@R=0pt@C=23pt{
        \mathcal{V}\textit{-Funct}^{\textit{lax}}(\textit{Mod}_b,\textit{Mod}_{b'}\!)
        \ar[r]
        &\textit{Comm}_{b\otimes b'}
        \\
        \mathcal{F}:\textit{Mod}_b\to \textit{Mod}_{b'}
        \ar@{|->}[r]
        &\bigl(b
        \xrightarrow[]{f(\mathcal{F})}
        \mathcal{F}(b)
        \xleftarrow[]{f'(\mathcal{F})}
        b'\bigr)
      }
    }}
  \end{equation}
  which sends each lax monoidal enriched functor
  $\mathcal{F}:\textit{Mod}_b\to \textit{Mod}_{b'}$
  to the following object
  $(f(\mathcal{F}),f'(\mathcal{F})):b\otimes b'\to \mathcal{F}(b)$
  in $\textit{Comm}_{b\otimes b'}$.
  \begin{enumerate}
    \item 
    We have a commutative monoid $\mathcal{F}(b)$
    in $\mathcal{V}$ whose unit, product morphisms are
    \begin{equation*}
      \begin{aligned}
        u_{\mathcal{F}(b)}
        &:\!\!
        \xymatrix@C=27pt{
          c
          \ar[r]^-{u_{b'}}
          &b'
          \ar[r]^-{u_{\mathcal{F}(b)/b'}}
          &\mathcal{F}(b)
          ,
        }
        \\
        m_{\mathcal{F}(b)}
        &:\!\!
        \xymatrix@C=35pt{
          \mathcal{F}(b)\otimes \mathcal{F}(b)
          \ar@{->>}[r]^(0.45){\textit{cq}'_{\mathcal{F}(b),\mathcal{F}(b)}}
          &\mathcal{F}(b)\tensor' \mathcal{F}(b)
          \ar[r]^-{m_{\mathcal{F}(b)/b'}}
          &\mathcal{F}(b)
          .
        }
      \end{aligned}
    \end{equation*}
    
    \item 
    We have a morphism
    $f'(\mathcal{F})=u_{\mathcal{F}(b)/b'}=\mathcal{F}^0:b'\to \mathcal{F}(b)$
    in $\textit{Comm}(\mathcal{V})$.
    
    \item
    We have a morphism
    $f(\mathcal{F}):b\to \mathcal{F}(b)$
    in $\textit{Comm}(\mathcal{V})$
    whose underlying morphism in $\mathcal{V}$ is
    \begin{equation} \label{eq MainResult f(F)}
      \xymatrix@C=30pt{
        f(\mathcal{F}):
        b
        \ar[r]^-{\jmath_b^{-1}}_-{\cong}
        &b\otimes c
        \ar[r]^-{\I_b\otimes u_{\mathcal{F}(b)}}
        &b\otimes \mathcal{F}(b)
        \ar[r]^-{\rho_{\mathcal{F}(b)}}
        &\mathcal{F}(b)
        .
      }
    \end{equation}
    Here,
    $\rho_{\mathcal{F}(b)}
    =\mathcal{F}(m_b)\circ t^{\mathcal{F}}_{b,b}
    :b\otimes \mathcal{F}(b)
    \to \mathcal{F}(b)$
    is the left $b$-action
    defined in (\ref{eq EnEWThm rightadj}).
  \end{enumerate}
\end{lemma}
\begin{proof}
  We will show that the functor (\ref{eq MainThm rightadj})
  is well-defined on objects,
  and leave for the readers to check the rest of the details.
  From the isomorphism of categories 
  $\textit{Comm}_{b'}\cong \textit{Comm}(\textit{Mod}_{b'}\!)$
  described in (\ref{eq SixFunctors Commb=CommMb}),
  the object
  $\mathcal{F}(b)$ in $\textit{Comm}(\mathcal{V})$
  and the morphism
  $f'(\mathcal{F}):b'\to \mathcal{F}(b)$
  in $\textit{Comm}(\mathcal{V})$ are well-defined.
  Let us denote the morphism in (\ref{eq MainResult f(F)}) as
  $f:b\to \mathcal{F}(b)$.
  Using the left $b$-action relations of $\rho_{\mathcal{F}(b)}$,
  one can check that
  \begin{equation*}
    f\circ u_b=u_{\mathcal{F}(b)}
    :c\to \mathcal{F}(b)
    ,
    \qquad
    f\circ m_b=\rho_{\mathcal{F}(b)}\circ (\I_b\otimes f)
    :b\otimes b\to \mathcal{F}(b)
    .
  \end{equation*}
  To conclude that $f:b\to\mathcal{F}(b)$
  is a morphism  in $\textit{Comm}(\mathcal{V})$,
  it suffices to show that
  \begin{equation*}
    m_{\mathcal{F}(b)}\circ (f\otimes \I_b)=\rho_{\mathcal{F}(b)}
    :b\otimes \mathcal{F}(b)\to \mathcal{F}(b).      
  \end{equation*}
  This is shown in \textsection\!~\ref{subsec Appendix MainThm rightadj}.
  Thus the functor (\ref{eq MainThm rightadj}) is well-defined on objects.
\end{proof}

\begin{proposition} \label{prop MainResult lambdaF tensornat}
  Let $\mathcal{F}:\textit{Mod}_b\to \textit{Mod}_{b'}$
  be a lax monoidal enriched functor
  and denote
  $(f,f')=(f(\mathcal{F}),f'(\mathcal{F})):b\otimes b'\to \mathcal{F}(b)$
  as the object in $\textit{Comm}_{b\otimes b'}$
  described in Lemma~\ref{lem MainResult MainThm rightadj}.
  Then the enriched natural transformation
  $\lambda^{\mathcal{F}}:-\tensor \mathcal{F}(b)\Rightarrow\mathcal{F}$
  defined in Proposition~\ref{prop EnEWThm defining lambda}
  becomes a monoidal enriched natural transformation
  \begin{equation*}
    \lambda^{\mathcal{F}}:
    f'_*f^*\Rightarrow \mathcal{F}
  \end{equation*}
  between lax monoidal enriched functors
  from $\textit{Mod}_b$ to $\textit{Mod}_{b'}$.
\end{proposition}
\begin{proof}
  The $(b,b')$-bimodule
  $\mathcal{F}(b)$
  described in (\ref{eq EnEWThm rightadj})
  is the same as the $(b,b')$-bimodule
  $\mathcal{F}(b)$
  whose actions are induced by $f$, $f'$.
  Thus we have
  $-\tensor \mathcal{F}(b)=f'_*f^*$
  as enriched functors from $\textit{Mod}_b$ to $\textit{Mod}_{b'}$.
  Moreover, we have the relation
  \begin{equation*}
    \vcenter{\hbox{
      \xymatrix@R=20pt@C=30pt{
        b'
        \ar@/^1pc/[dr]^-{\mathcal{F}^0}
        \ar[d]_-{(f'_*f^*)^0}
        \\
        b\tensor \mathcal{F}(b)
        \ar[r]^-{\lambda^{\mathcal{F}}_{b}}
        &\mathcal{F}(b)
      }
    }}
    \qquad\qquad
    \begin{aligned}
      &
      \lambda^{\mathcal{F}}_{b}
      \circ
      (f'_*f^*)^0
      \\
      &=
      \imath_{\mathcal{F}(b)}
      \circ
      (\imath_{\mathcal{F}(b)})^{-1}
      \circ
      u_{\mathcal{F}(b)/b'}
      \\
      &=
      u_{\mathcal{F}(b)/b'}
      =
      \mathcal{F}^0:
      b'\to \mathcal{F}(b)
      .
    \end{aligned}
  \end{equation*}
  We are left to show the following relation for all
  $x$, $y\in\textit{Mod}_b$.
  \begin{equation*}
    \vcenter{\hbox{
      \xymatrix@R=20pt@C=45pt{
        (x\tensor \mathcal{F}(b))\tensor' (y\tensor \mathcal{F}(b))
        \ar[d]_-{(f'_*f^*)^2_{x,y}}
        \ar[r]^-{\lambda^{\mathcal{F}}_{x}\tensor'\lambda^{\mathcal{F}}_{y}}
        &\mathcal{F}(x)\tensor' \mathcal{F}(y)
        \ar[d]^-{\mathcal{F}^2_{x,y}}
        \\
        (x\tensor y)\tensor \mathcal{F}(b)
        \ar[r]^-{\lambda^{\mathcal{F}}_{x\tensor y}}
        &\mathcal{F}(x\tensor y)
      }
    }}
  \end{equation*}
  This is shown in \textsection\!~\ref{subsec Appendix lambdaF tensornat}.
\end{proof}

\begin{theorem}\label{thm Intro MainThm}
  Let $b$, $b'$ be commutative monoids in $\mathcal{V}$.
  \begin{enumerate}
    \item 
    We have a fully faithful left adjoint functor
    \begin{equation*}
      \vcenter{\hbox{
        \xymatrix@R=0pt{
          \textit{Comm}_{b\otimes b'}
          \ar[r]
          &\mathcal{V}\text{-}\textit{Funct}^{\textit{lax}}(\textit{Mod}_b,\textit{Mod}_{b'}\!)
          \\
          \bigl(b
          \xrightarrow[]{f}
          \tilde{b}
          \xleftarrow[]{f'}
          b'\bigr)
          \ar@{|->}[r]
          &f'_*f^*:
          \textit{Mod}_b\to \textit{Mod}_{b'}
        }
      }}
    \end{equation*}
    from the category of commutative monoids
    in $\mathcal{V}$ over $b\otimes b'$
    to the category of lax monoidal enriched functors
    $\textit{Mod}_b\to \textit{Mod}_{b'}$.

    \item
    The essential image of this functor
    is the coreflective full subcategory
    $\mathcal{V}\text{-}\textit{Funct}^{\textit{lax}}_{\textit{cocon}}(\textit{Mod}_b,\textit{Mod}_{b'}\!)$
    of cocontinuous lax monoidal enriched functors
    $\textit{Mod}_b\to \textit{Mod}_{b'}$,
    and we have an adjoint equivalence of categories
    \begin{equation*}
      \vcenter{\hbox{
        \xymatrix@R=0pt{
          \textit{Comm}_{b\otimes b'}
          \ar@<0.5ex>[r]^-{\simeq}
          &\text{ }
          \mathcal{V}\text{-}\textit{Funct}^{\textit{lax}}_{\textit{cocon}}(\textit{Mod}_b,\textit{Mod}_{b'}\!)
          .
          \ar@<0.5ex>[l]^-{\simeq}
        }
      }}
    \end{equation*}
  \end{enumerate}
\end{theorem}
\begin{proof}
  We first show statement 1.
  We claim that the functor
  $\textit{Comm}_{b\otimes b'}
  \to \mathcal{V}\textit{-Funct}^{\textit{lax}}(\textit{Mod}_b,\textit{Mod}_{b'}\!)$
  in (\ref{eq SixFunctors MainThm leftadj})
  is left adjoint to the functor
  $\mathcal{V}\textit{-Funct}^{\textit{lax}}(\textit{Mod}_b,\textit{Mod}_{b'}\!)
  \to \textit{Comm}_{b\otimes b'}$
  of evaluating at $b$, described in (\ref{eq MainThm rightadj}).
  \begin{equation}\label{eq MainResult MainThm adjunction}
    \vcenter{\hbox{
      \xymatrix@C=40pt{
        \textit{Comm}_{b\otimes b'}
        \ar@<0.7ex>[r]^-{(\ref{eq SixFunctors MainThm leftadj})}
        &\mathcal{V}\textit{-Funct}^{\textit{lax}}(\textit{Mod}_b,\textit{Mod}_{b'}\!)
        \ar@<0.7ex>[l]^-{(\ref{eq MainThm rightadj})}
      }
    }}
  \end{equation}
  \begin{itemize}
    \item 
    Let $(f,f'):b\otimes b'\to \tilde{b}$
    be an object in  $\textit{Comm}_{b\otimes b'}$.
    The object 
    $b\tensor \tilde{b}$ in $\mathcal{V}$
    has a unique structure of a commutative monoid in $\mathcal{V}$
    such that the isomorphism
    $\imath_{\tilde{b}}:b\tensor \tilde{b}\xrightarrow[]{\cong}\tilde{b}$
    in $\mathcal{V}$
    becomes an isomorphism
    in $\textit{Comm}(\mathcal{V})$.
    One can check that the functor
    $\mathcal{V}\textit{-Funct}^{\textit{lax}}(\textit{Mod}_b,\textit{Mod}_{b'}\!)
    \to \textit{Comm}_{b\otimes b'}$
    in (\ref{eq MainThm rightadj}) sends
    the lax monoidal enriched functor
    $f'_*f^*:\textit{Mod}_b\to \textit{Mod}_{b'}$
    to the object
    $(\imath^{-1}_{\tilde{b}}\circ f,\imath^{-1}_{\tilde{b}}\circ f')
    :b\otimes b'\to b\tensor \tilde{b}$
    in $\textit{Comm}_{b\otimes b'}$.
    The component of the unit of the adjunction (\ref{eq MainResult MainThm adjunction})
    at each object $(f,f'):b\otimes b'\to \tilde{b}$ in $\textit{Comm}_{b\otimes b'}$
    is the isomorphism
    $\imath^{-1}_{\tilde{b}}:
    \tilde{b}\xrightarrow[]{\cong} b\tensor \tilde{b}$
    in $\textit{Comm}_{b\otimes b'}$.

    \item
    The component of the counit of the adjunction (\ref{eq MainResult MainThm adjunction})
    at each lax monoidal enriched functor
    $\mathcal{F}:\textit{Mod}_b\to\textit{Mod}_{b'}$
    is the monoidal enriched natural transformation
    $\lambda^{\mathcal{F}}:f'(\mathcal{F})_*f(\mathcal{F})^*\Rightarrow\mathcal{F}$
    described in Proposition~\ref{prop MainResult lambdaF tensornat}.
  \end{itemize}
  In the proof of Theorem~\ref{thm Intro EnEW},
  we already showed that $\lambda^{\mathcal{F}}$ is natural in variable $\mathcal{F}$,
  and that the unit, counit satisfy the triangle identities
  (\ref{eq EnEWProof triangle identities}).
  Therefore we obtain the adjunction
  (\ref{eq MainResult MainThm adjunction}) as we claimed,
  whose unit is a natural isomorphism.
  Hence the left adjoint 
  $\textit{Comm}_{b\otimes b'}
  \to \mathcal{V}\textit{-Funct}^{\textit{lax}}(\textit{Mod}_b,\textit{Mod}_{b'}\!)$
  is fully faithful.
  
  Next we show statement 2.
  By Proposition~\ref{prop EnEWThm defining lambda},
  $\lambda^{\mathcal{F}}$ is a monoidal enriched natural isomorphism
  if and only if
  $\mathcal{F}$ is cocontinuous as an enriched functor.
  We conclude that 
  the essential image of the left adjoint
  is the coreflective full subcategory
  $\mathcal{V}\textit{-Funct}^{\textit{lax}}_{\textit{cocon}}(\textit{Mod}_b,\textit{Mod}_{b'}\!)$
  of cocontinous lax monoidal enriched functors
  $\textit{Mod}_b\to \textit{Mod}_{b'}$.
\end{proof}

\begin{theorem} \label{thm Intro MainThm2}
  Let $b$, $b'$ be commutative monoids in $\mathcal{V}$.
  We have an equivalence of categories
  \begin{equation*}
    \xymatrix@R=0pt{
      \Hom_{\textit{Comm}(\mathcal{V})}(b,b')
      \ar[r]^-{\simeq}
      &\mathcal{V}\text{-}\textit{Funct}^{\textit{strong}}_{\textit{cocon}}(\textit{Mod}_b,\textit{Mod}_{b'}\!)
      \\
      f:b\to b'
      \text{ }
      \ar@{|->}[r]^-{\simeq}
      &\text{ }
      f^*:\textit{Mod}_b\to \textit{Mod}_{b'}
    }
  \end{equation*}
  from the discrete category of morphisms
  $b\to b'$ of commutative monoids in $\mathcal{V}$
  to the category of cocontinuous strong monoidal enriched functors
  $\textit{Mod}_b\to \textit{Mod}_{b'}$.
\end{theorem}
\begin{proof}
  Let us denote
  $\mathcal{V}\textit{-Funct}^{\textit{strong}}_{\textit{cocon}}(\textit{Mod}_b,\textit{Mod}_{b'}\!)$
  as the category of cocontinuous strong monoidal enriched functors
  $\textit{Mod}_b\to \textit{Mod}_{b'}$,
  which is a full subcategory of
  $\mathcal{V}\textit{-Funct}^{\textit{lax}}_{\textit{cocon}}(\textit{Mod}_b,\textit{Mod}_{b'}\!)$.
  Moreover, regard the set
  $\Hom_{\textit{Comm}(\mathcal{V})}(b,b')$
  as a discrete category.
  We claim that the functor
  \begin{equation}\label{eq MainResult MainCor functor}
    \xymatrix@R=0pt{
      \Hom_{\textit{Comm}(\mathcal{V})}(b,b')
      \ar[r]
      &\mathcal{V}\textit{-Funct}^{\textit{strong}}_{\textit{cocon}}(\textit{Mod}_b,\textit{Mod}_{b'}\!)
      \\
      f:b\to b'
      \ar@{|->}[r]
      &f^*:\textit{Mod}_b\to \textit{Mod}_{b'}
    }
  \end{equation}
  is an equivalence of categories.
  Let $f_1$, $f_2:b\to b'$ be morphisms in $\textit{Comm}(\mathcal{V})$.
  By Theorem~\ref{thm Intro MainThm}, we have
  \begin{equation*}
    \begin{aligned}
      \Hom_{\mathcal{V}\textit{-Funct}^{\textit{strong}}_{\textit{cocon}}(\textit{Mod}_b,\textit{Mod}_{b'}\!)}(f_1^*,f_2^*)
      &=
      \Hom_{\mathcal{V}\textit{-Funct}^{\textit{lax}}_{\textit{cocon}}(\textit{Mod}_b,\textit{Mod}_{b'}\!)}(f_1^*,f_2^*)
      \\
      &\cong
      \Hom_{\textit{Comm}_{b\otimes b'}}
      \bigl(
        b\otimes b'
        \xrightarrow[]{(f_1,I_{b'})}b'
        ,\text{ }
        b\otimes b'
        \xrightarrow[]{(f_2,I_{b'})}
        b'
      \bigr)
      \\
      &\cong
      \begin{cases}
        \{*\}
        &\text{ if }
        f_1=f_2
        ,
        \\
        \varnothing
        &\text{ otherwise}
        .
      \end{cases}
    \end{aligned}
  \end{equation*}
  This shows that the functor (\ref{eq MainResult MainCor functor}) is fully faithful.
  
  We are left to show that the functor (\ref{eq MainResult MainCor functor}) is essentially surjective.
  Let $\mathcal{F}:\textit{Mod}_b\to \textit{Mod}_{b'}$
  be a cocontinuous strong monoidal enriched functor.
  By Lemma~\ref{lem MainResult MainThm rightadj},
  we have an object
  $(f,f')=(f(\mathcal{F}),f'(\mathcal{F}))
  :b\otimes b'\to \mathcal{F}(b)$
  in $\textit{Comm}_{b\otimes b'}$,
  where $f'=\mathcal{F}^0:b'\xrightarrow{\cong}\mathcal{F}(b)$
  is an isomorphism in $\textit{Comm}(\mathcal{V})$.
  Moreover, 
  $\lambda^{\mathcal{F}}:f'_*f^*\cong \mathcal{F}$
  is a monoidal enriched natural isomorphism
  by Proposition~\ref{prop EnEWThm defining lambda} and Proposition~\ref{prop MainResult lambdaF tensornat}.
  Let us denote
  \begin{equation*}
    \tilde{f}:= f'^{-1}\circ f:b\to b'.
  \end{equation*}
  We have an object
  $(\tilde{f},\I_{b'}):b\otimes b'\to b'$
  in $\textit{Comm}_{b\otimes b'}$
  and an isomorphism
  $f':b'\xrightarrow[]{\cong} \mathcal{F}(b)$
  in $\textit{Comm}_{b\otimes b'}$.
  The image of this isomorphism
  under the functor (\ref{eq SixFunctors MainThm leftadj})
  is a monoidal enriched natural isomorphism
  $\I_{b'*}\tilde{f}^*\cong f'_*f^*$.
  Thus we have
  $\tilde{f}^*=\I_{b'*}\tilde{f}^*
  \cong f'_*f^*\cong \mathcal{F}$.
  This shows that the functor (\ref{eq MainResult MainCor functor}) is essentially surjective.
\end{proof}

\section{An approach via Day convolution}
\label{sec DayConvolution}
Throughout this section, we fix commutative monoids $b$, $b'$ in $\mathcal{V}$
and explain an alternative way of obtaining
the equivalence of categories
\begin{equation} \label{eq DayConvolution MainThmequivalence}
  \xymatrix@C=30pt{
    \textit{Comm}_{b\otimes b'}
    \ar[r]^-{\simeq}
    &\mathcal{V}\textit{-Funct}^{\textit{lax}}(\textit{Mod}_b,\textit{Mod}_{b'}\!) 
  }
\end{equation}
described in
Theorem~\ref{thm Intro MainThm}.
Consider the following equivalence of categories.
\begin{equation} \label{eq DayConvolution equivalence}
  \vcenter{\hbox{
    \xymatrix@R=0pt{
      \textit{Mod}_{b\otimes b'}
      \ar[r]^-{\cong}
      &\text{ }\!\!_b\textit{Mod}_{b'}
      \ar[r]^-{\simeq}
      &\mathcal{V}\textit{-Funct}_{\textit{cocon}}(\textit{Mod}_b,\textit{Mod}_{b'}\!)
      \\
      x\ar@{|->}[r]^-{\cong}
      &x\ar@{|->}[r]^-{\simeq}
      &-\tensor x:\textit{Mod}_b\to \textit{Mod}_{b'}
    }
  }}
\end{equation}
We have the symmetric monoidal category
of right $b\otimes b'$-modules
\begin{equation*}
  (\textit{Mod}_{b\otimes b'},\widetilde{\tensor},b\otimes b').
\end{equation*}
Meanwhile, the category
$\mathcal{V}\textit{-Funct}_{\textit{cocon}}(\textit{Mod}_b,\textit{Mod}_{b'}\!)$
also has a symmetric monoidal category structure,
whose monoidal product is given by Day convolution.
The goal of this section is to show that
the equivalence of categories (\ref{eq DayConvolution equivalence})
is compatible with symmetric monoidal category structures on both sides.
As a consequence, we obtain an equivalence
\begin{equation} \label{eq DayConvolution laxadjunction}
  \xymatrix@C=20pt{
    \textit{Comm}(\textit{Mod}_{b\otimes b'})
    \ar[r]^-{\simeq}
    &\mathcal{V}\textit{-Funct}^{\textit{lax}}_{\textit{cocon}}(\textit{Mod}_b,\textit{Mod}_{b'}\!)
  }
\end{equation}
by considering the categories of commutative monoids on both sides
of (\ref{eq DayConvolution equivalence}).
Together with the isomorphism of categories
$\textit{Comm}_{b\otimes b'}\cong
\textit{Comm}(\textit{Mod}_{b\otimes b'})$
explained in (\ref{eq SixFunctors Commb=CommMb}),
we deduce the desired equivalence of categories (\ref{eq DayConvolution MainThmequivalence}).

Recall that the Day convolution of enriched functors
$\mathcal{F}$, $\mathcal{G}:\textit{Mod}_b\to \textit{Mod}_{b'}$
is the left Kan extension
$\mathcal{H}:\textit{Mod}_b\to \textit{Mod}_{b'}$
of the enriched functor
$\tensor'\circ (\mathcal{F}\times \mathcal{G})
:\textit{Mod}_b\times \textit{Mod}_b\to \textit{Mod}_{b'}$
along $\tensor:\textit{Mod}_b\times \textit{Mod}_b\to \textit{Mod}_b$.
We denote the universal enriched natural transformation as
$\theta_{-,-}:\mathcal{F}(-)\tensor' \mathcal{G}(-)
\Rightarrow \mathcal{H}(-\tensor -)$.
Day convolutions and their properties were first studied by B. Day in
\cite{Day1970Thesis}, \cite{Day1970}.
\begin{equation*}
  \vcenter{\hbox{
    \xymatrix@C=15pt{
      \textit{Mod}_b\times \textit{Mod}_b
      \ar[rr]^-{\tensor}
      \ar[d]_-{\mathcal{F}\times \mathcal{G}}
      &\text{ }
      &\textit{Mod}_b
      \ar@{.>}[d]^-{\mathcal{H}}
      \\
      \textit{Mod}_{b'}\times \textit{Mod}_{b'}
      \ar[rr]_-{\tensor'}
      &\text{ }
      \xtwocell[u]{}<>{<0>{\theta\text{ }\text{ }}}
      &\textit{Mod}_{b'}
    }
  }}
\end{equation*}

Let $x$, $y$ be $(b,b')$-bimodules.
We are going to describe the Day convolution of enriched functors
$-\tensor x$, $-\tensor y:\textit{Mod}_b\to \textit{Mod}_{b'}$.
Consider $x$, $y$ as right $b\otimes b'$-modules
and denote their right $b\otimes b'$-actions as
$\widetilde{\gamma}_x:x\otimes b\otimes b'\to x$,
$\widetilde{\gamma}_y:y\otimes b\otimes b'\to y$.
The right $b\otimes b'$-module $x\widetilde{\tensor}y$
is given by the following coequalizer.
\begin{equation*}
  \xymatrix@C=35pt{
    x\otimes b\otimes b'\otimes y
    \ar@<0.5ex>[r]^-{\widetilde{\gamma}_x\otimes \I_y}
    \ar@<-0.5ex>[r]_-{\I_x\otimes \mathring{\widetilde{\gamma}}_y}
    &x\otimes y
    \ar@{->>}[r]^-{\widetilde{\textit{cq}}_{x,y}}
    &x\wtensor y
  }
\end{equation*}
We can also see $x\wtensor y$ as a $(b,b')$-bimodule
whose actions $\rho_{x\widetilde{\tensor} y}$,
$\gamma'_{x\widetilde{\tensor}y}$ of $b$, $b'$ are unique morphisms in $\mathcal{V}$
satisfying the following relations.
\begin{equation*}
  \xymatrix@C=50pt{
    b\otimes x\otimes y
    \ar@{->>}[d]_-{\I_b\otimes \widetilde{\textit{cq}}_{x,y}}
    \ar[r]^-{\rho_x\otimes I_y}
    &x\otimes y
    \ar@{->>}[d]|-{\widetilde{\textit{cq}}_{x,y}}
    &x\otimes y\otimes b'
    \ar@{->>}[d]^-{\widetilde{\textit{cq}}_{x,y}\otimes \I_{b'}}
    \ar[l]_-{\I_x\otimes \gamma'_y}
    \\
    b\otimes (x\wtensor y)
    \ar@{.>}[r]^-{\exists!\text{ }\rho_{x\widetilde{\tensor} y}}
    &x\wtensor y
    &(x \wtensor y)\otimes b'
    \ar@{.>}[l]_-{\exists!\text{ }\gamma'_{x\widetilde{\tensor}y}}
  }
\end{equation*}
Let us denote the coprojections
of the coproduct $b\otimes b'$ of $b$, $b'$ in $\textit{Comm}(\mathcal{V})$
as $p:b\to b\otimes b'$ and $p':b'\to b\otimes b'$.
By Lemma~\ref{lem SixFunctors f*f! Cadjuctions},
Proposition~\ref{prop SixFunctors laxtensorCadj}
we have lax monoidal enriched functors
$p_*:\textit{Mod}_{b\otimes b'}\to \textit{Mod}_b$
and
$p'_*:\textit{Mod}_{b\otimes b'}\to \textit{Mod}_{b'}$
whose coherence morphisms $(p_*)^2$, $(p'_*)^2$
satisfy the following relations in $\mathcal{V}$.
\begin{equation} \label{eq DayConvolution p*p'* formula}
  \vcenter{\hbox{
    \xymatrix@C=40pt{
      x\otimes y
      \ar@{->>}[r]^-{\textit{cq}_{x,y}}
      \ar@{->>}@/_1pc/[dr]_-{\widetilde{\textit{cq}}_{x,y}}
      &x\tensor y
      \ar[d]^-{(p_*)^2_{x,y}}
      &x\otimes y
      \ar@{->>}[r]^-{\textit{cq}'_{x,y}}
      \ar@{->>}@/_1pc/[dr]_-{\widetilde{\textit{cq}}_{x,y}}
      &x\tensor' y
      \ar[d]^-{(p'_*)^2_{x,y}}
      \\
      \text{ }
      &x\wtensor y
      &\text{ }
      &x\wtensor y
    }
  }}
\end{equation}

\begin{proposition} \label{prop DayConvolution monoidalproduct}
  Let $x$, $y$ be $(b,b')$-bimodules.
  The Day convolution of enriched functors
  $-\tensor x$, $-\tensor y:\textit{Mod}_b\to \textit{Mod}_{b'}$
  is the enriched functor
  $-\tensor (x \wtensor y):\textit{Mod}_b\to \textit{Mod}_{b'}$.
  The component of the universal enriched natural transformation
  \begin{equation*}
    \theta_{-,-}:
    (-\tensor x)\tensor' (-\tensor y)
    \Rightarrow
    (-\tensor -)\tensor (x\wtensor y)
  \end{equation*}
  at $w$, $z\in\textit{Mod}_b$
  is the unique morphism
  $\theta_{w,z}$ in $\textit{Mod}_{b'}$
  whose underlying morphism in $\mathcal{V}$
  satisfies the following relation.
  \begin{equation}\label{eq DayConvolution universaltheta}
    \vcenter{\hbox{
      \xymatrix@R=23pt@C=60pt{
        (w\tensor x)\otimes (z\tensor y)
        \ar@{->>}[d]_-{\textit{cq}_{w\tensor x,z\tensor y}}
        \ar@{->>}[r]^-{\textit{cq}'_{w\tensor x,z\tensor y}}
        &(w\tensor x)\tensor' (z\tensor y)
        \ar@{.>}[dd]^-{\theta_{w,z}}_-{\exists!}
        \\
        w\tensor x\tensor z\tensor y
        \ar[d]_-{\I_w\tensor s_{x,z}\tensor \I_y}^-{\cong}
        &\text{ }
        \\
        w\tensor z\tensor x\tensor y
        \ar[r]^-{\I_{w\tensor x}\tensor (p_*)^2_{x,y}}
        &w\tensor z\tensor (x\wtensor y)
      }
    }}
  \end{equation}
\end{proposition}
\begin{proof}
  We leave for the readers to check that
  the enriched natural transformation
  $\theta$ is well-defined.
  We verify the universal property of $\theta$ as follows.
  Suppose we are given an enriched functor
  $\mathcal{F}:\textit{Mod}_b\to \textit{Mod}_{b'}$
  and an enriched natural transformation
  $\alpha_{-,-}:
  (-\tensor x)\tensor' (-\tensor y)\Rightarrow \mathcal{F}(-\tensor -)$.
  We need to show that there exists a unique
  enriched natural transformation
  $\xi:-\tensor (x\wtensor y)\Rightarrow\mathcal{F}$
  which satisfies the relation
  $\alpha_{-,-}=\xi_{-\text{ }\!\tensor\text{ }\! -}\circ \theta_{-,-}$.
  \begin{equation*}
    \vcenter{\hbox{
      \xymatrix@R=20pt@C=5pt{
        \textit{Mod}_b\times \textit{Mod}_b
        \ar[rr]^-{\tensor}
        \ar[d]_-{(-\tensor x)\times (-\tensor y)}
        &\text{ }
        &\textit{Mod}_b
        \ar[d]^-{\mathcal{F}}
        \\
        \textit{Mod}_{b'}\times \textit{Mod}_{b'}
        \ar[rr]_-{\tensor'}
        &\text{ }
        \xtwocell[u]{}<>{<0>{\alpha\text{ }\text{ }}}
        &\textit{Mod}_{b'}
      }
    }}
    \quad=\quad
    \vcenter{\hbox{
      \xymatrix@R=20pt@C=15pt{
        \textit{Mod}_b\times \textit{Mod}_b
        \ar[rr]^-{\tensor}
        \ar[d]_-{(-\tensor x)\times (-\tensor y)}
        &\text{ }
        &\textit{Mod}_b
        \ar@<-2ex>@/_1.5pc/[d]|-{-\tensor (x\widetilde{\tensor} y)}
        \ar@<1ex>@/^1.5pc/[d]^-{\mathcal{F}}
        \\
        \textit{Mod}_{b'}\times \textit{Mod}_{b'}
        \ar[rr]_-{\tensor'}
        &\text{ }
        \xtwocell[u]{}<>{<-4>{\theta\text{ }\text{ }}}
        &\textit{Mod}_{b'}
        \xtwocell[u]{}<>{<0>{\exists!\text{ }\xi}}
      }
    }}
  \end{equation*}
  Consider the morphism
  $\phi:x\tensor' y\to \mathcal{F}(b)$
  in $\textit{Mod}_{b'}$ defined as
  \begin{equation} \label{eq DayConvolution monoidalproduct phidef}
    \xymatrix@C=25pt{
      \phi:
      x\tensor' y
      \ar[r]^-{\imath^{-1}_x\tensor'\imath^{-1}_y}_-{\cong}
      &(b\tensor x)\tensor' (b\tensor y)
      \ar[r]^-{\alpha_{b,b}}
      &\mathcal{F}(b\tensor b)
      \ar[r]^-{\mathcal{F}(\jmath_b)}_-{\cong}
      &\mathcal{F}(b)
      .
    }
  \end{equation}
  We have a well-defined morphism
  $\psi:x\wtensor y\to \mathcal{F}(b)$ of $(b,b')$-bimodules
  whose underlying morphism in $\mathcal{V}$ satisfies the following relation
  \begin{equation}\label{eq DayConvolution monoidalproduct psidef}
    \vcenter{\hbox{
      \xymatrix@R=20pt@C=40pt{
        x\otimes y
        \ar@{->>}[r]^-{\textit{cq}'_{x,y}}
        \ar@{->>}[d]_-{\widetilde{\textit{cq}}_{x,y}}
        &x\tensor' y
        \ar[d]^-{\phi}
        \\
        x\wtensor y
        \ar@{.>}[r]^-{\exists!\text{ }\psi}
        &\mathcal{F}(b)
      }
    }}
  \end{equation}
  and the composition
  $\xi:
  \xymatrix@C=20pt{
    -\tensor (x\wtensor y)
    \ar@2{->}[r]^-{-\tensor \psi}
    &-\tensor \mathcal{F}(b)
    \ar@2{->}[r]^-{\lambda^{\mathcal{F}}}
    &\mathcal{F}
  }$
  is the unique enriched natural transformation
  that we are looking for:
  see \textsection\!~\ref{subsec Appendix DayConvolution monoidalproduct} for the details.
  We conclude that the enriched functor
  $-\tensor (x\wtensor y):\textit{Mod}_b\to \textit{Mod}_{b'}$
  together with the enriched natural transformation $\theta$
  is the Day convolution of enriched functors
  $-\tensor x$, $-\tensor y:\textit{Mod}_b\to \textit{Mod}_{b'}$.
\end{proof}

As we have
Proposition~\ref{prop DayConvolution monoidalproduct}
and the adjoint equivalence of categories in
Corollary~\ref{cor EnEWthm bimodule},
we explicitly see that the Day convolution of
two cocontinuous enriched functors
from $\textit{Mod}_b$ to $\textit{Mod}_{b'}$
exists, and is also cocontinuous.
Therefore the category
\begin{equation*}
  \mathcal{V}\textit{-Funct}_{\textit{cocon}}(\textit{Mod}_b,\textit{Mod}_{b'}\!)
\end{equation*}
has a symmetric monoidal category structure,
whose unit object is
$\textit{Mod}_b(b,-)\otimes b':\textit{Mod}_b\to \textit{Mod}_{b'}$
and whose monoidal product is given by Day convolution.
The category of commutative monoids
in this symmetric monoidal category
is equivalent to the category
\begin{equation*}
  \mathcal{V}\textit{-Funct}^{\textit{lax}}_{\textit{cocon}}(\textit{Mod}_b,\textit{Mod}_{b'}\!)
\end{equation*}
of cocontinuous lax monoidal enriched functors
$\textit{Mod}_b\to \textit{Mod}_{b'}$.
The analogous results for enriched functors
from a small monoidal enriched category
to $\mathcal{V}$ 
are explained in 
\cite[\textsection\!~3]{Day1970Thesis},
\cite[\textsection\!~6.2]{Loregian2021}.

We claim that the equivalence of categories
$\textit{Mod}_{b\otimes b}
\xrightarrow[]{\simeq}
\mathcal{V}\textit{-Funct}_{\textit{cocon}}(\textit{Mod}_b,\textit{Mod}_{b'}\!)$
in (\ref{eq DayConvolution equivalence})
is compatible with symmetric monoidal category structures on both sides.
The unit object in $\textit{Mod}_{b\otimes b'}$
corresponds to the $(b,b')$-bimodule
$b\otimes b'$
whose actions of $b$, $b'$ are
$m_b\otimes \I_{b'}$ and $\I_b\otimes m_{b'}$.
For each $z\in\textit{Mod}_b$,
we have an isomorphism 
$e_{z,b,b'}:z\tensor (b\tensor b')\xrightarrow[]{\cong} z\otimes b'$
in $\textit{Mod}_{b'}$
which is the unique morphism satisfying the relation below.
\begin{equation*}
  \xymatrix@C=40pt{
    z\otimes (b\otimes b')
    \ar@{->>}[d]_-{\textit{cq}_{z,b\otimes b'}}
    \ar[r]^-{a'_{z,b,b'}}_-{\cong}
    &(z\otimes b)\otimes b'
    \ar[d]^-{\gamma_z\otimes \I_{b'}}
    \\
    z\tensor (b\otimes b')
    \ar@{.>}[r]^-{\exists!\text{ }e_{z,b,b'}}_-{\cong}
    &z\otimes b'
  }
\end{equation*}
We also have an isomorphism
$\bar{\gamma}_z:z\xrightarrow[]{\cong}\textit{Mod}_b(b,z)$
in $\mathcal{V}$ which corresponds to the morphism
$\gamma_z:z\otimes b\to z$ in $\textit{Mod}_b$.
One can check that we have an enriched natural isomorphism
\begin{equation} \label{eq DayConvolution units}
  \xymatrix@C=20pt{
    -\tensor (b\otimes b')
    \ar@2{->}[r]^-{\cong}
    &\textit{Mod}_b(b,-)\otimes b'
  }
\end{equation}
between enriched functors from $\textit{Mod}_b$ to $\textit{Mod}_{b'}$,
whose component at $z\in\textit{Mod}_b$ is
$(\bar{\gamma}_z\otimes \I_{b'})\circ e_{z,b,b'}$.
As we have Proposition~\ref{prop DayConvolution monoidalproduct} 
and the enriched natural isomorphism
(\ref{eq DayConvolution units}),
we conclude that the equivalence of categories
$\textit{Mod}_{b\otimes b}
\xrightarrow[]{\simeq}
\mathcal{V}\textit{-Funct}_{\textit{cocon}}(\textit{Mod}_b,\textit{Mod}_{b'}\!)$
is compatible with symmetric monoidal category structures on both sides.
Thus we have
\begin{equation*}
  \xymatrix{
    \textit{Comm}_{b\otimes b'}
    \ar[r]^-{(\ref{eq SixFunctors Commb=CommMb})}_-{\cong}
    &\textit{Comm}(\textit{Mod}_{b\otimes b'})
    \ar[r]^-{(\ref{eq DayConvolution laxadjunction})}_-{\simeq}
    &\mathcal{V}\textit{-Funct}^{\textit{lax}}_{\textit{cocon}}(\textit{Mod}_b,\textit{Mod}_{b'}\!)
  }
\end{equation*}
which is another way of obtaining
the equivalence of categories
(\ref{eq DayConvolution MainThmequivalence})
as explained in the beginning of this section.

\appendix
\section{Appendix}
\label{sec Appendix}
In this section, we exhibit details that we omitted.

\subsection{Proposition~\ref{prop SixFunctors MbtensorCcat}}
\label{subsec Appendix MbtensorCcat}
Let $b$ be a commutative monoid in $\mathcal{V}$
and let $x$, $y$, $z\in\textit{Mod}_b$.
The isomorphism $a_{z,x,y}$ in $\textit{Mod}_b$
is well-defined
as we can see from the diagram below,
where we leave for the readers to verify the relation
$\mathring{\gamma}_{x\tensor y}=
(\mathring{\gamma}_x\tensor\I_y)
\circ a_{b,x,y}
:b\otimes (x\tensor y)\to x\tensor y$.
\begin{equation*}
  \xymatrix@R=10pt@C=45pt{
    z\otimes (b\otimes (x\tensor y))
    \ar[d]_-{\I_z\otimes a_{b,x,y}}^-{\cong}
    \ar@<0.5ex>[rr]^-{(\gamma_z\otimes \I_{x\tensor y})\circ a_{z,b,x\tensor y}}
    \ar@<-0.5ex>[rr]_-{\I_z\otimes \mathring{\gamma}_{x\tensor y}}
    &\text{ }
    &z\otimes (x\tensor y)
    \ar[dd]^-{a_{z,x,y}}_-{\cong}
    \ar@{->>}[r]^-{\textit{cq}_{z,x\tensor y}}
    &z\tensor (x\tensor y)
    \ar@{.>}[dd]^-{\exists!\text{ }a_{z,x,y}}_-{\cong}
    \\
    z\otimes ((b\otimes x)\tensor y)
    \ar[d]_-{a_{z,b\otimes x,y}}^-{\cong}
    &\text{ }
    &\text{ }
    \\
    (z\otimes (b\otimes x))\tensor y
    \ar@<0.5ex>[rr]^-{\bigl((\gamma_z\otimes \I_{x})\circ a_{z,b,x}\bigr)\tensor\I_{y}}
    \ar@<-0.5ex>[rr]_-{(\I_z\otimes \mathring{\gamma}_{x})\tensor\I_{y}}
    &\text{ }
    &(z\otimes x)\tensor y
    \ar@{->>}[r]^-{\textit{cq}_{z,x}\tensor \I_{y}}
    &(z\tensor x)\tensor y
  }
\end{equation*}
The following diagram shows that
the morphism $s_{x,y}:x\tensor y\to y\tensor x$
in $\mathcal{V}$ is well-defined.
\begin{equation*}
  \xymatrix@R=12pt@C=40pt{
    (x\otimes b)\otimes y
    \ar[dd]_-{\gamma_x\otimes \I_y}
    \ar@{=}[r]
    &(x\otimes b)\otimes y
    \ar[d]^-{s_{x,b}\otimes \I_y}_-{\cong}
    \ar@{=}[r]
    &(x\otimes b)\otimes y
    \ar[d]^-{a_{x,b,y}^{-1}}_-{\cong}
    \ar@{=}[r]
    &(x\otimes b)\otimes y
    \ar[d]^-{a_{x,b,y}^{-1}}_-{\cong}
    \\
    \text{ }
    &(b\otimes x)\otimes y
    \ar@/_0.5pc/[dl]_(0.5){\mathring{\gamma}_x\otimes \I_y}
    \ar[d]^-{s_{b\otimes x,y}}_-{\cong}
    &x\otimes (b\otimes y)
    \ar[d]^-{s_{x,b\otimes y}}_-{\cong}
    \ar@{=}[r]
    &x\otimes (b\otimes y)
    \ar[d]^-{\I_x\otimes \mathring{\gamma}_y}
    \\
    x\otimes y
    \ar[dd]_-{s_{x,y}}^-{\cong}
    &y\otimes (b\otimes x)
    \ar[d]^-{a_{y,b,x}}_-{\cong}
    \ar@/_0.5pc/[ddl]_(0.35){\I_y\otimes \mathring{\gamma}_x}
    &(b\otimes y)\otimes x
    \ar@/^0.5pc/[dl]^(0.5){s_{b,y}\otimes \I_x}_-{\cong}
    \ar@/_0.5pc/[ddr]|-{\mathring{\gamma}_y\otimes \I_x}
    &x\otimes y
    \ar[dd]^-{s_{x,y}}_-{\cong}
    \\
    \text{ }
    &(y\otimes b)\otimes x
    \ar[d]^-{\gamma_y\otimes \I_x}
    &\text{ }
    &\text{ }
    \\
    y\otimes x
    \ar@{->>}[d]_-{\textit{cq}_{y,x}}
    &y\otimes x
    \ar@{->>}[d]^-{\textit{cq}_{y,x}}
    \ar@{=}[rr]
    &\text{ }
    &y\otimes x
    \ar@{->>}[d]^-{\textit{cq}_{y,x}}
    \\
    y\tensor x
    \ar@{=}[r]
    &y\tensor x
    \ar@{=}[rr]
    &\text{ }
    \ar@{=}[r]
    &y\tensor x
  }
\end{equation*}
Since $s_{x,y}:x\otimes y\xrightarrow[]{\cong}y\otimes x$ is an isomorphism in $\mathcal{V}$ with its inverse $s_{y,x}:y\otimes x\xrightarrow{\cong}y\otimes x$,
we obtain that
$s_{x,y}:x\tensor y\to y\tensor x$
is also an isomorphism in $\mathcal{V}$ with its inverse $s_{y,x}:y\tensor x\xrightarrow{\cong}x\tensor y$.
We see that
$s_{x,y}:x\tensor y\xrightarrow{\cong}y\tensor x$
is an isomorphism in $\textit{Mod}_b$,
by right-cancelling the epimorphism
$\textit{cq}_{x,y}\tensor \I_b$
in the diagram below.
\begin{equation*}
  \xymatrix@R=12pt{
    (x\otimes y)\otimes b
    \ar@{->>}[d]_-{\textit{cq}_{x,y}\otimes \I_b}
    \ar@{=}[r]
    &(x\otimes y)\otimes b
    \ar[dd]^-{\gamma_{x\otimes y}}
    \ar@{=}[r]
    &(x\otimes y)\otimes b
    \ar[d]^-{a^{-1}_{x,y,b}}_-{\cong}
    \ar@{=}[r]
    &(x\otimes y)\otimes b
    \ar[d]^-{s_{x,y}\otimes \I_b}_-{\cong}
    \ar@{=}[r]
    &(x\otimes y)\otimes b
    \ar@{->>}[d]^-{\textit{cq}_{x,y}\otimes \I_b}
    \\
    (x\tensor  y)\otimes b
    \ar[dd]_-{\gamma_{x\tensor y}}
    &\text{ }
    &x\otimes (y\otimes b)
    \ar@/^0.5pc/[dl]|-{\I_x\otimes \gamma_y}
    \ar[d]^-{s_{x,y\otimes b}}_-{\cong}
    &(y\otimes x)\otimes b
    \ar[d]^-{a^{-1}_{y,x,b}}_-{\cong}
    \ar@<1.5ex>@/^2pc/[ddd]|-{\gamma_{y\otimes x}}
    \ar@/^1pc/@{->>}[ddr]|(0.5){\textit{cq}_{y,x}\!\otimes \I_b}
    &(x\tensor  y)\otimes b
    \ar[dd]^-{s_{x,y}\otimes \I_b}_-{\cong}
    \\
    \text{ }
    &x\otimes y
    \ar@{->>}@/^0.5pc/[dl]|-{\textit{cq}_{x,y}}
    \ar[d]^-{s_{x,y}}_-{\cong}
    &(y\otimes b)\otimes x
    \ar@/^0.5pc/[dl]|-{\gamma_y\otimes \I_x}
    \ar[d]^-{a^{-1}_{y,b,x}}_-{\cong}
    &y\otimes (x\otimes b)
    \ar@/^0.5pc/[dl]^(0.5){\I_y\otimes s_{x,b}}_-{\cong}
    \ar[dd]^-{\I_y\otimes \gamma_x}
    &\text{ }
    \\
    x\tensor  y
    \ar[dd]_-{s_{x,y}}^-{\cong}
    &y\otimes x
    \ar@{->>}[dd]^-{\textit{cq}_{y,x}}
    &y\otimes (b\otimes x)
    \ar[d]^-{\I_y\otimes \mathring{\gamma}_x}
    &\text{ }
    &(y\tensor  x)\otimes b
    \ar[dd]^-{\gamma_{y\tensor x}}
    \\
    \text{ }
    &\text{ }
    &y\otimes x
    \ar@{->>}[d]^-{\textit{cq}_{y,x}}
    \ar@{=}[r]
    &y\otimes x
    \ar@{->>}[d]^-{\textit{cq}_{y,x}}
    &\text{ }
    \\
    y\tensor x
    \ar@{=}[r]
    &y\tensor x
    \ar@{=}[r]
    &y\tensor x
    \ar@{=}[r]
    &y\tensor x
    \ar@{=}[r]
    &y\tensor x
  }
\end{equation*}
This shows that the isomorphism
$s_{x,y}$ in $\textit{Mod}_b$ is well-defined.
Next we verify the symmetry, unital axioms
of coherence isomorphisms of $\textit{Mod}_b$.
We already showed that
$s_{y,x}=s_{x,y}^{-1}$.
We obtain the unital axiom by right-cancelling the epimorphisms
$\I_x\otimes \textit{cq}_{b,y}$
and
$\textit{cq}_{x,b\tensor y}$
in the diagram below.
\begin{equation*}
  \xymatrix@R=10pt@C=15pt{
    x\otimes (b\otimes y)
    \ar@{->>}[d]_-{\I_x\otimes \textit{cq}_{b,y}}
    \ar@{=}[rr]
    &\text{ }
    &x\otimes (b\otimes y)
    \ar[d]^-{a_{x,b,y}}_-{\cong}
    \ar@{=}[r]
    &x\otimes (b\otimes y)
    \ar[ddd]_-{\I_x\otimes \mathring{\gamma}_{y}}
    \ar@{=}[r]
    &x\otimes (b\otimes y)
    \ar@{->>}[d]^-{\I_x\otimes \textit{cq}_{b,y}}
    \\
    x\otimes (b\tensor y)
    \ar@{->>}[d]_-{\textit{cq}_{x,b\tensor y}}
    \ar@/^1pc/[dr]^-{a_{x,b,y}}_-{\cong}
    &\text{ }
    &(x\otimes b)\otimes y
    \ar@/_1pc/@{->>}[dl]|-{\textit{cq}_{x\otimes b,y}}
    \ar[dd]^-{\gamma_x\otimes \I_{y}}
    &\text{ }
    &x\otimes (b\tensor y)
    \ar@/_1pc/[ddl]_(0.4){\I_x\otimes \imath_{y}}^(0.55){\cong}
    \ar@{->>}[d]^-{\textit{cq}_{x,b\tensor y}}
    \\
    x\tensor (b\tensor y)
    \ar[d]_-{a_{x,b,y}}^-{\cong}
    &(x\otimes b)\tensor y
    \ar@{->>}@/^0.5pc/[dl]|-{\textit{cq}_{x,b}\tensor\I_{y}}
    \ar[dd]^-{\gamma_{x}\tensor\I_{y}}
    &\text{ }
    &\text{ }
    &x\tensor (b\tensor y)
    \ar[dd]^-{\I_{x}\tensor\imath_{y}}_-{\cong}
    \\
    (x\tensor b)\tensor y
    \ar[d]_-{\jmath_{x}\tensor\I_{y}}^-{\cong}
    &\text{ }
    &x\otimes y
    \ar@{->>}[d]^-{\textit{cq}_{x,y}}
    &x\otimes y
    \ar@{->>}[d]^-{\textit{cq}_{x,y}}
    &\text{ }
    \\
    x\tensor y
    \ar@{=}[r]
    &x\tensor y
    \ar@{=}[r]
    &x\tensor y
    \ar@{=}[r]
    &x\tensor y
    \ar@{=}[r]
    &x\tensor y
  }
\end{equation*}
\subsection{Proposition~\ref{prop SixFunctors laxtensorCadj}}
\label{subsec Appendix laxtensorCadj}
Let $f:b\to b'$ be a morphism in $\textit{Comm}(\mathcal{V})$.
Let $w\in \mathcal{V}$ and $x$, $y$, $z\in\textit{Mod}_{b'}$.
We show that $(f_*)^2_{x,y}$ is enriched natural in variable $x\in\textit{Mod}_{b'}$.
It suffices to show that
$(f_*)^2_{w\otimes x,y}\circ a_{w,x,y}
=a'_{w,x,y}\circ (\I_w\otimes (f_*)^2_{x,y})
:w\otimes (x\tensor y)\to (w\otimes x)\tensor' y$
which we obtain by
by right-cancelling the epimorphism
$\I_w\otimes \textit{cq}_{x,y}$ in the diagram below.
\begin{equation*}
  \xymatrix@R=10pt@C=50pt{
    w\otimes (x\otimes y)
    \ar@{->>}[d]_-{\I_w\otimes \textit{cq}_{x,y}}
    \ar@{=}[r]
    &w\otimes (x\otimes y)
    \ar[d]^-{a_{w,x,y}}_-{\cong}
    \ar@/^0.5pc/@{->>}[ddl]^(0.6){\I_w\otimes \textit{cq}'_{x,y}}
    \ar@{=}[r]
    &w\otimes (x\otimes y)
    \ar@{->>}[d]^-{\I_w\otimes \textit{cq}_{x,y}}
    \\
    w\otimes (x\tensor y)
    \ar[d]_-{\I_w\otimes (f_*)^2_{x,y}}
    &(w\otimes x)\otimes y
    \ar@{->>}[dd]^-{\textit{cq}'_{w\otimes x,y}}
    \ar@/^0.5pc/@{->>}[dr]^(0.5){\textit{cq}_{w\otimes x,y}}
    &w\otimes (x\tensor y)
    \ar[d]^-{a_{w,x,y}}_-{\cong}
    \\
    w\otimes (x\tensor' y)
    \ar[d]_-{a'_{w,x,y}}^-{\cong}
    &\text{ }
    &(w\otimes x)\tensor y
    \ar[d]^-{(f_*)^2_{w\otimes x,y}}
    \\
    (w\otimes x)\tensor' y
    \ar@{=}[r]
    &(w\otimes x)\tensor' y
    \ar@{=}[r]
    &(w\otimes x)\tensor' y
  }
\end{equation*}
Next we show that 
$(f_*)^2$
is compatible with coherence isomorphisms
$a$ of $\textit{Mod}_b$ and $a'$ of $\textit{Mod}_{b'}$.
We obtain this by right-cancelling the epimorphisms
$\I_z\otimes \textit{cq}_{x,y}$
and
$\textit{cq}_{z,x\tensor y}$
in the following continued diagrams in $\mathcal{V}$.
\begin{equation*}
  \vcenter{\hbox{
    \xymatrix@R=10pt@C=25pt{
      z\otimes (x\otimes y)
      \ar@{->>}[d]_-{\I_z\otimes \textit{cq}_{x,y}}
      \ar@{=}[r]
      &z\otimes (x\otimes y)
      \ar@{->>}[dd]^-{\I_z\otimes \textit{cq}'_{x,y}}
      \ar@{=}[rr]
      &\text{ }
      &z\otimes (x\otimes y)
      \ar[d]^-{a_{z,x,y}}_-{\cong}
      \\
      z\otimes (x\tensor y)
      \ar@{->>}[d]_-{\textit{cq}_{z,x\tensor y}}
      \ar@/^0.5pc/[dr]|-{\I_z\otimes (f_*)^2_{x,y}}
      &\text{ }
      &\text{ }
      &(z\otimes x)\otimes y
      \ar@/_1pc/@{->>}[ddl]_-{\textit{cq}'_{z\otimes x,y}}
      \ar@{->>}[dd]^-{\textit{cq}'_{z,x}\otimes \I_y}
      \\
      z\tensor (x\tensor y)
      \ar[d]_-{\I_z\tensor(f_*)^2_{x,y}}
      &z\otimes (x\tensor' y)
      \ar@{->>}@/^0.5pc/[dl]|-{\textit{cq}_{z,x\tensor' y}}
      \ar[dd]^-{\textit{cq}'_{z,x\tensor' y}}
      \ar@/^0.5pc/[dr]^-{a'_{z,x,y}}_-{\cong}
      &\text{ }
      &\text{ }
      \\
      z\tensor (x\tensor' y)
      \ar[d]_-{(f_*)^2_{z,x\tensor' y}}
      &\text{ }
      &(z\otimes x)\tensor' y
      \ar@{->>}[dd]^-{\textit{cq}'_{z,x}\tensor'\I_y}
      &(z\tensor' x)\otimes y
      \ar@{->>}[dd]^-{\textit{cq}'_{z\tensor' x,y}}
      \\
      z\tensor' (x\tensor' y)
      \ar[d]_-{a'_{z,x,y}}^-{\cong}
      \ar@{=}[r]
      &z\tensor' (x\tensor' y)
      \ar[d]^-{a'_{z,x,y}}_-{\cong}
      &\text{ }
      &\text{ }
      \\
      (z\tensor' x)\tensor' y
      \ar@{=}[r]
      &(z\tensor' x)\tensor' y
      \ar@{=}[r]
      &(z\tensor' x)\tensor' y
      \ar@{=}[r]
      &(z\tensor' x)\tensor' y
    }      
  }}
  \cdots
\end{equation*}
\begin{equation*}
  \cdots
  \vcenter{\hbox{
    \xymatrix@R=10pt@C=25pt{
      z\otimes (x\otimes y)
      \ar[d]^-{a_{z,x,y}}_-{\cong}
      \ar@{=}[rrr]
      &\text{ }
      &\text{ }
      &z\otimes (x\otimes y)
      \ar@{->>}[d]^-{\I_z\otimes \textit{cq}_{x,y}}
      \\
      (z\otimes x)\otimes y
      \ar@{->>}[dd]^-{\textit{cq}'_{z,x}\otimes \I_y}
      \ar@{=}[r]
      &(z\otimes x)\otimes y
      \ar@{->>}[d]^-{\textit{cq}_{z,x}\otimes \I_y}
      \ar@{->>}@/^0.5pc/[dr]^-{\textit{cq}_{z\otimes x,y}}
      &\text{ }
      &z\otimes (x\tensor y)
      \ar@{->>}[d]^-{\textit{cq}_{z,x\tensor y}}
      \ar@/_0.5pc/[dl]_-{a_{z,x,y}}^-{\cong}
      \\
      \text{ }
      &(z\tensor x)\otimes y
      \ar[d]^-{(f_*)^2_{z,x}\otimes \I_y}
      \ar@{->>}@/^0.5pc/[dr]|-{\textit{cq}_{z\tensor x,y}}
      &(z\otimes x)\tensor y
      \ar@{->>}[d]^-{\textit{cq}_{z,x}\tensor\I_y}
      &z\tensor (x\tensor y)
      \ar[d]^-{a_{z,x,y}}_-{\cong}
      \\
      (z\tensor' x)\otimes y
      \ar@{->>}[dd]^-{\textit{cq}'_{z\tensor' x,y}}
      \ar@{=}[r]
      &(z\tensor' x)\otimes y
      \ar@/_0.5pc/@{->>}[dr]_-{\textit{cq}_{z\tensor' x,y}}
      &(z\tensor x)\tensor y
      \ar[d]^-{(f_*)^2_{z,x}\tensor\I_y}
      \ar@{=}[r]
      &(z\tensor x)\tensor y
      \ar[d]^-{(f_*)^2_{z,x}\tensor\I_y}
      \\
      \text{ }
      &\text{ }
      &(z\tensor' x)\tensor y
      \ar@{=}[r]
      &(z\tensor' x)\tensor y
      \ar[d]^-{(f_*)^2_{z\tensor' x,y}}
      \\
      (z\tensor' x)\tensor' y
      \ar@{=}[rrr]
      &\text{ }
      &\text{ }
      &(z\tensor' x)\tensor' y
    }      
  }}
\end{equation*}
Now let $x$, $y\in\textit{Mod}_b$.
We show that the morphisms
$(\mathring{f}^*)^2$,
$(f^*)^2$
in $\mathcal{V}$
are inverse to each other.
The morphism
$(\mathring{f}^*)^2_{x,y}$
in $\mathcal{V}$ satisfies the relation
\begin{equation} \label{eq f*op2 property}
  \vcenter{\hbox{
    \xymatrix@R=10pt@C=40pt{
      x\otimes y\tensor b'
      \ar@{->>}[d]_-{\textit{cq}_{x,y}\tensor \I_{b'}}
      \ar[r]^-{\eta^*_x\tensor \I_{y\tensor b'}}
      &(x\tensor b')\otimes (y\tensor b')
      \ar@{->>}[d]^-{\textit{cq}'_{x\tensor b', y\tensor b'}}
      \\
      x\tensor y\tensor b'
      \ar[r]^-{(\mathring{f}^*)^2_{x,y}}
      &(x\tensor b')\tensor' (y\tensor b')
    }
  }}
\end{equation}
as we can see from the diagram below.
\begin{equation*}
  \xymatrix@R=12pt@C=30pt{
    x\otimes y\tensor b'
    \ar@{->>}[d]_-{\textit{cq}_{x,y}\tensor\I_{b'}}
    \ar@{=}[r]
    &x\otimes y\tensor b'
    \ar[dd]^-{(\eta^*_{x}\otimes\eta^*_{y})\tensor\I_{b'}}
    \ar@{=}[r]
    &x\otimes y\tensor b'
    \ar[d]^-{\eta^*_{x}\otimes \I_{y\tensor b'}}
    \\
    x\tensor y\tensor b'
    \ar[dd]^-{(\eta^*_{x}\tensor \eta^*_{y})\tensor\I_{b'}}
    \ar@{.>}@/_2.5pc/@<-4ex>[dddd]|(0.35){(\mathring{f}^*)^2_{x,y}}
    &\text{ }
    &(x\tensor b')\otimes (y\tensor b')
    \ar@{=}[ddd]
    \ar@/^0.8pc/[dl]|-{\I_{x\tensor b'}\otimes (\eta^*_{y}\tensor\I_{b'})}
    \\
    \text{ }
    &(x\tensor b')\otimes (y\tensor b')\tensor b'
    \ar@/^0.5pc/@{->>}[dl]|-{\textit{cq}_{x\tensor b',y\tensor b'}\tensor\I_{b'}}
    \ar@{=}[d]
    &\text{ }
    \\
    x\tensor b'\tensor y\tensor b'\tensor b'
    \ar[d]^-{(f_*)^2_{x\tensor b',y\tensor b'}\tensor\I_{b'}}
    &(x\tensor b')\otimes (y\tensor b')\tensor b'
    \ar@/^0.5pc/@{->>}[dl]|-{\textit{cq}'_{x\tensor b',y\tensor b'}\tensor\I_{b'}}
    \ar[d]^-{\varepsilon^*_{(x\tensor b')\otimes (y\tensor b')}}
    \ar@/^0.8pc/[dr]|-{\I_{x\tensor b'}\otimes \varepsilon^*_{y\tensor b'}}
    &\text{ }
    \\
    ((x\tensor b')\tensor' (y\tensor b'))\tensor b'
    \ar[d]^-{\varepsilon^*_{(x\tensor b')\tensor' (y\tensor b')}}
    &(x\tensor b')\otimes (y\tensor b')
    \ar@{->>}[d]^-{\textit{cq}'_{x\tensor b',y\tensor b'}}
    \ar@{=}[r]
    &(x\tensor b')\otimes (y\tensor b')
    \ar@{->>}[d]^-{\textit{cq}'_{x\tensor b',y\tensor b'}}
    \\
    (x\tensor b')\tensor' (y\tensor b')
    \ar@{=}[r]
    &(x\tensor b')\tensor' (y\tensor b')
    \ar@{=}[r]
    &(x\tensor b')\tensor' (y\tensor b')
  }
\end{equation*}
We also have
$(\mathring{f}^*)^2_{x,y}
=(f_*)^2_{x\tensor b',y\tensor b'}\circ (\eta^*_x\tensor \I_{y\tensor b'})
:x\tensor y\tensor b'\to (x\tensor b')\tensor' (y\tensor b')$
by right-cancelling the epimorphism
$\textit{cq}_{x,y}\tensor\I_{b'}$
in the diagram below.
We use (\ref{eq f*op2 property}).
\begin{equation*}
  \xymatrix@R=10pt@C=30pt{
    x\otimes y\tensor b'
    \ar@{->>}[d]_-{\textit{cq}_{x,y}\tensor\I_{b'}}
    \ar@{=}[r]
    &x\otimes y\tensor b'
    \ar[d]^-{\eta^*_{x}\otimes \I_{y\tensor b'}}
    \ar@{=}[r]
    &x\otimes y\tensor b'
    \ar@{->>}[d]^-{\textit{cq}_{x,y\tensor b'}=\text{ }\!\textit{cq}_{x,y}\tensor\I_{b'}}
    \\
    x\tensor y\tensor b'
    \ar[dd]_-{(\mathring{f}^*)^2_{x,y}}
    &(x\tensor b')\otimes (y\tensor b')
    \ar@{->>}[dd]^-{\textit{cq}'_{x\tensor b',y\tensor b'}}
    \ar@{->>}@/^0.5pc/[dr]|-{\textit{cq}_{x\tensor b',y\tensor b'}}
    &x\tensor y\tensor b'
    \ar[d]^-{\eta^*_{x}\tensor\I_{y\tensor b'}}
    \\
    \text{ }
    &\text{ }
    &x\tensor b'\tensor y\tensor b'
    \ar[d]^-{(f_*)^2_{x\tensor b',y\tensor b'}}
    \\
    (x\tensor b')\tensor' (y\tensor b')
    \ar@{=}[r]
    &(x\tensor b')\tensor' (y\tensor b')
    \ar@{=}[r]
    &(x\tensor b')\tensor' (y\tensor b')
  }
\end{equation*}
Now we show that
$(\mathring{f}^*)^2$,
$(f^*)^2$
are inverse to each other.
We obtain
$(f^*)^2\circ (\mathring{f}^*)^2=\I$
by right-cancelling the epimorphism
$\textit{cq}_{x,y}\tensor\I_{b'}$
in the diagram below.
We use (\ref{eq f*op2 property}).
\begin{equation*}
  \xymatrix@R=10pt@C=20pt{
    x\otimes y\tensor b'
    \ar@{->>}[d]_-{\textit{cq}_{x,y}\tensor\I_{b'}}
    \ar@{=}[r]
    &x\otimes y\tensor b'
    \ar[dd]_-{\eta^*_{x}\otimes \I_{y\tensor b'}}
    \ar@{=}[r]
    &x\otimes y\tensor b'
    \ar@{->>}[d]^-{\textit{cq}_{x,y\tensor b'}}
    \ar@{=}[r]
    &x\otimes y\tensor b'
    \ar@{->>}[d]^-{\textit{cq}_{x,y}\tensor\I_{b'}}
    \\
    x\tensor y\tensor b'
    \ar[dd]_-{(\mathring{f}^*)^2_{x,y}}
    &\text{ }
    &x\tensor y\tensor b'
    \ar[d]^-{\jmath_{x}^{-1}\tensor\I_{y\tensor b'}}_-{\cong}
    \ar@<-2.2ex>@/_2pc/[dd]|(0.3){\eta^*_{x}\tensor \I_{y\tensor b'}}
    \ar@{=}[r]
    &x\tensor y\tensor b'
    \ar@{=}[dddd]
    \\
    \text{ }
    &(x\tensor b')\otimes (y\tensor b')
    \ar@/^0.5pc/@{->>}[dl]|-{\textit{cq}'_{x\tensor b',y\tensor b'}}
    \ar@/_0.5pc/@{->>}[dr]|-{\textit{cq}_{x\tensor b',y\tensor b'}}
    &x\tensor b\tensor y\tensor b'
    \ar[d]^-{\I_x\tensor f\tensor\I_{y\tensor b'}}
    &\text{ }
    \\
    (x\tensor b)\tensor' (y\tensor b')
    \ar[dd]_-{(f^*)^2_{x,y}}
    &\text{ }
    &x\tensor b'\tensor y\tensor b'
    \ar[d]^-{\I_x\tensor s_{b',y}\tensor\I_{b'}}_-{\cong}
    &\text{ }
    \\
    \text{ }
    &\text{ }
    &x\tensor y\tensor b'\tensor b'
    \ar[d]^-{\I_{x\tensor y}\tensor m_{b'/b}}
    &\text{ }
    \\
    x\tensor y\tensor b'
    \ar@{=}[rr]
    &\text{ }
    &x\tensor y\tensor b'
    \ar@{=}[r]
    &x\tensor y\tensor b'
  }
\end{equation*}
We also obtain
$(\mathring{f}^*)^2\circ (f^*)^2=\I$
by right-cancelling
$\textit{cq}_{x,b'}\otimes \I_{y\tensor b'}$
and
$\textit{cq}'_{x\tensor b',y\tensor b'}$
in the following continued diagrams.
\begin{equation*}
  \vcenter{\hbox{
    \xymatrix@R=12pt@C=30pt{
      x\otimes b'\otimes y\tensor b'
      \ar@{->>}[d]_-{\textit{cq}_{x,b'}\otimes \I_{y\tensor b'}}
      \ar@{=}[r]
      &x\otimes b'\otimes y\tensor b'
      \ar@{->>}[d]^-{\I_x\otimes \textit{cq}_{b',y\tensor b'}}
      \ar@{=}[r]
      &x\otimes b'\otimes y\tensor b'
      \ar@{->>}[d]^-{\I_x\otimes \textit{cq}_{b',y\tensor b'}}
      \ar@{.>}@/^3.5pc/@<5ex>[ddd]|-{\I_x\otimes \mathring{\gamma}'_{y\tensor b'}}
      \\
      (x\tensor b')\otimes (y\tensor b')
      \ar@{->>}[d]_-{\textit{cq}'_{x\tensor b',y\tensor b'}}
      \ar@{->>}@/^0.5pc/[dr]|-{\textit{cq}_{x\tensor b',y\tensor b'}}
      &x\otimes b'\tensor y\tensor b'
      \ar@{->>}[d]^-{\textit{cq}_{x,b'\tensor y\tensor b'}}
      \ar@{=}[r]
      &x\otimes b'\tensor y\tensor b'
      \ar[d]^-{\I_x\otimes s_{b',y}\tensor\I_{b'}}_-{\cong}
      \\
      (x\tensor b')\tensor' (y\tensor b')
      \ar[dd]_-{(f^*)^2_{x,y}}
      &x\tensor b'\tensor y\tensor b'
      \ar[d]^-{\I_x\tensor s_{b',y}\tensor\I_{b'}}_-{\cong}
      &x\otimes y\tensor b'\tensor b'
      \ar[d]^-{\I_x\otimes \I_y\tensor m_{b'/b}}
      \ar@/^0.5pc/@{->>}[dl]|(0.45){\textit{cq}_{x,y\tensor b'\tensor b'}}
      \\
      \text{ }
      &x\tensor y\tensor b'\tensor b'
      \ar[d]^-{\I_{x\tensor y}\tensor m_{b'/b}}
      &x\otimes y\tensor b'
      \ar@/^0.5pc/@{->>}[dl]|-{\textit{cq}_{x,y\tensor b'}}
      \ar[d]^-{\eta^*_{x}\otimes \I_{y\tensor b'}}
      \\
      x\tensor y\tensor b'
      \ar[dd]_-{(\mathring{f}^*)^2_{x,y}}
      \ar@{=}[r]
      &x\tensor y\tensor b'
      \ar[d]^-{\eta^*_{x}\tensor\I_{y\tensor b'}}
      &(x\tensor b')\otimes (y\tensor b')
      \ar@/^0.5pc/@{->>}[dl]|-{\textit{cq}_{x\tensor b',y\tensor b'}}
      \ar@{->>}[dd]^-{\textit{cq}'_{x\tensor b',y\tensor b'}}
      \\
      \text{ }
      &(x\tensor b')\tensor (y\tensor b')
      \ar[d]^-{(f_*)^2_{x\tensor b',y\tensor b'}}
      &\text{ }
      \\
      (x\tensor b)\tensor' (y\tensor b')
      \ar@{=}[r]
      &(x\tensor b)\tensor' (y\tensor b')
      \ar@{=}[r]
      &(x\tensor b)\tensor' (y\tensor b')
    }
  }}
  \!\!\!\!
  \!\!\!\!
  \!\!\!\!
  \!\!\!\!
  \cdots
\end{equation*}
\begin{equation*}
  \cdots
  \vcenter{\hbox{
    \xymatrix@R=12pt@C=20pt{
      x\otimes b'\otimes y\tensor b'
      \ar[d]^-{\I_x\otimes \mathring{\gamma}'_{y\tensor b'}}
      \ar@{=}[r]
      &x\otimes b'\otimes y\tensor b'
      \ar[d]^-{\eta^*_{x}\otimes \I_{b'}\otimes \I_{y\tensor b'}}
      \ar@{=}[r]
      &x\otimes b'\otimes y\tensor b'
      \ar@{->>}[dd]^-{\textit{cq}_{x,b'}\otimes \I_{y\tensor b'}}
      \\
      x\otimes y\tensor b'
      \ar[d]^-{\eta^*_{x}\otimes \I_{y\tensor b'}}
      &(x\tensor b')\otimes b'\otimes (y\tensor b')
      \ar@/^1pc/[dl]|-{\I_{x\tensor b'}\otimes \mathring{\gamma}'_{y\tensor b'}}
      \ar@/_1pc/[dr]|-{\gamma'_{x\tensor b'}\otimes \I_{y\tensor b'}}
      &\text{ }
      \\
      (x\tensor b')\otimes (y\tensor b')
      \ar@{->>}[d]^-{\textit{cq}'_{x\tensor b',y\tensor b'}}
      &\text{ }
      &(x\tensor b')\otimes (y\tensor b')
      \ar@{->>}[d]^-{\textit{cq}'_{x\tensor b',y\tensor b'}}
      \\
      (x\tensor b')\tensor' (y\tensor b')
      \ar@{=}[rr]
      &\text{ }
      &(x\tensor b')\tensor' (y\tensor b')
    }
  }}
\end{equation*}
We leave for the readers to check the last equality
$\gamma'_{x\tensor b'}\circ (\eta^*_{x}\otimes \I_{b'})
=\textit{cq}_{x,b'}$
using the relation (\ref{eq gamma'opgamma'}).
This shows that
the morphisms
$(\mathring{f}^*)^2$,
$(f^*)^2$
in $\mathcal{V}$
are inverse to each other.

Lastly, we show that $f^*\dashv f_*$
satisfies the projection formula:
see (\ref{eq Preliminary projformuladef}).
Let $x\in\textit{Mod}_b$, $y\in\textit{Mod}_{b'}$
and recall the $(b,b')$-bimodule $b'$.
We can see $y$ as a $(b,b')$-bimodule
whose left $b$-action is the opposite action $\mathring{\gamma}_y$
of $f_*(y)\in\textit{Mod}_b$.
We can also see $b'\tensor'y$
as a $(b,b')$-bimodule in an analogous way,
and $\imath'_y:b'\tensor'y\xrightarrow{\cong}y$
becomes an isomorphism in $\text{ }\!\!_b\textit{Mod}_{b'}$.
We leave for the readers to check that the following isomorphism in $\mathcal{V}$
exists.
\begin{equation*}
  \xymatrix@R=14pt@C=60pt{
    x\otimes (b'\tensor' y)
    \ar@{->>}[d]_-{\textit{cq}_{x,b'\tensor'y}}
    \ar[r]^-{a'_{x,b',y}}_-{\cong}
    &(x\otimes b')\tensor' y
    \ar@{->>}[d]^-{\textit{cq}_{x,b'}\tensor' \I_{y'}}
    \\
    x\tensor (b'\tensor' y)
    \ar@{.>}[r]^-{\exists!\text{ }a_{x,b',y}}_-{\cong}
    &(x\tensor b')\tensor' y
  }
\end{equation*}
We have the following description
as morphisms in $\mathcal{V}$
\begin{equation*}
  \xymatrix@R=14pt@C=50pt{
    x\tensor f_*(y)
    \ar@{=}[d]
    \ar[r]^-{\eta^*_{x}\tensor\I_{f_*(y)}}
    &f_*f^*(x)\tensor f_*(y)
    \ar[r]^-{(f_*)^2_{f^*(x),y}}
    &f_*(f^*(x)\tensor' y)
    \ar@{=}[d]
    \\
    x\tensor y
    \ar[r]^-{\I_x\tensor (\imath'_y)^{-1}}_-{\cong}
    &x\tensor (b'\tensor' y)
    \ar[r]^-{a_{x,b',y}}_-{\cong}
    &(x\tensor b')\tensor' y
  }
\end{equation*}
which we obtain by right-cancelling the epimorphism
$\textit{cq}_{x,y}$
in the diagram below.
\begin{equation*}
  \xymatrix@R=16pt@C=25pt{
    x\otimes y
    \ar@{->>}[dd]_-{\textit{cq}_{x,y}}
    \ar@{=}[r]
    &x\otimes y
    \ar@{.>}@<-3ex>@/_3pc/[dddd]|-{\eta^*_{x}\otimes\I_y}
    \ar[dd]_-{\jmath^{-1}_x\otimes \I_y}^-{\cong}
    \ar@{=}[r]
    &x\otimes y
    \ar[d]^-{\I_x\otimes \imath^{-1}_y}_-{\cong}
    \ar@{=}[r]
    &x\otimes y
    \ar@{.>}@/^1.3pc/[dddl]|-{\I_x\otimes (\imath'_y)^{-1}}
    \ar@{->>}[dd]^-{\textit{cq}_{x,y}}
    \\
    \text{ }
    &\text{ }
    &x\otimes (c\otimes y)
    \ar@/^0.3pc/[dl]_(0.5){a_{x,c,y}}^-{\cong}
    \ar[d]^-{\I_x\otimes (u_{b'}\otimes \I_y)}
    &\text{ }
    \\
    x\tensor y
    \ar[ddd]_-{\eta^*_{x}\tensor\I_y}
    &(x\otimes c)\otimes y
    \ar[d]_-{(\I_x\otimes u_{b'})\otimes \I_y}
    &x\otimes (b'\otimes y)
    \ar@/^0.3pc/[dl]_(0.5){a_{x,b',y}}^-{\cong}
    \ar@{->>}[d]^-{\I_x\otimes \textit{cq}'_{b',y}}
    &x\tensor y
    \ar[dd]^-{\I_x\tensor(\imath'_y)^{-1}}_-{\cong}
    \\
    \text{ }
    &(x\otimes b')\otimes y
    \ar@{->>}[d]_-{\textit{cq}_{x,b'}\otimes \I_y}
    \ar@/^0.5pc/@{->>}[ddr]|-{\textit{cq}'_{x\otimes b',y}}
    &x\otimes (b'\tensor' y)
    \ar[dd]^-{a'_{x,b',y}}_-{\cong}
    \ar@/^0.5pc/@{->>}[dr]|-{\textit{cq}_{x,b'\tensor'y}}
    &\text{ }
    \\
    \text{ }
    &(x\tensor b')\otimes y
    \ar@/_0.5pc/@{->>}[dl]|-{\textit{cq}_{x\tensor b',y}}
    \ar@{->>}[dd]_-{\textit{cq}'_{x\tensor b',y}}
    &\text{ }
    &x\tensor (b'\tensor' y)
    \ar[dd]^-{a_{x,b',y}}_-{\cong}
    \\
    (x\tensor b')\tensor y
    \ar[d]_-{(f_*)^2_{f^*(x),y}}
    &\text{ }
    &(x\otimes b')\tensor' y
    \ar@{->>}[d]^-{\textit{cq}_{x,b'}\tensor' \I_y}
    &\text{ }
    \\
    (x\tensor b')\tensor' y
    \ar@{=}[r]
    &(x\tensor b')\tensor' y
    \ar@{=}[r]
    &(x\tensor b')\tensor' y
    \ar@{=}[r]
    &(x\tensor b')\tensor' y
  }
\end{equation*}
In particular, $f^*\dashv f_*$ satisfies the projection formula.

\subsection{Lemma~\ref{lem MainResult MainThm rightadj}}
\label{subsec Appendix MainThm rightadj}
Let $b$, $b'$ be commutative monoids in $\mathcal{V}$
and let $\mathcal{F}:\textit{Mod}_b\to \textit{Mod}_{b'}$
be a lax monoidal enriched functor.
Consider
$f=f(\mathcal{F}):b\to \mathcal{F}(b)$
and
$\rho_{\mathcal{F}(b)}:b\otimes \mathcal{F}(b)\to \mathcal{F}(b)$
described in Lemma~\ref{lem MainResult MainThm rightadj}.
We show the relation
$m_{\mathcal{F}(b)}\circ (f\otimes \I_b)=\rho_{\mathcal{F}(b)}
:b\otimes \mathcal{F}(b)\to \mathcal{F}(b)$
as follows.

\begin{equation*}
  \xymatrix@R=25pt@C=23pt{
    b\otimes \mathcal{F}(b)
    \ar[d]^-{\jmath^{-1}_b\otimes \I_{\mathcal{F}(b)}}_-{\cong}
    \ar@{.>}@/_2pc/@<-2ex>[ddddd]|(0.3){f\otimes \I_{\mathcal{F}(b)}}
    \ar@{=}[r]
    &b\otimes \mathcal{F}(b)
    \ar[d]^-{\I_b\otimes \imath^{-1}_{\mathcal{F}(b)}}_-{\cong}
    \ar@{=}[r]
    &b\otimes \mathcal{F}(b)
    \ar[ddd]^-{\I_b\otimes \imath_{\mathcal{F}(b)}^{-1}}_-{\cong}
    \ar@{=}[r]
    &b\otimes \mathcal{F}(b)
    \ar@{=}[dddd]
    \\
    b\otimes c\otimes \mathcal{F}(b)
    \ar[d]^-{(\I_b\otimes u_{b'})\otimes \I_{\mathcal{F}(b)}}
    \ar@{=}[r]
    &b\otimes c\otimes \mathcal{F}(b)
    \ar[d]^-{\I_b\otimes (u_{b'}\otimes \I_{\mathcal{F}(b)})}
    &\text{ }
    &\text{ }
    \\
    b\otimes b'\otimes \mathcal{F}(b)
    \ar[d]^-{(\I_b\otimes \mathcal{F}^0)\otimes \I_{\mathcal{F}(b)}}
    \ar@/^0.8pc/@{->>}[dr]|(0.65){\textit{cq}'_{b\otimes b',\mathcal{F}(b)}}
    \ar@{=}[r]
    &b\otimes b'\otimes \mathcal{F}(b)
    \ar@/^0.8pc/@{->>}[dr]|(0.65){\I_b\otimes \textit{cq}'_{b',\mathcal{F}(b)}}
    &\text{ }
    &\text{ }
    \\
    b\otimes \mathcal{F}(b)\otimes \mathcal{F}(b)
    \ar[d]^-{t^{\mathcal{F}}_{b,b}\otimes \I_{\mathcal{F}(b)}}
    \ar@/^0.8pc/@{->>}[dr]|(0.65){\textit{cq}'_{b\otimes \mathcal{F}(b),\mathcal{F}(b)}}
    &b\otimes b'\tensor' \mathcal{F}(b)
    \ar[d]^-{(\I_b\otimes \mathcal{F}^0)\tensor'\I_{\mathcal{F}(b)}}
    \ar@{=}[r]
    &b\otimes b'\tensor' \mathcal{F}(b)
    \ar[d]^-{\I_b\otimes (\mathcal{F}^0\tensor'\I_{\mathcal{F}(b)})}
    &\text{ }
    \\
    \mathcal{F}(b\otimes b)\otimes \mathcal{F}(b)
    \ar[d]^-{\mathcal{F}(m_b)\otimes \I_{\mathcal{F}(b)}}
    \ar@/^0.8pc/@{->>}[dr]|(0.65){\textit{cq}'_{\mathcal{F}(b\otimes b),\mathcal{F}(b)}}
    &b\otimes \mathcal{F}(b)\tensor' \mathcal{F}(b)
    \ar[d]^-{t^{\mathcal{F}}_{b,b}\tensor'\I_{\mathcal{F}(b)}}
    \ar@{=}[r]
    &b\otimes \mathcal{F}(b)\tensor' \mathcal{F}(b)
    \ar[d]^-{\I_b\otimes \mathcal{F}^2_{b,b}}
    &b\otimes \mathcal{F}(b)
    \ar@/^1pc/[dl]_(0.25){\I_b\otimes \mathcal{F}(\imath_{b})^{-1}}^{\cong}
    \ar[d]^-{t^{\mathcal{F}}_{b,b}}
    \ar@{.>}@/^2pc/@<1ex>[dddd]|-{\rho_{\mathcal{F}(b)}}
    \\
    \mathcal{F}(b)\otimes \mathcal{F}(b)
    \ar@{.>}@/_2pc/@<-1ex>[ddd]|-{m_{\mathcal{F}(b)}}
    \ar@{->>}[d]^-{\textit{cq}'_{\mathcal{F}(b),\mathcal{F}(b)}}
    &\mathcal{F}(b\otimes b)\tensor' \mathcal{F}(b)
    \ar@/^0.5pc/[dl]|-{\mathcal{F}(m_b)\tensor'\I_{\mathcal{F}(b)}}
    \ar[d]^-{\mathcal{F}^2_{b\otimes b,b}}
    &b\otimes \mathcal{F}(b\tensor b)
    \ar[d]^-{t^{\mathcal{F}}_{b,b\tensor b}}
    &\mathcal{F}(b\otimes b)
    \ar@/^1pc/[dl]^(0.48){\mathcal{F}(\I_b\otimes \imath_{b})^{-1}}_{\cong}
    \ar@{=}[dd]
    \\
    \mathcal{F}(b)\tensor' \mathcal{F}(b)
    \ar[d]^-{\mathcal{F}^2_{b,b}}
    &\mathcal{F}(b\otimes b\tensor b)
    \ar@/^0.5pc/[dl]|-{\mathcal{F}(m_b\tensor\I_{b})}
    \ar[d]^-{\mathcal{F}(\jmath_{b\otimes b})}_-{\cong}
    \ar@{=}[r]
    &\mathcal{F}(b\otimes b\tensor b)
    \ar[d]^-{\mathcal{F}(\I_b\otimes \imath_{b})}_-{\cong}
    &\text{ }
    \\
    \mathcal{F}(b\tensor b)
    \ar[d]^-{\mathcal{F}(\jmath_{b})}_-{\cong}
    &\mathcal{F}(b\otimes b)
    \ar[d]^-{\mathcal{F}(m_b)}
    \ar@{=}[r]
    &\mathcal{F}(b\otimes b)
    \ar[d]^-{\mathcal{F}(m_b)}
    \ar@{=}[r]
    &\mathcal{F}(b\tensor b)
    \ar[d]^-{\mathcal{F}(m_b)}
    \\
    \mathcal{F}(b)
    \ar@{=}[r]
    &\mathcal{F}(b)
    \ar@{=}[r]
    &\mathcal{F}(b)
    \ar@{=}[r]
    &\mathcal{F}(b)
  }
\end{equation*}
We omit obvious brackets when describing objects,
and keep the information of the omitted brackets when describing morphisms.
We used the relation (\ref{eq MainResult Fxbyb Cnaturalxb})
and we leave for the readers to check that
$\jmath_{b\tensor b}
\circ a_{b,b,b}
=\I_b\otimes \imath_b
:b\otimes (b\tensor b)
\xrightarrow{\cong}
b\otimes b$.

\subsection{Proposition~\ref{prop MainResult lambdaF tensornat}}
\label{subsec Appendix lambdaF tensornat}
Let $b$, $b'$ be commutative monoids in $\mathcal{V}$
and let 
$\mathcal{F}:\textit{Mod}_b\to\textit{Mod}_{b'}$
be a lax monoidal enriched functor.
For each pair $x$, $y\in\textit{Mod}_b$,
we show the following relation.
\begin{equation} \label{eq MainResult lambdaF tensor}
  \vcenter{\hbox{
    \xymatrix@R=16pt@C=40pt{
      (x\tensor \mathcal{F}(b))\tensor' (y\tensor \mathcal{F}(b))
      \ar[d]_-{(f'_*f^*)^2_{x,y}}
      \ar[r]^-{\lambda^{\mathcal{F}}_{x}\tensor'\lambda^{\mathcal{F}}_{y}}
      &\mathcal{F}(x)\tensor' \mathcal{F}(y)
      \ar[d]^-{\mathcal{F}^2_{x,y}}
      \\
      (x\tensor y)\tensor \mathcal{F}(b)
      \ar[r]^-{\lambda^{\mathcal{F}}_{x\tensor y}}
      &\mathcal{F}(x\tensor y)
    }
  }}
\end{equation}
It suffices to verify the relation (\ref{eq MainResult lambdaF tensor})
as morphisms in $\mathcal{V}$.
We
pre-compose epimorphisms
$\textit{cq}'_{x\otimes\mathcal{F}(b),y\otimes\mathcal{F}(b)}$,
$\textit{cq}_{x,\mathcal{F}(b)}\tensor'\textit{cq}_{y,\mathcal{F}(b)}$
and begin our calculation.

\begin{equation*}
  \hspace*{-0.5cm}
  \vcenter{\hbox{
    \xymatrix@R=13pt@C=8pt{
      x\otimes \mathcal{F}(b)\otimes y\otimes \mathcal{F}(b)
      \ar@{->>}[d]_-{\textit{cq}'_{x\otimes\mathcal{F}(b),y\otimes\mathcal{F}(b)}}
      \ar@{=}[r]
      &x\otimes \mathcal{F}(b)\otimes y\otimes \mathcal{F}(b)
      \ar@{->>}[d]^-{\textit{cq}'_{x\otimes\mathcal{F}(b),y\otimes\mathcal{F}(b)}}
      \ar@{=}[r]
      &x\otimes \mathcal{F}(b)\otimes y\otimes \mathcal{F}(b)
      \ar@{->>}[d]^-{\I_x\otimes \textit{cq}'_{\mathcal{F}(b),y\otimes \mathcal{F}(b)}}
      \\
      (x\otimes \mathcal{F}(b))\tensor' (y\otimes \mathcal{F}(b))
      \ar@{->>}[d]_-{\I_{x\otimes \mathcal{F}(b)}\tensor'\textit{cq}_{y,\mathcal{F}(b)}}
      \ar@{=}[r]
      &x\otimes \mathcal{F}(b)\tensor' (y\otimes \mathcal{F}(b))
      \ar[d]^-{\I_{x\otimes\mathcal{F}(b)}\tensor't^{\mathcal{F}}_{y,b}}
      \ar@{=}[r]
      &x\otimes \mathcal{F}(b)\tensor' (y\otimes \mathcal{F}(b))
      \ar[d]^-{\I_x\otimes (\I_{\mathcal{F}(b)}\tensor't^{\mathcal{F}}_{y,b})}
      \\
      (x\otimes \mathcal{F}(b))\tensor' (y\tensor \mathcal{F}(b))
      \ar@{->>}[d]_-{\textit{cq}_{x,\mathcal{F}(b)}\tensor'\I_{y\tensor\mathcal{F}(b)}}
      \ar@/^0.5pc/[dr]|-{\I_{x\otimes\mathcal{F}(b)}\tensor'\lambda^{\mathcal{F}}_{y}}
      &x\otimes \mathcal{F}(b)\tensor' \mathcal{F}(y\otimes b)
      \ar[d]^-{\I_{x\otimes\mathcal{F}(b)}\tensor'\mathcal{F}(\gamma_{y})}
      \ar@{=}[r]
      &x\otimes \mathcal{F}(b)\tensor' \mathcal{F}(y\otimes b)
      \ar[d]^-{\I_x\otimes (\I_{\mathcal{F}(b)}\tensor'\mathcal{F}(\gamma_{y}))}
      \\
      (x\tensor \mathcal{F}(b))\tensor' (y\tensor \mathcal{F}(b))
      \ar[d]_-{\I_{x\tensor\mathcal{F}(b)}\tensor'\lambda^{\mathcal{F}}_{y}}
      &x\otimes \mathcal{F}(b)\tensor' \mathcal{F}(y)
      \ar@{->>}@/^0.5pc/[dl]|-{\textit{cq}_{x,\mathcal{F}(b)}\tensor'\I_{\mathcal{F}(y)}}
      \ar[d]^-{t^{\mathcal{F}}_{x,b}\tensor'\I_{\mathcal{F}(y)}}
      \ar@{=}[r]
      &x\otimes \mathcal{F}(b)\tensor' \mathcal{F}(y)
      \ar[d]^-{\I_x\otimes \mathcal{F}^2_{b,y}}
      \\
      (x\tensor \mathcal{F}(b))\tensor' \mathcal{F}(y)
      \ar[d]_-{\lambda^{\mathcal{F}}_{x}\tensor'\I_{\mathcal{F}(y)}}
      &\mathcal{F}(x\otimes b)\tensor' \mathcal{F}(y)
      \ar@/^0.5pc/[dl]|-{\mathcal{F}(\gamma_{x})\tensor'\I_{\mathcal{F}(y)}}
      \ar[d]^-{\mathcal{F}^2_{x\otimes b,y}}
      &x\otimes \mathcal{F}(b\tensor y)
      \ar[d]^-{t^{\mathcal{F}}_{x,b\tensor y}}
      \\
      \mathcal{F}(x)\tensor' \mathcal{F}(y)
      \ar[dd]_-{\mathcal{F}^2_{x,y}}
      &\mathcal{F}(x\otimes b\tensor y)
      \ar[dd]^-{\mathcal{F}(\gamma_{x}\tensor \I_y)}
      \ar@{=}[r]
      &\mathcal{F}(x\otimes b\tensor y)
      \ar[d]^-{\mathcal{F}(\I_x\otimes \imath_{y})}_-{\cong}
      \\
      \text{ }
      &\text{ }
      &\mathcal{F}(x\otimes y)
      \ar[d]^-{\mathcal{F}(\textit{cq}_{x,y})}
      \\
      \mathcal{F}(x\tensor y)
      \ar@{=}[r]
      &\mathcal{F}(x\tensor y)
      \ar@{=}[r]
      &\mathcal{F}(x\tensor y)
    }
  }}
  \!\!\!\!
  \!\!\!\!
  \cdots
\end{equation*}
We omit obvious brackets when describing objects,
and keep the information of the omitted brackets when describing morphisms.
We maintain the position of brackets in each of the vertical compositions.
We used the defining relations (\ref{eq EnEWThm defining lambda}) of
$\lambda^{\mathcal{F}}_{x}$,
$\lambda^{\mathcal{F}}_{y}$
as well as (\ref{eq MainResult Fxbyb Cnaturalxb}).
Readers should check that
$(\gamma_{x}\tensor\I_y)
\circ a_{x,b,y}
=\textit{cq}_{x,y}
\circ (\I_x\otimes \imath_{y})$.
We continue as follows.
\begin{equation*}
  \cdots
  \!\!\!\!
  \!\!\!\!
  \vcenter{\hbox{
    \xymatrix@R=10pt@C=12pt{
      x\otimes \mathcal{F}(b)\otimes y\otimes \mathcal{F}(b)
      \ar@{->>}[ddd]^-{\I_x\otimes \textit{cq}'_{\mathcal{F}(b),y\otimes \mathcal{F}(b)}}
      \ar@{=}[rr]
      &\text{ }
      &x\otimes \mathcal{F}(b)\otimes y\otimes \mathcal{F}(b)
      \ar[d]^-{\I_x\otimes s_{\mathcal{F}(b),y}\otimes \I_{\mathcal{F}(b)}}_-{\cong}
      \\
      \text{ }
      &\text{ }
      &x\otimes y\otimes \mathcal{F}(b)\otimes \mathcal{F}(b)
      \ar@{->>}[d]^-{\I_x\otimes (\I_y\otimes\textit{cq}'_{\mathcal{F}(b),\mathcal{F}(b)})}
      \\
      \text{ }
      &\text{ }
      &x\otimes y\otimes \mathcal{F}(b)\tensor' \mathcal{F}(b)
      \ar@/_1pc/[dll]_(0.4){\I_x\otimes t^{\mathcal{F}(b)\tensor'-}_{y,\mathcal{F}(b)}}^(0.4){\cong}
      \ar[d]^-{\I_x\otimes (\I_y\otimes \mathcal{F}^2_{b,b})}
      \\
      x\otimes \mathcal{F}(b)\tensor' (y\otimes \mathcal{F}(b))
      \ar[d]^-{\I_x\otimes (\I_{\mathcal{F}(b)}\tensor't^{\mathcal{F}}_{y,b})}
      &\text{ }
      &x\otimes y\otimes \mathcal{F}(b\tensor b)
      \ar[d]^-{\I_x\otimes t^{\mathcal{F}}_{y,b\tensor b}}
      \\
      x\otimes \mathcal{F}(b)\tensor' \mathcal{F}(y\otimes b)
      \ar[d]^-{\I_x\otimes (\I_{\mathcal{F}(b)}\tensor'\mathcal{F}(\gamma_{y}))}
      \ar@/^1pc/[dr]^(0.7){\I_x\otimes \mathcal{F}^2_{b,y\otimes b}}
      &\text{ }
      &x\otimes \mathcal{F}(y\otimes b\tensor b)
      \ar@/_1pc/[dl]_-{\I_x\otimes \mathcal{F}(t^{b\tensor-}_{y,b})}^-{\cong}
      \ar[dd]^-{\I_x\otimes \mathcal{F}(\I_y\otimes \jmath_{b})}_-{\cong}
      \\
      x\otimes \mathcal{F}(b)\tensor' \mathcal{F}(y)
      \ar[d]^-{\I_x\otimes \mathcal{F}^2_{b,y}}
      &x\otimes \mathcal{F}(b\tensor (y\otimes b))
      \ar[d]^-{\I_x\otimes \mathcal{F}(\I_b\tensor\gamma_{y})}
      \ar@/^1pc/[dr]^(0.6){\I_x\otimes \mathcal{F}(\imath_{y\otimes b})}_-{\cong}
      &\text{ }
      \\
      x\otimes \mathcal{F}(b\tensor y)
      \ar[d]^-{t^{\mathcal{F}}_{x,b\tensor y}}
      \ar@{=}[r]
      &x\otimes \mathcal{F}(b\tensor y)
      \ar[d]^-{\I_x\otimes \mathcal{F}(\imath_{y})}_-{\cong}
      &x\otimes \mathcal{F}(y\otimes b)
      \ar[d]^-{\I_x\otimes \mathcal{F}(\gamma_{y})}
      \\
      \mathcal{F}(x\otimes b\tensor y)
      \ar[d]^-{\mathcal{F}(\I_x\otimes \imath_{y})}_-{\cong}
      &x\otimes \mathcal{F}(y)
      \ar[d]^-{t^{\mathcal{F}}_{x,y}}
      \ar@{=}[r]
      &x\otimes \mathcal{F}(y)
      \ar[d]^-{t^{\mathcal{F}}_{x,y}}
      \\
      \mathcal{F}(x\otimes y)
      \ar[d]^-{\mathcal{F}(\textit{cq}_{x,y})}
      \ar@{=}[r]
      &\mathcal{F}(x\otimes y)
      \ar[d]^-{\mathcal{F}(\textit{cq}_{x,y})}
      \ar@{=}[r]
      &\mathcal{F}(x\otimes y)
      \ar[d]^-{\mathcal{F}(\textit{cq}_{x,y})}
      \\
      \mathcal{F}(x\tensor y)
      \ar@{=}[r]
      &\mathcal{F}(x\tensor y)
      \ar@{=}[r]
      &\mathcal{F}(x\tensor y)
    }      
  }}
  \!\!\!\!
  \!\!\!\!
  \cdots
\end{equation*}
We used the defining relation 
(\ref{eq SixFunctors zbtensorb- strength})
of the tensorial strength
$t^{\mathcal{F}(b)\tensor'-}$
associated to 
the enriched functor
$\mathcal{F}(b)\tensor'-:
\textit{Mod}_{b'}\to \textit{Mod}_{b'}$,
as well as (\ref{eq MainResult Fxbyb Cnaturalyb}).
Readers should check that
$\imath_{y\otimes b}\circ t^{b\tensor-}_{y,b}
=\I_y\otimes \jmath_{b}
:y\otimes (b\tensor b)
\xrightarrow{\cong}y\otimes b$.
We further continue as follows.
\begin{equation*}
  \cdots
  \!\!\!\!
  \vcenter{\hbox{
    \xymatrix@R=14pt@C=15pt{
      x\otimes \mathcal{F}(b)\otimes y\otimes \mathcal{F}(b)
      \ar[d]^-{\I_x\otimes s_{\mathcal{F}(b),y}\otimes \I_{\mathcal{F}(b)}}_-{\cong}
      \ar@{=}[r]
      &x\otimes \mathcal{F}(b)\otimes y\otimes \mathcal{F}(b)
      \ar[d]^-{\I_x\otimes s_{\mathcal{F}(b),y}\otimes \I_{\mathcal{F}(b)}}_-{\cong}
      \ar@{=}[r]
      &x\otimes \mathcal{F}(b)\otimes y\otimes \mathcal{F}(b)
      \ar[d]^-{\I_x\otimes s_{\mathcal{F}(b),y}\otimes \I_{\mathcal{F}(b)}}_-{\cong}
      \\
      x\otimes y\otimes \mathcal{F}(b)\otimes \mathcal{F}(b)
      \ar@{->>}[d]^-{\I_x\otimes (\I_y\otimes\textit{cq}'_{\mathcal{F}(b),\mathcal{F}(b)})}
      \ar@{=}[r]
      &x\otimes y\otimes \mathcal{F}(b)\otimes \mathcal{F}(b)
      \ar[ddd]^-{\I_x\otimes (\I_y\otimes m_{\mathcal{F}(b)})}
      \ar@{=}[r]
      &x\otimes y\otimes \mathcal{F}(b)\otimes \mathcal{F}(b)
      \ar[ddd]^-{\I_{x\otimes y}\otimes m_{\mathcal{F}(b)}}
      \\
      x\otimes y\otimes \mathcal{F}(b)\tensor' \mathcal{F}(b)
      \ar[d]^-{\I_x\otimes (\I_y\otimes \mathcal{F}^2_{b,b})}
      &\text{ }
      &\text{ }
      \\
      x\otimes y\otimes \mathcal{F}(b\tensor b)
      \ar[d]^-{\I_x\otimes t^{\mathcal{F}}_{y,b\tensor b}}
      \ar@/^0.5pc/[dr]^(0.6){\I_x\otimes (\I_y\otimes \mathcal{F}(\jmath_{b}))}_-{\cong}
      &\text{ }
      &\text{ }
      \\
      x\otimes \mathcal{F}(y\otimes b\tensor b)
      \ar[d]^-{\I_x\otimes \mathcal{F}(\I_y\otimes \jmath_{b})}_-{\cong}
      &x\otimes y\otimes \mathcal{F}(b)
      \ar[d]^-{\I_x\otimes t^{\mathcal{F}}_{y,b}}
      \ar@{=}[r]
      &x\otimes y\otimes \mathcal{F}(b)
      \ar[dd]^-{t^{\mathcal{F}}_{x\otimes y,b}}
      \\
      x\otimes \mathcal{F}(y\otimes b)
      \ar[d]^-{\I_x\otimes \mathcal{F}(\gamma_{y})}
      \ar@{=}[r]
      &x\otimes \mathcal{F}(y\otimes b)
      \ar[d]^-{t^{\mathcal{F}}_{x,y\otimes b}}
      &\text{ }
      \\
      x\otimes \mathcal{F}(y)
      \ar[d]^-{t^{\mathcal{F}}_{x,y}}
      &\mathcal{F}(x\otimes y\otimes b)
      \ar[d]^-{\mathcal{F}(\I_x\otimes \gamma_{y})}
      \ar@{=}[r]
      &\mathcal{F}(x\otimes y\otimes b)
      \ar[d]^-{\mathcal{F}(\gamma_{x\otimes y})}
      \\
      \mathcal{F}(x\otimes y)
      \ar[d]^-{\mathcal{F}(\textit{cq}_{x,y})}
      \ar@{=}[r]
      &\mathcal{F}(x\otimes y)
      \ar[d]^-{\mathcal{F}(\textit{cq}_{x,y})}
      \ar@{=}[r]
      &\mathcal{F}(x\otimes y)
      \ar[d]^-{\mathcal{F}(\textit{cq}_{x,y})}
      \\
      \mathcal{F}(x\tensor y)
      \ar@{=}[r]
      &\mathcal{F}(x\tensor y)
      \ar@{=}[r]
      &\mathcal{F}(x\tensor y)
    }      
  }}
  \!\!\!\!
  \!\!\!\!
  \cdots
\end{equation*}
Upto this point, we treated
$\mathcal{F}(b)$ as the underlying object in $\mathcal{V}$
of the right $b'$-module.
From now on, we will treat
$\mathcal{F}(b)$ as the underlying object of the right $b$-module,
whose right $b$-action is induced from the morphism
$f:b\to \mathcal{F}(b)$
in $\textit{Comm}(\mathcal{V})$.
We further continue.
\begin{equation*}
  \cdots
  \!\!\!\!
  \!\!\!\!
  \!\!\!\!
  \vcenter{\hbox{
    \xymatrix@R=15pt@C=15pt{
      x\otimes \mathcal{F}(b)\otimes y\otimes \mathcal{F}(b)
      \ar[d]^-{\I_x\otimes s_{\mathcal{F}(b),y}\otimes \I_{\mathcal{F}(b)}}_-{\cong}
      \ar@{=}[r]
      &x\otimes \mathcal{F}(b)\otimes y\otimes \mathcal{F}(b)
      \ar[d]^-{\I_x\otimes s_{\mathcal{F}(b),y}\otimes \I_{\mathcal{F}(b)}}_-{\cong}
      \ar@{=}[r]
      &x\otimes \mathcal{F}(b)\otimes y\otimes \mathcal{F}(b)
      \ar@{->>}[d]^-{\I_{x\otimes \mathcal{F}(b)}\otimes \textit{cq}_{y,\mathcal{F}(b)}}
      \\
      x\otimes y\otimes \mathcal{F}(b)\otimes \mathcal{F}(b)
      \ar[dd]^-{\I_{x\otimes y}\otimes m_{\mathcal{F}(b)}}
      \ar@{=}[r]
      &x\otimes y\otimes \mathcal{F}(b)\otimes \mathcal{F}(b)
      \ar@{->>}[d]^-{\I_{x\otimes y}\otimes \textit{cq}_{\mathcal{F}(b),\mathcal{F}(b)}}
      &x\otimes \mathcal{F}(b)\otimes (y\tensor \mathcal{F}(b))
      \ar@{->>}[d]^-{\textit{cq}_{x,\mathcal{F}(b)}\otimes \I_{y\tensor\mathcal{F}(b)}}
      \\
      \text{ }
      &x\otimes y\otimes \mathcal{F}(b)\tensor \mathcal{F}(b)
      \ar[dl]|-{\I_{x\otimes y}\otimes m_{\mathcal{F}(b)/b}}
      \ar@{->>}[d]^-{\textit{cq}_{x,y}\otimes \I_{\mathcal{F}(b)\tensor\mathcal{F}(b)}}
      &(x\tensor \mathcal{F}(b))\otimes (y\tensor \mathcal{F}(b))
      \ar@{->>}[d]^-{\textit{cq}_{x\tensor\mathcal{F}(b),y\tensor\mathcal{F}(b)}}
      \\
      x\otimes y\otimes \mathcal{F}(b)
      \ar[d]^-{t^{\mathcal{F}}_{x\otimes y,b}}
      \ar@{->>}[dr]|(0.5){\textit{cq}_{x,y}\otimes \I_{\mathcal{F}(b)}}
      &(x\tensor y)\otimes (\mathcal{F}(b)\tensor \mathcal{F}(b))
      \ar[d]^-{\I_{x\tensor y}\otimes m_{\mathcal{F}(b)/b}}
      \ar@/^0.5pc/@{->>}[dr]|-{\textit{cq}_{x\tensor y,\mathcal{F}(b)\tensor\mathcal{F}(b)}}
      &x\tensor \mathcal{F}(b)\tensor y\tensor \mathcal{F}(b)
      \ar[d]^-{\I_x\tensor s_{\mathcal{F}(b),y}\tensor\I_{\mathcal{F}(b)}}_-{\cong}
      \\
      \mathcal{F}(x\otimes y\otimes b)
      \ar[d]^-{\mathcal{F}(\gamma_{x\otimes y})}
      \ar[dr]|-{\mathcal{F}(\textit{cq}_{x,y}\otimes \I_b)}
      &(x\tensor y)\otimes \mathcal{F}(b)
      \ar[d]^-{t^{\mathcal{F}}_{x\tensor y,b}}
      \ar@/^0.5pc/@{->>}[dr]|-{\textit{cq}_{x\tensor y,\mathcal{F}(b)}}
      &x\tensor y\tensor \mathcal{F}(b)\tensor \mathcal{F}(b)
      \ar[d]^-{\I_{x\tensor y}\tensor m_{\mathcal{F}(b)/b}}
      \\
      \mathcal{F}(x\otimes y)
      \ar[d]^-{\mathcal{F}(\textit{cq}_{x,y})}
      &\mathcal{F}((x\tensor y)\otimes b)
      \ar[d]^-{\mathcal{F}(\gamma_{x\tensor y})}
      &x\tensor y\tensor \mathcal{F}(b)
      \ar[d]^-{\lambda^{\mathcal{F}}_{x\tensor y}}
      \\
      \mathcal{F}(x\tensor y)
      \ar@{=}[r]
      &\mathcal{F}(x\tensor y)
      \ar@{=}[r]
      &\mathcal{F}(x\tensor y)
    }
  }}
  \!\!\!\!
  \cdots
\end{equation*}
We used the defining relation
(\ref{eq EnEWThm defining lambda})
of $\lambda^{\mathcal{F}}_{x\tensor y}$.
Finally, we finish our calculation as follows.
\begin{equation*}
  \cdots
  \vcenter{\hbox{
    \xymatrix@R=10pt@C=50pt{
      x\otimes \mathcal{F}(b)\otimes y\otimes \mathcal{F}(b)
      \ar@{->>}[d]^-{\I_{x\otimes \mathcal{F}(b)}\otimes \textit{cq}_{y,\mathcal{F}(b)}}
      \ar@{=}[r]
      &x\otimes \mathcal{F}(b)\otimes y\otimes \mathcal{F}(b)
      \ar@{->>}[d]^-{\textit{cq}'_{x\otimes \mathcal{F}(b),y\otimes \mathcal{F}(b)}}
      \\
      x\otimes \mathcal{F}(b)\otimes (y\tensor \mathcal{F}(b))
      \ar@{->>}[d]^-{\textit{cq}_{x,\mathcal{F}(b)}\otimes \I_{y\tensor\mathcal{F}(b)}}
      \ar@{->>}@/^0.5pc/[dr]|-{\textit{cq}'_{x\otimes \mathcal{F}(b),y\tensor \mathcal{F}(b)}}
      &(x\otimes \mathcal{F}(b))\tensor' (y\otimes \mathcal{F}(b))
      \ar@{->>}[d]^-{\I_{x\otimes \mathcal{F}(b)}\tensor'\textit{cq}_{y,\mathcal{F}(b)}}
      \\
      (x\tensor \mathcal{F}(b))\otimes (y\tensor \mathcal{F}(b))
      \ar@{->>}[d]^-{\textit{cq}_{x\tensor\mathcal{F}(b),y\tensor\mathcal{F}(b)}}
      \ar@{->>}@/^0.5pc/[dr]|-{\textit{cq}'_{x\tensor \mathcal{F}(b),y\tensor \mathcal{F}(b)}}
      &(x\otimes \mathcal{F}(b))\tensor' (y\tensor \mathcal{F}(b))
      \ar@{->>}[d]^-{\textit{cq}_{x,\mathcal{F}(b)}\tensor'\I_{y\tensor \mathcal{F}(b)}}
      \\
      x\tensor \mathcal{F}(b)\tensor y\tensor \mathcal{F}(b)
      \ar[d]^-{\I_x\tensor s_{\mathcal{F}(b),y}\tensor\I_{\mathcal{F}(b)}}_-{\cong}
      &(x\tensor \mathcal{F}(b))\tensor' (y\tensor \mathcal{F}(b))
      \ar[dd]^-{(f'_*f^*)^2_{x,y}}
      \\
      x\tensor y\tensor \mathcal{F}(b)\tensor \mathcal{F}(b)
      \ar[d]^-{\I_{x\tensor y}\tensor m_{\mathcal{F}(b)/b}}
      &\text{ }
      \\
      x\tensor y\tensor \mathcal{F}(b)
      \ar[d]^-{\lambda^{\mathcal{F}}_{x\tensor y}}
      \ar@{=}[r]
      &x\tensor y\tensor \mathcal{F}(b)
      \ar[d]^-{\lambda^{\mathcal{F}}_{x\tensor y}}
      \\
      \mathcal{F}(x\tensor y)
      \ar@{=}[r]
      &\mathcal{F}(x\tensor y)
    }
  }}
\end{equation*}
We used the description (\ref{eq SixFunctors f'*f* tensor coherence})
of the morphism $(f'_*f^*)^2_{x,y}$.
After right-cancelling the epimorphisms
$\textit{cq}'_{x\otimes \mathcal{F}(b),y\otimes \mathcal{F}(b)}$,
$\textit{cq}_{x,\mathcal{F}(b)}\tensor'\textit{cq}_{y,\mathcal{F}(b)}$
in our calculation,
we obtain (\ref{eq MainResult lambdaF tensor}).

\subsection{Proposition~\ref{prop DayConvolution monoidalproduct}}
\label{subsec Appendix DayConvolution monoidalproduct}
Let $b$, $b'$ be commutative monoids in $\mathcal{V}$
and let $x$, $y$ be $(b,b')$-bimodules.
We verify the universal property of the enriched natural transformation
$\theta$ described in Proposition~\ref{prop DayConvolution monoidalproduct}.
Suppose we are given an enriched functor
$\mathcal{F}:\textit{Mod}_b\to \textit{Mod}_{b'}$
and an enriched natural transformation
$\alpha_{-,-}:
(-\tensor x)\tensor' (-\tensor y)\Rightarrow \mathcal{F}(-\tensor -)$.
Let $v\in \mathcal{V}$, $w$, $z\in\textit{Mod}_b$.
We have the following relations
since $\alpha_{w,z}$
is enriched natural in variables $w$, $z$.
\begin{equation}\label{eq DayConvolution alphawbCnat}
  \vcenter{\hbox{
    \xymatrix@R=10pt@C=13pt{
      v\!\otimes\! (w\!\tensor\! x)\!\tensor'\! (z\tensor y)
      \ar@{=}[d]
      \ar[r]^-{\I_v\otimes \alpha_{w,z}}
      &v\!\otimes\! \mathcal{F}(w\!\tensor\! z)
      \ar[d]^-{t^{\mathcal{F}}_{v,w\tensor z}}
      \\
      (v\!\otimes\! w\!\tensor\! x)\!\tensor'\! (z\!\tensor\! y)
      \ar[r]^-{\alpha_{v\otimes w,z}}
      &\mathcal{F}(v\!\otimes\! w\!\tensor\! z)
    }
  }}
  \text{ }
  \vcenter{\hbox{
    \xymatrix@R=8pt@C=13pt{
      v\!\otimes\! (w\!\tensor\! x)\!\tensor'\! (z\!\tensor\! y)
      \ar[dd]_-{t^{(w\tensor x)\tensor'-}_{v,z\tensor y}}^-{\cong}
      \ar[r]^-{\I_v\otimes \alpha_{w,z}}
      &v\!\otimes\! \mathcal{F}(w\!\tensor\! z)
      \ar[d]^-{t^{\mathcal{F}}_{v,w\tensor z}}
      \\
      \text{ }
      &\mathcal{F}(v\!\otimes\! w\!\tensor\! z)
      \ar[d]^-{\mathcal{F}(t^{w\tensor -}_{v,z})}_-{\cong}
      \\
      (w\!\tensor\! x)\!\tensor'\! (v\!\otimes\! z\!\tensor\! y)
      \ar[r]^-{\alpha_{w,v\otimes z}}
      &\mathcal{F}(w\!\tensor\! (v\!\otimes\! z))
    }
  }}
\end{equation}
Recall the morphism
$\phi:x\tensor' y\to \mathcal{F}(b)$
in $\textit{Mod}_{b'}$
defined in (\ref{eq DayConvolution monoidalproduct phidef}).
Using the relations
(\ref{eq DayConvolution alphawbCnat}),
one can check that the underlying morphism of
$\phi$ in $\mathcal{V}$
satisfies the relations below.
\begin{equation} \label{eq DayConvolution alphawbbb}
  \vcenter{\hbox{
    \xymatrix@R=9pt@C=17pt{
      w\otimes x\otimes y
      \ar@{->>}[dd]_-{\I_w\otimes \textit{cq}'_{x,y}}
      \ar@{->>}[r]^-{\textit{cq}_{w,x}\otimes \I_y}
      &(w\tensor x)\otimes y
      \ar@{->>}[d]^-{\textit{cq}'_{w\tensor x,y}}
      &x\otimes z\otimes y
      \ar[d]_-{s_{x,z}\otimes \I_y}^-{\cong}
      \ar@{->>}[r]^-{\I_x\otimes \textit{cq}_{z,y}}
      &x\otimes (z\tensor y)
      \ar@{->>}[d]^-{\textit{cq}'_{x,z\tensor y}}
      \\
      \text{ }
      &(w\tensor x)\tensor' y
      \ar[d]^-{\I_{w\tensor x}\tensor' \imath^{-1}_y}_-{\cong}
      &z\otimes x\otimes y
      \ar@{->>}[d]_-{\I_z\otimes \textit{cq}'_{x,y}}
      &x\tensor' (z\tensor y)
      \ar[d]^-{\imath^{-1}_x\tensor' \I_{z\tensor y}}_-{\cong}
      \\
      w\otimes x\tensor' y
      \ar[d]_-{\I_w\otimes \phi}
      &(w\tensor x)\tensor' (b\tensor y)
      \ar[d]^-{\alpha_{w,b}}
      &z\otimes x\tensor' y
      \ar[d]_-{\I_z\otimes \phi}
      &(b\tensor x)\tensor' (z\tensor y)
      \ar[d]^-{\alpha_{b,z}}
      \\
      w\otimes \mathcal{F}(b)
      \ar[d]_-{t^{\mathcal{F}}_{w,b}}
      &\mathcal{F}(w\tensor b)
      \ar[d]^-{\mathcal{F}(\jmath_{w})}_-{\cong}
      &z\otimes \mathcal{F}(b)
      \ar[d]_-{t^{\mathcal{F}}_{z,b}}
      &\mathcal{F}(b\tensor z)
      \ar[d]^-{\mathcal{F}(\imath_{z})}_-{\cong}
      \\
      \mathcal{F}(w\otimes b)
      \ar[r]^-{\mathcal{F}(\gamma_{w})}
      &\mathcal{F}(w)
      &\mathcal{F}(z\otimes b)
      \ar[r]^-{\mathcal{F}(\gamma_{z})}
      &\mathcal{F}(z)
    }
  }}
\end{equation}
We claim that there exists a unique morphism 
$\psi:x\wtensor y\to \mathcal{F}(b)$
of $(b,b')$-bimodules, whose underlying morphism 
in $\mathcal{V}$ satisfies the relation
(\ref{eq DayConvolution monoidalproduct psidef}).
The following continued diagrams imply that
the morphism
$\psi:x\wtensor y\to \mathcal{F}(b)$
in $\mathcal{V}$ is well-defined.
\begin{equation*}
  \vcenter{\hbox{
    \xymatrix@R=7pt@C=8pt{
      x\otimes b\otimes y
      \ar[dd]_-{\mathring{\rho}_x\otimes \I_y}
      \ar@{=}[r]
      &x\otimes b\otimes y
      \ar[d]^-{s_{x,b}\otimes \I_y}_-{\cong}
      \ar@{=}[r]
      &x\otimes b\otimes y
      \ar[d]^-{s_{x,b}\otimes \I_y}_-{\cong}
      \ar@{=}[r]
      &x\otimes b\otimes y
      \ar[d]^-{s_{x,b}\otimes \I_y}_-{\cong}
      \\
      \text{ }
      &b\otimes x\otimes y
      \ar@/_0.5pc/[dl]|-{\rho_x\otimes \I_y}
      \ar@{->>}[d]^-{\textit{cq}'_{b\otimes x,y}}
      \ar@{=}[r]
      &b\otimes x\otimes y
      \ar@{->>}[d]^-{\I_b\otimes \textit{cq}'_{x,y}}
      \ar@{=}[r]
      &b\otimes x\otimes y
      \ar@{->>}[d]^-{\I_b\otimes \textit{cq}'_{x,y}}
      \\
      x\otimes y
      \ar@{->>}[d]_-{\textit{cq}'_{x,y}}
      &b\otimes x\tensor' y
      \ar@/_0.5pc/[dl]|-{\rho_{x}\tensor' \I_y}
      \ar[d]^-{(\I_b\otimes \imath^{-1}_x)\tensor' \imath^{-1}_y}_-{\cong}
      \ar@{=}[r]
      &b\otimes x\tensor' y
      \ar[d]^-{\I_b\otimes (\imath^{-1}_x\tensor' \imath^{-1}_y)}_-{\cong}
      \ar@{=}[r]
      &b\otimes x\tensor' y
      \ar[ddd]^-{\I_b\otimes \phi}
      \\
      x\tensor' y
      \ar[d]^-{\imath^{-1}_x\tensor'\imath^{-1}_y}_-{\cong}
      \ar@{.>}@/_2pc/@<-2ex>[dddd]|-{\phi}
      &(b\otimes b\tensor x)\tensor' (b\tensor y)
      \ar@/^0.5pc/[dl]|-{(m_b\tensor\I_x)\tensor' \I_{b\tensor y}}
      \ar[dd]^-{\alpha_{b\otimes b,b}}
      \ar@{=}[r]
      &b\otimes (b\tensor x)\tensor' (b\tensor y)
      \ar[d]^-{\I_b\otimes \alpha_{b,b}}
      &\text{ }
      \\
      (b\tensor x)\tensor' (b\tensor y)
      \ar[dd]^-{\alpha_{b,b}}
      &\text{ }
      &b\otimes \mathcal{F}(b\tensor b)
      \ar[d]^-{t^{\mathcal{F}}_{b,b\tensor b}}
      \ar@/^0.5pc/[dr]^-{\I_b\otimes \mathcal{F}(\jmath_{b})}_-{\cong}
      &\text{ }
      \\
      \text{ }
      &\mathcal{F}(b\otimes b\tensor b)
      \ar@/^0.5pc/[dl]|-{\mathcal{F}(m_b\tensor\I_{b})}
      \ar[d]^-{\mathcal{F}(\jmath_{b\otimes b})}_-{\cong}
      \ar@{=}[r]
      &\mathcal{F}(b\otimes b\tensor b)
      \ar[d]^-{\mathcal{F}(\I_b\otimes \jmath_{b})}_-{\cong}
      &b\otimes \mathcal{F}(b)
      \ar[d]^-{t^{\mathcal{F}}_{b,b}}
      \\
      \mathcal{F}(b\tensor b)
      \ar[d]^-{\mathcal{F}(\jmath_{b})}_-{\cong}
      &\mathcal{F}(b\otimes b)
      \ar[d]^-{\mathcal{F}(m_b)}
      \ar@{=}[r]
      &\mathcal{F}(b\otimes b)
      \ar[d]^-{\mathcal{F}(m_b)}
      \ar@{=}[r]
      &\mathcal{F}(b\otimes b)
      \ar[d]^-{\mathcal{F}(m_b)}
      \\
      \mathcal{F}(b)
      \ar@{=}[r]
      &\mathcal{F}(b)
      \ar@{=}[r]
      &\mathcal{F}(b)
      \ar@{=}[r]
      &\mathcal{F}(b)
    }
  }}
  \!\!\!\!
  \!\!
  \cdots
\end{equation*}
We used the left relation in (\ref{eq DayConvolution alphawbCnat})
for $v=b$, $w=z=b$.
We continue.
\begin{equation*}
  \cdots
  \vcenter{\hbox{
    \xymatrix@R=7pt@C=40pt{
      x\otimes b\otimes y
      \ar[d]^-{s_{x,b}\otimes \I_y}_-{\cong}
      \ar@{=}[r]
      &x\otimes b\otimes y
      \ar@{->>}[d]^-{\I_x\otimes \textit{cq}_{b,y}}
      \ar@{=}[r]
      &x\otimes b\otimes y
      \ar[dd]^-{\I_x\otimes \rho_y}
      \\
      b\otimes x\otimes y
      \ar@{->>}[d]^-{\I_b\otimes \textit{cq}'_{x,y}}
      &x\otimes b\tensor y
      \ar@{->>}[d]^-{\textit{cq}'_{x,b\tensor y}}
      \ar@/^0.5pc/[dr]^-{\I_x\otimes \imath_{y}}_-{\cong}
      &\text{ }
      \\
      b\otimes x\tensor' y
      \ar[d]^-{\I_b\otimes \phi}
      &x\tensor' (b\tensor y)
      \ar[d]^-{\imath^{-1}_x\tensor' \I_{b\tensor y}}_-{\cong}
      \ar@/^0.5pc/[dr]^-{\I_x\tensor' \imath_y}_-{\cong}
      &x\otimes y
      \ar@{->>}[d]^-{\textit{cq}'_{x,y}}
      \\
      b\otimes \mathcal{F}(b)
      \ar[d]^-{t^{\mathcal{F}}_{b,b}}
      &(b\tensor x)\tensor' (b\tensor y)
      \ar[d]^-{\alpha_{b,b}}
      &x\tensor' y
      \ar[dd]^-{\phi}
      \\
      \mathcal{F}(b\otimes b)
      \ar[d]^-{\mathcal{F}(m_b)}
      &\mathcal{F}(b\tensor b)
      \ar[d]^-{\mathcal{F}(\jmath_{b})}_-{\cong}
      &\text{ }
      \\
      \mathcal{F}(b)
      \ar@{=}[r]
      &\mathcal{F}(b)
      \ar@{=}[r]
      &\mathcal{F}(b)
    }
  }}
\end{equation*}
We used the right relation in (\ref{eq DayConvolution alphawbbb}) for $z=b$.
Readers should check that
$(m_b\tensor\I_{x})\circ a'_{b,b,x}\circ(\I_b\otimes \imath^{-1}_x)
=\textit{cq}_{b,x}
=\imath_x^{-1}\circ \rho_x$.
One can check that $\psi$ is compatible with right $b'$-actions.
By right-cancelling the epimorphism
$\I_b\otimes \widetilde{\textit{cq}}_{x,y}$
in the diagram below
\begin{equation*}
  \xymatrix@R=7pt@C=17pt{
    b\otimes x\otimes y
    \ar@{->>}[d]_-{\I_b\otimes \widetilde{\textit{cq}}_{x,y}}
    \ar@{=}[r]
    &b\otimes x\otimes y
    \ar@{->>}[d]^-{\I_b\otimes \textit{cq}'_{x,y}}
    \ar@{=}[r]
    &b\otimes x\otimes y
    \ar@{->>}[d]^-{\textit{cq}_{b,x}\otimes \I_y}
    \ar@{=}[r]
    &b\otimes x\otimes y
    \ar[dd]^-{\rho_x\otimes \I_y}
    \ar@{=}[r]
    &b\otimes x\otimes y
    \ar@{->>}[d]^-{\I_b\otimes \widetilde{\textit{cq}}_{x,y}}
    \\
    b\otimes x\wtensor y
    \ar[dd]_-{\I_b\otimes \psi}
    &b\otimes x\tensor' y
    \ar[dd]^-{\I_b\otimes \phi}
    &(b\tensor x)\otimes y
    \ar@{->>}[d]^-{\textit{cq}'_{b\tensor x ,y}}
    \ar@/^0.5pc/[dr]^-{\imath_{x}\otimes \I_y}_-{\cong}
    &\text{ }
    &b\otimes x\wtensor y
    \ar[dd]^-{\rho_{x\widetilde{\tensor} y}}
    \\
    \text{ }
    &\text{ }
    &(b\tensor x)\tensor' y
    \ar[d]^-{\I_{b\tensor x}\tensor'\imath^{-1}_y}_-{\cong}
    \ar@/^0.5pc/[dr]^-{\imath_x\tensor' \I_y}_-{\cong}
    &x\otimes y
    \ar@{->>}[d]^-{\textit{cq}'_{x,y}}
    \ar@{->>}@/^0.5pc/[dr]|-{\widetilde{\textit{cq}}_{x,y}}
    &\text{ }
    \\
    b\otimes \mathcal{F}(b)
    \ar[dd]_-{\rho_{\mathcal{F}(b)}}
    \ar@{=}[r]
    &b\otimes \mathcal{F}(b)
    \ar[d]^-{t^{\mathcal{F}}_{b,b}}
    &(b\tensor x)\tensor' (b\tensor y)
    \ar[d]^-{\alpha_{b,b}}
    &x\tensor' y
    \ar[dd]^-{\phi}
    &x\wtensor y
    \ar[dd]^-{\psi}
    \\
    \text{ }
    &\mathcal{F}(b\otimes b)
    \ar[d]^-{\mathcal{F}(m_b)}
    &\mathcal{F}(b\tensor b)
    \ar[d]^-{\mathcal{F}(\jmath_{b})}_-{\cong}
    &\text{ }
    &\text{ }
    \\
    \mathcal{F}(b)
    \ar@{=}[r]
    &\mathcal{F}(b)
    \ar@{=}[r]
    &\mathcal{F}(b)
    \ar@{=}[r]
    &\mathcal{F}(b)
    \ar@{=}[r]
    &\mathcal{F}(b)
  }
\end{equation*}
we obtain that $\psi$ is compatible with left $b$-actions.
We used the left relation in (\ref{eq DayConvolution alphawbbb})
for $w=b$
and the defining relation (\ref{eq DayConvolution monoidalproduct psidef}) of $\psi$.
Thus
$\psi:x\wtensor y\to \mathcal{F}(b)$
is a well-defined morphism of $(b,b')$-bimodules.

Let $w$, $z\in\textit{Mod}_b$.
We show that the enriched natural transformation
$\xi=\lambda^{\mathcal{F}}\circ (-\tensor \psi):
-\tensor (x\wtensor y)\Rightarrow\mathcal{F}$
satisfies the relation
$\alpha_{w,z}=\xi_{w\tensor z}\circ \theta_{w,z}$.
We begin our calculation.
\begin{equation*}
  \vcenter{\hbox{
    \xymatrix@R=16pt@C=15pt{
      w\otimes x\otimes z\otimes y
      \ar@{->>}[dd]_-{\textit{cq}_{w,x}\otimes \textit{cq}_{z,y}}
      \ar@{=}[rr]
      &\text{ }
      &w\otimes x\otimes z\otimes y
      \ar[d]^-{\I_w\otimes s_{x,z}\otimes \I_y}_-{\cong}
      \ar@{=}[r]
      &w\otimes x\otimes z\otimes y
      \ar[d]^-{\I_w\otimes s_{x,z}\otimes \I_y}_-{\cong}
      \\
      \text{ }
      &\text{ }
      &w\otimes z\otimes x\otimes y
      \ar@{->>}[d]^-{\I_{w\otimes z}\otimes \textit{cq}_{x,y}}
      \ar@{=}[r]
      &w\otimes z\otimes x\otimes y
      \ar@{->>}[dd]^-{\I_{w\otimes z}\otimes \widetilde{\textit{cq}}_{x,y}}
      \\
      (w\tensor x)\otimes (z\tensor y)
      \ar@{->>}[d]_-{\textit{cq}'_{w\tensor x,z\tensor y}}
      \ar@/^0.5pc/@{->>}[dr]|-{\textit{cq}_{w\tensor x,z\tensor y}}
      &\text{ }
      &w\otimes z\otimes (x\tensor y)
      \ar@{->>}[d]^-{\textit{cq}_{w,z}\otimes \I_{x\tensor y}}
      \ar@/^0.5pc/[dr]|-{\I_{w\otimes z}\otimes (p_*)^2_{x,y}}
      &\text{ }
      \\
      (w\tensor x)\tensor' (z\tensor y)
      \ar[dd]_-{\theta_{w,z}}
      &w\tensor x\tensor z\tensor y
      \ar@/_0.5pc/[dr]_-{\I_w\tensor s_{x,z}\tensor\I_y}^-{\cong}
      &(w\tensor z)\otimes (x\tensor y)
      \ar@{->>}[d]^-{\textit{cq}_{w\tensor z,x\tensor y}}
      \ar@/^0.5pc/[dr]|-{\I_{w\tensor z}\otimes (p_*)^2_{x,y}}
      &w\otimes z\otimes (x\wtensor y)
      \ar@{->>}[d]^-{\textit{cq}_{w,z}\otimes \I_{x\widetilde{\tensor} y}}
      \\
      \text{ }      
      &\text{ }
      &w\tensor z\tensor x\tensor y
      \ar[d]^-{\I_{w\tensor z}\tensor(p_*)^2_{x,y}}
      &(w\tensor z)\otimes (x\wtensor y)
      \ar@{->>}[d]^-{\textit{cq}_{w\tensor z,x\widetilde{\tensor} y}}
      \\
      w\tensor z\tensor (x\wtensor y)
      \ar[d]_-{\xi_{w\tensor z}}
      \ar@{=}[rr]
      &\text{ }
      &w\tensor z\tensor (x\wtensor y)
      \ar[d]^-{\xi_{w\tensor z}}
      \ar@{=}[r]
      &w\tensor z\tensor (x\wtensor y)
      \ar[d]^-{\xi_{w\tensor z}}
      \\
      \mathcal{F}(w\tensor z)
      \ar@{=}[rr]
      &\text{ }
      &\mathcal{F}(w\tensor z)
      \ar@{=}[r]
      &\mathcal{F}(w\tensor z)
    }
  }}
  \cdots
\end{equation*}
We used the defining relation
(\ref{eq DayConvolution universaltheta})
of $\theta$
as well as the defining relation
(\ref{eq DayConvolution p*p'* formula})
of $(p_*)^2$.
We continue our calculation as follows.
\begin{equation*}
  \cdots
  \!\!\!\!
  \vcenter{\hbox{
    \xymatrix@R=16pt@C=15pt{
      w\otimes x\otimes z\otimes y
      \ar[d]^-{\I_w\otimes s_{x,z}\otimes \I_y}_-{\cong}
      \ar@{=}[rr]
      &\text{ }
      &w\otimes x\otimes z\otimes y
      \ar[d]^-{\I_w\otimes s_{x,z}\otimes \I_y}_-{\cong}
      \ar@{=}[r]
      &w\otimes x\otimes z\otimes y
      \ar[d]^-{\I_w\otimes s_{x,z}\otimes \I_y}_-{\cong}
      \\
      w\otimes z\otimes x\otimes y
      \ar@{->>}[d]^-{\I_{w\otimes z}\otimes \widetilde{\textit{cq}}_{x,y}}
      \ar@{=}[rr]
      &\text{ }
      &w\otimes z\otimes x\otimes y
      \ar@{->>}[d]^-{\I_{w\otimes z}\otimes \textit{cq}'_{x,y}}
      \ar@{=}[r]
      &w\otimes z\otimes x\otimes y
      \ar@{->>}[d]^-{\I_w\otimes \I_z\otimes \textit{cq}'_{x,y}}
      \\
      w\otimes z\otimes (x\wtensor y)
      \ar@{->>}[dd]^-{\textit{cq}_{w,z}\otimes \I_{x\widetilde{\tensor} y}}
      \ar@/^0.5pc/[dr]|-{\I_{w\otimes z}\otimes \psi}
      &\text{ }
      &w\otimes z\otimes (x\tensor' y)
      \ar[d]^-{\I_{w\otimes z}\otimes \phi}
      \ar@{=}[r]
      &w\otimes z\otimes (x\tensor' y)
      \ar[d]^-{\I_w\otimes \I_z\otimes \phi}
      \\
      \text{ }
      &w\otimes z\otimes \mathcal{F}(b)
      \ar@{->>}[dd]^-{\textit{cq}_{w,z}\otimes \I_{\mathcal{F}(b)}}
      \ar@{=}[r]
      &w\otimes z\otimes \mathcal{F}(b)
      \ar[dd]^-{t^{\mathcal{F}}_{w\otimes z,b}}
      \ar@{=}[r]
      &w\otimes z\otimes \mathcal{F}(b)
      \ar[d]^-{\I_w\otimes t^{\mathcal{F}}_{z,b}}
      \\
      (w\tensor z)\otimes (x\wtensor y)
      \ar@{->>}[d]^-{\textit{cq}_{w\tensor z,x\widetilde{\tensor} y}}
      \ar@/^0.5pc/[dr]|-{\I_{w\tensor z}\otimes \psi}
      &\text{ }
      &\text{ }
      &w\otimes \mathcal{F}(z\otimes b)
      \ar[d]^-{t^{\mathcal{F}}_{w,z\otimes b}}
      \\
      w\tensor z\tensor (x\wtensor y)
      \ar[dd]^-{\xi_{w\tensor z}}
      \ar@/^0.5pc/[dr]|-{\I_{w\tensor z}\tensor\psi}
      &(w\tensor z)\otimes \mathcal{F}(b)
      \ar@{->>}[d]^-{\textit{cq}_{w\tensor z,\mathcal{F}(b)}}
      \ar@/^0.5pc/[dr]|-{t^{\mathcal{F}}_{w\tensor z,b}}
      &\mathcal{F}(w\otimes z\otimes b)
      \ar[d]^-{\mathcal{F}(\textit{cq}_{w,z}\otimes \I_{b})}
      \ar@/^0.5pc/[dr]|(0.65){\mathcal{F}(\gamma_{w\otimes z})}
      \ar@{=}[r]
      &\mathcal{F}(w\otimes z\otimes b)
      \ar[d]^-{\mathcal{F}(\I_w\otimes \gamma_{z})}
      \\
      \text{ }
      &w\tensor z\tensor \mathcal{F}(b)
      \ar[d]^-{\lambda^{\mathcal{F}}_{w\tensor z}}
      &\mathcal{F}((w\tensor z)\otimes b)
      \ar[d]^-{\mathcal{F}(\gamma_{w\tensor z})}
      &\mathcal{F}(w\otimes z)
      \ar[d]^-{\mathcal{F}(\textit{cq}_{w,z})}
      \\
      \mathcal{F}(w\tensor z)
      \ar@{=}[r]
      &\mathcal{F}(w\tensor z)
      \ar@{=}[r]
      &\mathcal{F}(w\tensor z)
      \ar@{=}[r]
      &\mathcal{F}(w\tensor z)
    }
  }}
  \!\!\!\!
  \cdots
\end{equation*}
We used the definition
of $\xi_{w\tensor z}$,
the defining relation
(\ref{eq EnEWThm defining lambda})
of $\lambda^{\mathcal{F}}_{w\tensor z}$
and the defining relation
(\ref{eq DayConvolution monoidalproduct psidef})
of $\psi$.
When moving from the 3rd vertical composition to the last,
the position of the brackets changed.
We further continue our calculation as follows.
\begin{equation*}
  \cdots
  \vcenter{\hbox{
    \xymatrix@R=18pt@C=25pt{
      w\otimes x\otimes z\otimes y
      \ar[d]^-{\I_w\otimes s_{x,z}\otimes \I_y}_-{\cong}
      \ar@{=}[r]
      &w\otimes x\otimes z\otimes y
      \ar@{->>}[d]^-{\I_w\otimes \I_x\otimes \textit{cq}_{z,y}}
      \ar@{=}[r]
      &w\otimes x\otimes z\otimes y
      \ar@{->>}[d]^-{\I_{w\otimes x}\otimes \textit{cq}_{z,y}}
      \\
      w\otimes z\otimes x\otimes y
      \ar@{->>}[d]^-{\I_w\otimes \I_z\otimes \textit{cq}'_{x,y}}
      &w\otimes x\otimes (z\tensor y)
      \ar@{->>}[d]^-{\I_w\otimes \textit{cq}'_{x,z\tensor y}}
      \ar@{=}[r]
      &w\otimes x\otimes (z\tensor y)
      \ar@{->>}[d]^-{\textit{cq}'_{w\otimes x,z\tensor y}}
      \\
      w\otimes z\otimes (x\tensor' y)
      \ar[d]^-{\I_w\otimes \I_z\otimes \phi}
      &w\otimes x\tensor' (z\tensor y)
      \ar[d]^-{\I_w\otimes \imath^{-1}_x \tensor' \I_{z\tensor y}}_-{\cong}
      \ar@{=}[r]
      &w\otimes x\tensor' (z\tensor y)
      \ar[d]^-{(\I_w\otimes \imath^{-1}_x) \tensor' \I_{z\tensor y}}_-{\cong}
      \\
      w\otimes z\otimes \mathcal{F}(b)
      \ar[d]^-{\I_w\otimes t^{\mathcal{F}}_{z,b}}
      &w\otimes (b\tensor x)\tensor' (z\tensor y)
      \ar[d]^-{\I_w\otimes \alpha_{b,z}}
      \ar@{=}[r]
      &(w\otimes b\tensor x)\tensor' (z\tensor y)
      \ar[dd]^-{\alpha_{w\otimes b,z}}
      \\
      w\otimes \mathcal{F}(z\otimes b)
      \ar[d]^-{t^{\mathcal{F}}_{w,z\otimes b}}
      \ar@/^0.5pc/[dr]|-{\I_w\otimes \mathcal{F}(\gamma_{z})}
      &w\otimes \mathcal{F}(b\tensor z)
      \ar[d]^-{\I_w\otimes \mathcal{F}(\imath_{z})}_-{\cong}
      \ar@/^0.5pc/[dr]|-{t^{\mathcal{F}}_{w,b\tensor z}}
      &\text{ }
      \\
      \mathcal{F}(w\otimes z\otimes b)
      \ar[d]^-{\mathcal{F}(\I_w\otimes \gamma_{z})}
      &w\otimes \mathcal{F}(z)
      \ar[d]^-{t^{\mathcal{F}}_{w,z}}
      &\mathcal{F}(w\otimes b\tensor z)
      \ar@/^0.5pc/[dl]^-{\mathcal{F}(\I_w\otimes \imath_{z})}_-{\cong}
      \ar[dd]^-{\mathcal{F}(\gamma_{w}\tensor\I_{z})}
      \\
      \mathcal{F}(w\otimes z)
      \ar[d]^-{\mathcal{F}(\textit{cq}_{w,z})}
      \ar@{=}[r]
      &\mathcal{F}(w\otimes z)
      \ar[d]^-{\mathcal{F}(\textit{cq}_{w,z})}
      &\text{ }
      \\
      \mathcal{F}(w\tensor z)
      \ar@{=}[r]
      &\mathcal{F}(w\tensor z)
      \ar@{=}[r]
      &\mathcal{F}(w\tensor z)
    }
  }}
  \cdots
\end{equation*}
We used the right relation in (\ref{eq DayConvolution alphawbbb})
and the left relation in (\ref{eq DayConvolution alphawbCnat}).
We finish as follows.
\begin{equation*}
  \cdots
  \vcenter{\hbox{
    \xymatrix@R=16pt@C=20pt{
      w\otimes x\otimes z\otimes y
      \ar@{->>}[d]^-{\I_{w\otimes x}\otimes \textit{cq}_{z,y}}
      \ar@{=}[r]
      &w\otimes x\otimes z\otimes y
      \ar@{->>}[d]^-{\I_{w\otimes x}\otimes \textit{cq}_{z,y}}
      \ar@{=}[r]
      &w\otimes x\otimes z\otimes y
      \ar@{->>}[dd]^-{\textit{cq}_{w,x}\otimes \textit{cq}_{z,y}}
      \\
      w\otimes x\otimes (z\tensor y)
      \ar@{->>}[d]^-{\textit{cq}'_{w\otimes x,z\tensor y}}
      \ar@{=}[r]
      &w\otimes x\otimes (z\tensor y)
      \ar@{->>}[d]^-{\textit{cq}'_{w\otimes x,z\tensor y}}
      \ar@/^0.5pc/@{->>}[dr]|-{\textit{cq}_{w,x}\otimes \I_{z\tensor y}}
      &\text{ }
      \\
      w\otimes x\tensor' (z\tensor y)
      \ar[d]^-{(\I_w\otimes \imath^{-1}_x)\tensor' \I_{z\tensor y}}_-{\cong}
      \ar@{=}[r]
      &w\otimes x\tensor' (z\tensor y)
      \ar@{->>}[dd]^-{\textit{cq}_{w,x}\tensor'\I_{z\tensor y}}
      &(w\tensor x)\otimes (z\tensor y)
      \ar@{->>}[dd]^-{\textit{cq}'_{w\tensor x,z\tensor y}}
      \\
      (w\otimes b\tensor x)\tensor' (z\tensor y)
      \ar[d]^-{\alpha_{w\otimes b,z}}
      \ar@/^0.5pc/[dr]|-{(\gamma_{w}\tensor\I_x)\tensor'\I_{z\tensor y}}
      &\text{ }
      &\text{ }
      \\
      \mathcal{F}(w\otimes b\tensor z)
      \ar[d]^-{\mathcal{F}(\gamma_{w}\tensor\I_{z})}
      &(w\tensor x)\tensor' (z\tensor y)
      \ar[d]^-{\alpha_{w,z}}
      \ar@{=}[r]
      &(w\tensor x)\tensor' (z\tensor y)
      \ar[d]^-{\alpha_{w,z}}
      \\
      \mathcal{F}(w\tensor z)
      \ar@{=}[r]
      &\mathcal{F}(w\tensor z)
      \ar@{=}[r]
      &\mathcal{F}(w\tensor z)
    }
  }}
\end{equation*}
We leave for the readers to check the relation
$(\gamma_{w}\tensor\I_x)\circ a'_{w,b,x}
\circ (\I_w\otimes \imath^{-1}_x)
=\textit{cq}_{w,x}:w\otimes x
\to w\tensor x$
in $\textit{Mod}_{b'}$.
After right-cancelling the epimorphisms
$\textit{cq}_{w,x}\otimes \textit{cq}_{z,y}$
and
$\textit{cq}'_{w\tensor x,z\tensor y}$
in our calculation,
we obtain that
$\alpha_{w,z}=\xi_{w\tensor z}\circ \theta_{w,z}$.
This shows the existence of 
the enriched natural transformation
$\xi:-\tensor(x\wtensor y)\Rightarrow \mathcal{F}$.

Now we show the uniqueness of $\xi$.
The component $\theta_{b,b}$ of $\theta$ at $b$, $b$
satisfies the relation
\begin{equation*}
  \xymatrix@R=12pt@C=20pt{
    (b\tensor x)\tensor' (b\tensor y)
    \ar[rr]^-{\imath_x\tensor'\imath_y}_-{\cong}
    \ar[d]_-{\theta_{b,b}}
    &\text{ }
    &x\tensor' y
    \ar[d]^-{(p'_*)^2_{x,y}}
    \\
    b\tensor b\tensor (x\wtensor y)
    \ar[r]^-{\I_b\tensor\imath_{x\widetilde{\tensor} y}}_-{\cong}
    &b\tensor (x\wtensor y)
    \ar[r]^-{\imath_{x\widetilde{\tensor} y}}_-{\cong}
    &x\wtensor y
  }
\end{equation*}
which we obtain by right-cancelling the epimorphism
$\textit{cq}'_{b\tensor x,b\tensor y}$
in the diagram below.
\begin{equation*}
  \xymatrix@R=18pt@C=10pt{
    (b\tensor x)\otimes (b\tensor y)
    \ar@{->>}[d]_-{\textit{cq}'_{b\tensor x,b\tensor y}}
    \ar@{=}[r]
    &(b\tensor x)\otimes (b\tensor y)
    \ar@{->>}[d]^-{\textit{cq}_{b\tensor x,b\tensor y}}
    \ar@{=}[r]
    &(b\tensor x)\otimes (b\tensor y)
    \ar[dd]^-{\imath_{x}\otimes \imath_{y}}_-{\cong}
    \ar@{=}[r]
    &(b\tensor x)\otimes (b\tensor y)
    \ar@{->>}[d]^-{\textit{cq}'_{b\tensor x,b\tensor y}}
    \\
    (b\tensor x)\tensor' (b\tensor y)
    \ar[dd]_-{\theta_{b,b}}
    &b\tensor x\tensor b\tensor y
    \ar[d]^-{\I_b\tensor s_{x,y}\tensor\I_y}_-{\cong}
    &\text{ }
    &(b\tensor x)\tensor' (b\tensor y)
    \ar[dd]^-{\imath_x\tensor' \imath_y}_-{\cong}
    \\
    \text{ }
    &b\tensor b\tensor x\tensor y
    \ar[d]^-{\I_b\tensor\imath_{x\tensor y}}_-{\cong}
    \ar@/^0.5pc/[dl]|-{\I_{b\tensor b}\tensor(p_*)^2_{x,y}}
    &x\otimes y
    \ar@{->>}@/^0.5pc/[ddl]|-{\textit{cq}_{x,y}}
    \ar@{->>}[ddd]^-{\widetilde{\textit{cq}}_{x,y}}
    \ar@{->>}@/_0.5pc/[dr]|-{\textit{cq}'_{x,y}}
    &\text{ }
    \\
    b\tensor b\tensor (x\wtensor y)
    \ar[d]_-{\I_b\tensor\imath_{x\widetilde{\tensor} y}}^-{\cong}
    &b\tensor x\tensor y
    \ar[d]^-{\imath_{x\tensor y}}_-{\cong}
    \ar@/^0.5pc/[dl]|-{\I_b\tensor(p_*)^2_{x,y}}
    &\text{ }
    &x\tensor' y
    \ar@{->>}[dd]^-{(p'_*)^2_{x,y}}
    \\
    b\tensor (x\wtensor y)
    \ar[d]_-{\imath_{x\widetilde{\tensor} y}}^-{\cong}
    &x\tensor y
    \ar[d]^-{(p_*)^2_{x,y}}
    &\text{ }
    &\text{ }
    \\
    x\wtensor y
    \ar@{=}[r]
    &x\wtensor y
    \ar@{=}[r]
    &x\wtensor y
    \ar@{=}[r]
    &x\wtensor y
  }
\end{equation*}
The morphism
$(p'_*)_{x,y}$
is an epimorphism,
which follows from the defining relation (\ref{eq DayConvolution p*p'* formula}).
Thus $\theta_{b,b}$ is also an epimorphism.
Suppose we are given another enriched natural transformation
$\zeta:-\tensor(x\wtensor y)\Rightarrow \mathcal{F}$
which satisfies the relation
$\alpha_{w,z}=\zeta_{w\tensor z}\circ \theta_{w,z}$
for all $w$, $z\in\textit{Mod}_b$.
By right-cancelling the epimorphism $\theta_{b,b}$ in the relation
$\xi_{b\tensor b}\circ \theta_{b,b}
=\alpha_{b,b}
=\zeta_{b\tensor b}\circ \theta_{b,b}$
we obtain that
$\xi_{b\tensor b}
=\zeta_{b\tensor b}$.
Using the isomorphism
$\imath_{b}:b\tensor b\xrightarrow[]{\cong}b$
in $\textit{Mod}_b$,
we further obtain that $\xi_{b}=\zeta_{b}$.
Then we have $\xi_{z}=\zeta_{z}$
for all $z\in\textit{Mod}_b$
as we can see from the diagram below.
\begin{equation*}
  \xymatrix@R=12pt@C=30pt{
    z\tensor (x\wtensor y)
    \ar[ddd]_-{\xi_{z}}
    \ar@{=}[r]
    &z\tensor (x\wtensor y)
    \ar[d]^-{\bigl(\lambda^{-\tensor(x\widetilde{\tensor} y)}_{z}\bigr)^{-1}}_-{\cong}
    \ar@{=}[r]
    &z\tensor (x\wtensor y)
    \ar[ddd]^-{\zeta_{z}}
    \\
    \text{ }
    &z\tensor (b\tensor (x\wtensor y))
    \ar[d]^-{\I_{z}\tensor\xi_{b}\text{ }=\text{ }\I_{z}\tensor\zeta_{b}}
    &\text{ }
    \\
    \text{ }
    &z\tensor \mathcal{F}(b)
    \ar[d]^-{\lambda^{\mathcal{F}}_{z}}
    &\text{ }
    \\
    \mathcal{F}(z)
    \ar@{=}[r]
    &\mathcal{F}(z)
    \ar@{=}[r]
    &\mathcal{F}(z)
  }
\end{equation*}
We used the relation (\ref{eq EnEWProof lambda naturalinF})
for $\xi$ and $\zeta$.
Note that
$\lambda^{-\tensor(x\widetilde{\tensor} y)}$
is an enriched natural isomorphism,
since the enriched functor
$-\tensor (x\wtensor y):\textit{Mod}_b\to \textit{Mod}_{b'}$
is cocontinuous: see Proposition~\ref{prop EnEWThm defining lambda}.
This shows that $\xi=\zeta$ and that such $\xi$ is unique.
This completes the proof of the universal property of $\theta$.



\end{document}